\documentclass[a4paper, 12pt, reqno]{article}
\usepackage[cp866]{inputenc}
\usepackage[english]{babel}

\usepackage{amssymb, amsmath, amsthm, enumerate}
\usepackage{verbatim}
\usepackage{indentfirst}

\usepackage{amsfonts,amssymb}
 \newtheorem{thm}{Theorem}[subsection]
 \newtheorem{cor}[thm]{Corollary}
 \newtheorem{lem}[thm]{\bf Lemma}
 \newtheorem{prop}[thm]{\bf Proposition}
 \theoremstyle{definition}
 \newtheorem{defn}[thm]{Definition}
 \theoremstyle{remark}
 \newtheorem{rem}[thm]{\bf Remark}
  \newtheorem{assump}[thm]{\bf Assumption}
 \newtheorem{notat}[thm]{\bf Notation}
  \newtheorem{property}[thm]{\bf Property}
    \newtheorem*{ex}{\bf Example}
 \numberwithin{equation}{subsection}

\begin{document}

\begin{center}{\sc\large Geometry of Carnot--Carath\'{e}odory Spaces,\\
 Differentiability and
Coarea Formula}
\end{center}

\begin{center}
{\bf Maria Karmanova, Sergey Vodopyanov}
\end{center}

\begin{abstract}
We give a simple proof of  Gromov's Theorem  on nilpotentization of vector fields,
and exhibit a new method for  obtaining  quantitative estimates of comparing geometries of two different local Carnot groups in Carnot--Carath\'{e}odory
spaces with $C^{1,\alpha}$-smooth basis vector fields, $\alpha\in[0,1]$.
From here we obtain the similar estimates for comparing geometries of a
Carnot--Carath\'{e}odory space and a local Carnot group. These
two theorems imply basic results of the theory:
Gromov type Local Approximation Theorems,
and for  $\alpha>0$ Rashevski\v{\i}-Chow
Theorem and Ball--Box Theorem, etc.
We apply the obtained results for proving $hc$-differentiability
of mappings of Carnot--Carath\'{e}odory spaces with continuous
horizontal derivatives. The latter is used in proving the coarea
formula for some classes of contact mappings of Carnot--Carath\'{e}odory spaces.
\end{abstract}

\tableofcontents

\section{Introduction}

The geometry of Carnot--Carath\'eodory spaces naturally arises in
the theory of subelliptic equations, contact geometry, optimal
control theory, nonholonomic mechanics, neurobiology and other
areas (see {works by A.~A.~Agra\-chev~{\cite{agr}},
A.~A.~Ag\-ra\-chev and J.-P.~Gauthier~{\cite{agrg}},
A.~A.~Agrachev and A.~Mari\-go~\cite{AM}, A.~A.~Agrachev and
A.~V.~Sarychev~\cite{agrs1, agrs2, agrs3, agrs4, agrs5},
A.~Bella\"{\i}che~\cite{B}, A. Bonfiglioli, E. Lanconelli and F.
Uguzzoni~{\cite{blu}}, S.~Buckley, P.~Koskela and
G.~Lu~{\cite{buckleykl1}}, L.~Capogna{\cite{cap1, cap2}},
G.~Citti, N.~Garofalo and E.~Lanconelli~{\cite{cgl}}, L.~Capogna,
D.~Danielli and N.~Garofalo~{\cite{capognadg2, capognadg2.5,
capognadg4, capognadg5, capognadg6}}, Ya.~Eliashberg~{\cite{el1,
el2, el3, el4}}, G.~B.~Folland~{\cite{Foll, Foll1}}, G.~B.~Folland
and E.~M.~Stein~{\cite{fs}}, B.~Franchi, R.~Serapioni, F.~Serra
Cassano~{\cite{fssc1, fssc2, fssc3,  fssc4}},
N.~Garofalo~{\cite{garofalobook}}, N.~Garofalo and
D.-M.~Nhieu~{\cite{GN1,GN2}}, R.~W.~Goodman~{\cite{good}},
M.~Gromov~\cite{G, G2001}, L.~H\"ormander~\cite{hor},
F.~Jean~\cite{J}, V.~Jurdjevic~{\cite{jur}}, G.~P.~Leonardi,
S.~Rigot~{\cite{LeoRig}}, W.~Liu and H.~J.~Sussman~{\cite{LS}},
G.~Lu~{\cite{Lu}}, G.~A.~Margulis and G.~D.~Mostow~\cite{MM,MM1},
G.~Me\-ti\-vier~\cite{Me}, J.~Mitchell~{\cite{Mit}},
R.~Montgomery~\cite{Montgom1, M}, R.~Monti~{\cite{Mon, Mon1}},
A.~Nagel, F.~Ricci, E.~M.~Stein~{\cite{NRS, NRS1}}, A.~Nagel,
E.~M.~Stein and S.~Wain\-ger~{\cite{nsw}}, P.~Pansu~{\cite{P,
Pan01, Pan1, Pan}}, L.~P.~Rothschild and E.~M.~Stein~\cite{rs}, R.
S. Strichartz~{\cite{Stri}}, A.~M.~Vershik and
V.~Ya.~Gershkovich~\cite{VerG}, S.~K.~Vodopyanov~{\cite{V4, vod1,
vod2, vod3, V3}}, S.~K.~Vodopyanov and A.~V.~Greshnov~{\cite{vg}},
C.~J.~Xu and C.~Zuily~{\cite{XZ}} for an introduction to this
theory and some its applications).

A Carnot--Carath\'eodory space  (below referred to as a  {\it
Carnot manifolds}) $\mathbb{M}$ (see, for example, \cite{G,VerG})
is a connected Riemannian manifold with a distinguished horizontal
subbundle $H\mathbb{M}$ in the tangent bundle $T\mathbb{M}$ that
meets some algebraic conditions on the commutators of vector
fields $\{X_1,\dots,X_n\}$ constituting a local basis in
$H\mathbb{M}$, $n=\dim H\mathbb{M}$.

The distance $d_c$ (the intrinsic Carnot--Carath\'eodory metric)
between points $x,y\in\mathbb{M}$ is defined as the infimum of the
lengths of {\it horizontal} curves joining $x$ and $y$ and is
non-Riemannian if $H\mathbb{M}$ is a proper subbundle (a piecewise
smooth curve $\gamma$ is called horizontal if
 $\dot\gamma(t)\in H_{\gamma(t)}\mathbb{M}$). See results on properties of this metric in
 the monograph by D.~Yu.~Burago, Yu.~D.~Burago nd S.~V.~Ivanov~{\cite{bbi}}.

The Carnot--Carath\'eodory metric is applied in the study of
hypoelliptic operators, see C.~Fefferman and D.~H.~Phong
\cite{feffermanp}, L.~H\"ormander \cite{hor}, D.~Jerison
\cite{jerison}, A.~Nagel, E.~M.~Stein and S.~Wainger \cite{nsw},
L.~P.~Rothschild and E.~M.~Stein \cite{rs}, A.~S\'anchez-Calle
\cite{sanchez}. Also, this metric and its properties are
essentially used in theory of PDE's (see papers by M.~Biroli and
U.~Mosco~\cite{birolim1, birolim4}, S.~M.~Buckley, P.~Koskela and
G.~Lu~\cite{buckleykl1}, L.~Capogna, D.~Danielli and
N.~Garofalo~\cite{capognadg2, capognadg2.5, capognadg4,
capognadg5, capognadg6}, V.~M.~Chernikov and
S.~K.~Vodop'yanov~\cite{chernikovv}, D.~Danielli, N.~Garofalo and
D.-M.~Nhieu~\cite{daniellign}, B.~Franchi~\cite{franchi1},
B.~Franchi, S.~Gallot and R.~Wheeden~\cite{franchigaw},
B.~Franchi, C.~E.~Guti\'errez and R.~L.~Wheeden~\cite{franchigw},
B.~Franchi and E.~Lanconelli~\cite{franchil1,franchil2},
B.~Franchi, G.~Lu and R.~Wheeden~\cite{franchilw1, franchilw3},
B.~Franchi and R.~Serapioni~\cite{franchis1},
R.~Garattini~\cite{garattini}, N.~Garofalo and
E.~Lanconelli~\cite{garofalol}, P.~Haj\l{}asz and
P.~Strzelecki~\cite{hajlaszs}, J.~Jost~\cite{jost1, jost2,jost3,
jost4}, J.~Jost and C.~J.~Xu~\cite{jostx},
S.~Marchi~\cite{marchi}, K.~T.~Sturm~\cite{sturm4}).

The following results are usually regarded as foundations of the
geometry of Carnot manifolds:

\begin{enumerate}
\item Rashevski\v{\i}--Chow  Theorem {\cite{chow,rash}}  on
connection of two points by a horizontal path;

\item Ball--Box Theorem {\cite{nsw}} (saying that a ball in
 Carnot--Ca\-ra\-th\'{e}\-o\-do\-ry metric contains a ``box'' and is a subset of a
``box'' with controlled ``radii'');

\item Mitchell's Theorem {\cite{Mit}} on convergence of rescaled
Carnot--Carath\'eo\-dory spaces around $g\in\mathbb M$ to a
nilpotent tangent cone;

\item\label{convgr} Gromov's Theorem {\cite{G}} on convergence of
``rescaled'' with respect to $g\in\mathbb M$ basis vector fields
to {\it nilpotentized} ({\it at} $g$) vector fields generating a
graded nilpotent Lie algebra (the corresponding connected and
simply connected Lie group is called the {\it nilpotent tangent
cone at $g$}); here $g\in\mathbb M$ is an arbitrary point;

\item Gromov Approximation Theorem {\cite{G}} on local comparison
of Carnot--Ca\-ra\-th\'{e}\-o\-do\-ry metrics in the initial space
and in the nilpotent tangent cone, and its improvements due to A.
Bella\"{\i}che {\cite{B}}.
\end{enumerate}

The goal of the paper is both to give a new approach to the
geometry of Carnot manifolds and to establish some basic results
of geometric measure theory on these metric structures including
an appropriate differentiability and a coarea formula.

New results in the geometry of Carnot manifolds  contains
essentially new quantitative estimates of closeness of geometries
of different tangent cones located one near another. One of the
peculiarities of the paper is that we solve all problems under
minimal assumption on smoothness of the basis vector fields (they
are $C^{1,\alpha}$-smooth, $0\leq\alpha\leq 1$), although all the
basic results are new even for $C^{\infty}$-vector fields. In some
parts of this paper, the symbol $C^{1,\alpha}$ means that the
derivatives of the basis vector fields are $H^{\alpha}$-continuous
with respect to some nonnegative symmetric function $\mathfrak
d:U\times U\to\mathbb R$, $U\Subset\mathbb M$, such that
$\mathfrak d\geq C\rho$, $0<C<\infty$, where $C$ depends only on
$U$, and $\rho$ is Riemannian distance. Some additional properties
of $\mathfrak d$ are described below when it is necessary. Note
that from the very beginning it is unknown whether
Rashevski\v{\i}--Chow Theorem is true for $C^{1,\alpha}$-smooth
basis vector fields. Therefore Carnot--Carath\'eodory distance can
not be well defined. We use the quasimetric $d_{\infty}$ instead
of $d_c$, which  is defined as follows: if
$y=\exp\Bigl(\sum\limits_{i=1}^Ny_iX_i\Bigr)(x)$, then
$d_{\infty}(x,y)=\max\limits_{i=1,\ldots,
N}\{|y_i|^{\frac{1}{\operatorname{deg}X_i}}\}$, and in smooth case
is equivalent to $d_c$ \cite{nsw, G}. One of the main results is
the following (see below Theorem {\ref{est_eps}} for sharp
statement).

\begin{thm}
Suppose that $d_{\infty}(u,u^{\prime})=C\varepsilon$,
$d_{\infty}(u, v)={\mathcal C}\varepsilon$ for some $C,{\mathcal
C}<\infty$,
$$
w_{\varepsilon}=\exp\Bigl(\sum\limits_{i=1}^Nw_i\varepsilon^{\operatorname{deg}
X_i}\widehat{X}^{u}_i\Bigr)(v) \text{ and }
w_{\varepsilon}^{\prime}=\exp\Bigl(\sum\limits_{i=1}^Nw_i\varepsilon^{\operatorname{deg}
X_i}\widehat{X}^{u^{\prime}}_i\Bigr)(v).
$$
Then, for $\alpha>0$, we have
$$
\max\{d_{\infty}^u(w_{\varepsilon},w^{\prime}_{\varepsilon}),
d_{\infty}^{u^{\prime}}(w_{\varepsilon},w^{\prime}_{\varepsilon})\}\leq
L\varepsilon\rho(u,u^{\prime})^{\frac{\alpha}{M}}
$$
where $L$ is uniformly bounded in $u, u^{\prime}, v\in
U\subset\mathbb M$, and in $\{w_i\}_{i=1}^N$ belonging to some
compact neighborhood of $0$.

In the case of $\alpha=0$, we have
$$
\max\{d_{\infty}^u(w_{\varepsilon},w^{\prime}_{\varepsilon}),
d_{\infty}^{u^{\prime}}(w_{\varepsilon},w^{\prime}_{\varepsilon})\}\leq\varepsilon
o(1)
$$
where $o(1)$ is uniform in $u, u^{\prime}, v\in U\subset\mathbb
M$, and in $\{w_i\}_{i=1}^N$ belonging to some compact
neighborhood of $0$ as $\varepsilon\to0$.
\end{thm}

Here we assume that $U\subset\mathbb M$ is a compact neighborhood
small enough and $\rho$ is a Riemannian metric. The symbol
$\widehat{X}_i^u$ ($\widehat{X}_i^{u^{\prime}}$) denotes
nilpotentized at $u$ ($u^{\prime}$) vector fields (see item
{\ref{convgr}} above). These vector fields constitute Lie algebra
of the nilpotent tangent cone at $u$ ($u^{\prime}$).

Further, in Theorem {\ref{est_chain}} we extend this result to the
case of a ''chain`` consisting of several points.

The content of obtained estimates is very profound: they imply
both new properties of Carnot manifolds and the above-mentioned
ones.

The investigation of sub-Riemannian geometry under minimal
smoothness of the basis vector fields is motivated by the recently
constructed by G.~Citti and A.~Sarti, and R.~K.~Hladky and
S.~D.~Pauls model of visualization~{\cite{CS, HP}}. More exactly,
the model of a brain perception of a black-and-white plain image
is constructed in these papers. This model makes possible the
interpretation on a computer of a human brain's work during the
visualization of information. In particular, it is shown how the
human brain completes the image part of which is closed. The
geometry of this model is based on a roto-translation group which
is a three-dimensional Carnot manifold with a tangent cone being a
Heisenberg group $\mathbb H^1$ at each point. Since by now there
are no theorems on regularity of the image created by a human
brain, any reduction of smoothness of vector fields is essential
for the construction of sharp visualization models.

The main result  concerning the geometry of Carnot manifolds is
proved in Section~{\ref{geomcarnot}}. The method of proving is
new, and it essentially uses H\"{o}lder dependence of solutions to
ordinary differential equations on parameter (see Theorem
{\ref{ODE}}). Probably, this dependence is not a new result. For
reader's convenience
 we give its independent proof in Section~{\ref{proofdiff}}. In Subsection~{\ref{prelim}}, all other auxiliary result are formulated.

In Subsection~{\ref{GrTh}} we prove, in particular, the following
statements
\smallskip

{\bf A:}  \textit{Let  $X_j\in C^{1}$ on a Carnot manifold
$\mathbb M$. On  $\operatorname{Box}(g,\varepsilon r_g)$, consider
the vector fields $\{^{\varepsilon}X_i\}{=}\{\varepsilon^{\deg
X_i}X_i\}$, $i=1,\dots,N$. Then the uniform convergence
$$
X_i^{\varepsilon}=\bigl(\Delta^g_{\varepsilon^{-1}}\bigr)_*{^{\varepsilon}X_i}
\to \widehat{X}^g_i\quad \text{as $\varepsilon\to 0$},\quad
i=1,\dots,N,
$$
holds at the points of the box~$\operatorname{Box}(g,r_g)$ and
this convergence is uniform in  $g$ belonging to some compact set,
where the collection $\{\widehat{X}^g_i\}$, $i=1,\ldots,N$, of
vector fields around $g$ constitutes a basis of  a graded
nilpotent Lie algebra};
\smallskip

{\bf B:}  \textit{There exists a constant $Q=Q(U)$, $U$ is a
compact domain in $\mathbb M$, such that the inequality
$$
d_\infty(u,v)\leq Q(d_\infty(u,w)+d_\infty(w,v))
$$
holds for every triple of points $u$, $w$, $v\in U$ where $Q(U)$
depends on $U$}.
\smallskip

{\bf C:}  \textit{Given points $u,v\in \mathbb{M}$,
 $d_{\infty}(u,
v)={\mathcal C}\varepsilon$ for some ${\mathcal C}<\infty$,
$$
w_{\varepsilon}=\exp\Bigl(\sum\limits_{i=1}^Nw_i\varepsilon^{\operatorname{deg}
X_i}{X}_i\Bigr)(v) \text{ and }
w_{\varepsilon}^{\prime}=\exp\Bigl(\sum\limits_{i=1}^Nw_i\varepsilon^{\operatorname{deg}
X_i}\widehat{X}^{u}_i\Bigr)(v),
$$
we have
$$
\max\bigl\{d_{\infty}(w_{\varepsilon},w^{\prime}_{\varepsilon}),
d_{\infty}^{u}(w_{\varepsilon},w^{\prime}_{\varepsilon})\bigr\}\leq\varepsilon
o(1)
$$
where $o(1)$ is uniform in $u, v$ belonging to a compact
neighborhood $U\subset\mathbb M$, and in $\{w_i\}_{i=1}^N$
belonging to some compact neighborhood of $0$ as
$\varepsilon\to0$}.
\smallskip

Statement {\bf A} is just Gromov's Theorem \cite{G}  on the
nilpotentization of vector fields. Gromov has formulated it for
$C^1$-smooth fields, however, Example \ref{cex} by  Valeri\v{\i}
Berestovski\v{\i}
 makes evident that arguments of the proof given in {\cite[pp.~128--133]{G}} have to be corrected. In  Corollary \ref{GTheorem}
we give a new proof of this assertion based on an another idea.

Statement {\bf B} says that the quasimetric  $d_{\infty}$ meets
the generalised triangle inequality. The implication $\text{\bf
A}\Longrightarrow\text{\bf B }$ is proved in Corollary
\ref{trin}.

Statement {\bf C} gives an estimate of divergence of integral
lines of the given vector fields and the nilpotentized vector
fieldes.

 The implication $\text{\bf B}\Longrightarrow\text{\bf C}$ is a particular case of
 Theorem \ref{mainresult}.

In theory developed in Subsection~{\ref{geomconecone1}}, is based
on the generalized triangle inequality as a starting point.

 In Subsection~{\ref{geomconecone1}}, we prove one of the basic results of Section~{\ref{geomcarnot}}, namely, Theorem~{\ref{est_eps}} which compares local geometries of two different local Carnot groups. It is essentially based on the main theorem of Subsection~{\ref{geomconeconeglob}} which compares ''global`` geometries of different tangent cones (i.~e., it looks like Theorem~{\ref{est_eps}} with $\varepsilon=1$). Subsection~{\ref{apprths}} is devoted to approximation theorems.
In particular, we compare  metrics of two tangent cones, and the
metric of a tangent cone with the initial one. There we give their
proofs and the proofs of some auxiliary properties of the
geometry. Further, in Subsection~{\ref{geomconecone}}, we prove
Theorem~{\ref{est_chain}}, which is the ''continuation`` of
Theorem~{\ref{est_eps}}. In  Subsection~{\ref{geomcarnotcone}}, we
compare the geometry of a Carnot manifold with the one of a
tangent cone. In Subsection~{\ref{applications}}, we give
applications of our results to investigation of the sub-Riemannian
geometry. We prove Gromov type theorem on the nilpotentization of
vector fields \cite{G},   a new statement implying
Rashevski\v{\i}--Chow Theorem, Ball--Box Theorem, Mitchell Theorem
on Hausdorff dimension of Carnot manifolds and many other
corollaries.

Main results  of Section~{\ref{geomcarnot}}  are formulated in
short communications {\cite{vk1, vk2}}.

Section~{\ref{differentiabilitycarnot}} is devoted to
differentiability of mappings in the category of Carnot manifolds.

 We recall the classical definition of differentiability for a
mapping $f:\mathbb M\to \mathbb N$ of two Riemannian manifolds:
$f$  is differentiable at $x\in \mathbb M$ if there exists a
linear mapping $L:T_x\mathbb M\to T_{f(x)}\mathbb N$ of the
tangent spaces such that
$$
\rho_{\mathbb N}(f(\exp_x h), \exp_{f(x)}Lh)=o(\|h\|_x),\quad h\in
T_x\mathbb M,
$$
where $\exp_x:T_x\mathbb M\to\mathbb M$ and
$\exp_{f(x)}:T_{f(x)}\mathbb N\to\mathbb N$ are the exponential
mappings, and $\rho_{\mathbb N}$ is the Riemannian metric in
$\mathbb N$, $\|h\|_x$ is the length of $h\in T_x \mathbb M$.

It is known (see \cite{G,M}) that the local geometry of a Carnot
manifold at a point $g\in\mathbb{M}$ can be modelled as a~graded
nilpotent Lie group ${\mathbb G}_g\mathbb M$. It means that the
tangent space $T_g\mathbb M$ has an additional structure of
a~graded nilpotent Lie group. If~$\mathbb M$ and~$\mathbb N$ are
two Carnot manifolds and $f:\mathbb M\to\mathbb N$ is a mapping
then a suitable concept of differentiability can be obtained from
the previous concept in the following way: $f$  is
$hc$-differentiable at $x\in \mathbb M$ if there exists a
horizontal homomorphism   $L:\mathbb G_x\mathbb M\to \mathbb
G_{f(x)}\mathbb N$ of the nilpotent tangent cones such that
$$
\tilde d_c(f(\exp_xh),\exp_{f(x)}Lh)=o(|h|_x), \quad h\in\mathbb
G_x\mathbb M,
$$
where $\tilde d_c$ is the Carnot--Carath{\'e}odory distance in
$\mathbb N$ and $|\cdot|_x$ is an homogeneous norm in $\mathbb G_x
\mathbb M$.

For us, it is convenient to regard some neighborhood of a point
$g$ as a subspace of the metric space $(\mathbb{M},d_{c})$ and as
a neighborhood of unity of the local Carnot group  ${\mathcal
G}^g\mathbb{M}$ with Carnot--Carath\'eodory metric $d^{g}_c$ (see
Definition~1.2). In the sense explained below,
$\exp^{-1}:{\mathcal G}^g\mathbb{M}\to \mathbb G_g\mathbb M$ is an
isometrical  monomorphism of the Lie structures. Then the last
definition of $hc$-differentiability can be reformulated as
follows. Given two Carnot manifolds $(\mathbb M,d_{c})$ and
$(\mathbb N, \tilde d_{c})$ and a set $E\subset\mathbb{M}$,
a~mapping $f:E\to \mathbb N$ is called {\it $hc$-differentiable}
at a point $g\in E$ (see the paper by S.~K.~Vodopyanov and
A.~V.~Greshnov~\cite{vg}, and also \cite{vod1,vod2,vod3,V3}) if
there exists a horizontal homomorphism $L:\bigl({\mathcal
G}^g\mathbb M, d^g\bigr)\to \bigl({\mathcal G}^{f(g)}\mathbb N,
d_c^{f(g}\bigr)$ of the local Carnot groups such that
\begin{equation}\label{defdiffer}
d^{f(g)}_c(f(w),L(w))=o(d^g_c(g,w))\quad\text{as $E\cap {\mathcal
G^g}\ni w\to g$}.
\end{equation}
The given definition of $hc$-differentiability of mappings for
Carnot manifolds
 can be treated as a straightforward generalization of the classical definition
 of differentiability.
Clearly, if the Carnot manifolds are Carnot groups then this
definition of $hc$-differentiability is equivalent to the
definition of $\mathcal P$-differentiability introduced by
P.~Pansu in~\cite{Pan} for an open set $E\subset\mathbb G$. For an
arbitrary $E\subset\mathbb G$, the last concept was investigated
by S.~K.~Vodop$'$yanov~{\cite{V4}} and by S.~K.~Vodopyanov and
A.~D.~Ukhlov~\cite{VU1} (see also the paper by
V.~Magnani~\cite{M1}).

In Section~{\ref{differentiabilitycarnot}}, we introduce  the
notion of $hc$-differentiability, which is adequate to the
geometry of Carnot manifold, and study its properties. Moreover,
in this section, we prove the $hc$-differentiability of a
composition of $hc$-differentiable mappings.

In the same section we prove the $hc$-differentiability of
rectifiable curves. In the case of curves, the definition of the
$hc$-differentiability is interpreted as follows: a mapping
$\mathbb R\supset E\ni t\mapsto \gamma(t)\in\mathbb{N}$ is
$hc$-differentiable at a point $s\in E$ in a Carnot manifold
$\mathbb{N}$ if the relation
\begin{equation}\label{defdiffercurves}
d^{\gamma(s)}_c\Bigl(\gamma(s+\tau), \exp(\tau a)(\gamma(s))\Bigr)
=o(\tau)\quad\text{as $\tau\to0$, $s+\tau\in E$,}
\end{equation}
holds, where $a\in H_{\gamma(s)}\mathbb N$
(\eqref{defdiffercurves} agrees with \eqref{defdiffer} when
$\mathbb M=\mathbb R$, see also {\cite{MM}}).
 On Carnot groups, relation~\eqref{defdiffercurves} is
equivalent to the $\mathcal P$-differentiability of curves in the
sense of Pansu~\cite{Pan}. Our proof of differentiability is new
even for Carnot groups.
 We prove step by step the $hc$-differentiability
 of the absolutely differentiable curves, the Lipschitz mappings of
 subsets of $\mathbb R$ into $\mathbb M$, and the rectifiable curves.
 Here we generalize a classical result and obtain
 the following assertion: {\it
the continuity of horizontal derivatives of a contact mapping
defined on an open set implies its pointwise
$hc$-differentiability} (Theorem~\ref{theorem3.3}).

As an important corollary to these assertion, we infer that the
nilpotent tangent cone is defined by the horizontal subbundle of
the Carnot manifold: {\it tangent cones found from different
collections of basis vector fields are isomorphic as local Carnot
groups} (Corollary~\ref{cor3.5}). Thus, the correspondence ``local
basis $\mapsto$ nilpotent tangent cone'' is functorial. In the
case of $C^\infty$-vector fields, this result was established by
A.~Agrachev and A.~Marigo~\cite{AM}, and G.~A.~Margulis and
G.~D.~Mostow~{\cite{MM1}} where a coordinate-free definition of
the tangent cone was given.

Main results  of Section~{\ref{differentiabilitycarnot}}  are
formulated in short communications by
S.~K.~Vo\-do\-pyanov~{\cite{vod1, vod2}} (see some details and
more general results on this subject including
Rademacher--Stepanov Theorem in \cite{vod3,V3}).

Section {\ref{coareacarnot}} is dedicated to such application of
results on $hc$-differentiability as the sub-Riemannian analog of
the coarea formula. It is well known that the coarea formula
\begin{equation}\label{coarea_e}
\int\limits_{U}{\mathcal J}_k(\varphi,x)\,dx=\int\limits_{\mathbb
R^k}\,dz\int\limits_{\varphi^{-1}(z)}\,d{\mathcal H}^{n-k}(u),
\end{equation}
where ${\mathcal
J}_k(\varphi,x)=\sqrt{\det(D\varphi(x)D\varphi^*(x))}$, has many
applications in analysis on Euclidean spaces. Here we assume that
$\varphi\in C^1(U,\mathbb R^k)$, $U\subset\mathbb R^n$, $n\geq k$.
For the first time, it was established by
A.~S.~Kronrod~{\cite{kr}} for the case of a function
$\varphi:\mathbb R^2\to\mathbb R$. Next, it was generalized by
H.~Federer first for mappings of Riemannian manifolds
$\varphi:\mathcal M^n\to\mathcal N^k$, $n\geq k$, in {\cite{fd1}},
and then,  for mappings of rectifiable sets in Euclidean spaces
$\varphi:\mathcal M^n\to\mathcal N^k$, $n\geq k$, in {\cite{F}}.
Next, in the paper~{\cite{ot}}, M.~Ohtsuka generalized the coarea
formula~{\eqref{coarea_e}} for mappings $\varphi:\mathbb
R^n\to\mathbb R^m$, $n,m\geq k$, with $\mathcal
H^k$-$\sigma$-finite image $\varphi(\mathbb R^n)$. An
infinite-dimensional analog of the coarea formula was proved by
H.~Airault and P.~Malliavin in 1988 {\cite{AMal}} for the case of
Wiener spaces. This result can be found in the monograph by
P.~Malliavin {\cite{Mal}}. See other proofs and applications of
the coarea formula in the monographs by L.~C.~Evans and
R.~F.~Gariepy~{\cite{eg}}, M.~Giaquinta, G.~Modica and
J.~Sou\v{c}ek~{\cite{gq}}, F.~Lin and X.~Yang~ {\cite{fx}}.

Formula {\eqref{coarea_e}} can be applied in the theory of
exterior forms, currents, in minimal surfaces problems (see, for
example, paper by H.~Federer and W.~H.~Fleming~{\cite{FF}}). Also,
Stokes formula can be easily obtained by using the coarea formula
(see, for instance, lecture notes by
S.~K.~Vodopyanov~{\cite{V1}}). Because of the development of
analysis on more general structures, a natural question arise to
extend the coarea formula on objects of more general geometry in
comparison with Euclidean spaces, especially on metric spaces and
structures on sub-Riemannian geometry. In 1999, L.~Ambrosio and
B.~Kirchheim {\cite{AK}} proved the analog of the coarea formula
for Lipschitz mappings defined on ${\mathcal H}^n$-rectifiable
metric space with values in $\mathbb R^k$, $n\geq k$. In 2004,
this formula was proved for Lipschitz mappings defined on
${\mathcal H}^n$-rectifiable metric space with values in
${\mathcal H}^k$-rectifiable metric space, $n\geq k$, by
M.~Karmanova~{\cite{Km1, Km3}}. Moreover, necessary and sufficient
conditions on the image and the preimage of a Lipschitz mapping
defined on $\mathcal H^n$-rectifiable metric space with values in
an {\it arbitrary} metric space for the validity of the coarea
formula were found. Independently of this result, the level sets
of such mappings were investigated, and the metric analog of
Implicit Function Theorem was proved~by M.~Karmanova~{\cite{Km2,
Km3, Km4}}.

All the above results are connected with rectifiable metric
spaces. Note that, their structure is similar to the one of
Riemannian manifolds. But there are also {\it non-rectifiable}
metric spaces which geometry is not comparable with the Riemannian
one. {\it Carnot manifolds} are of special interest. The problem
of the sub-Riemannian coarea formula is one of well-known
intrinsic unsolved problems.

A Heisenberg group and a Carnot group are well-known particular
cases of a Carnot manifold. In 1982, P.~Pansu proved the coarea
formula for functions defined on a Heisenberg group {\cite{P}}.
Next, in {\cite{H}}, J.~Heinonen extended this formula to smooth
functions defined on a Carnot group. In~{\cite{MS}}, R.~Monti and
F.~Serra Cassano proved the analog of the coarea formula for
$BV$-functions defined on a two-step Carnot--Carath\'{e}odory
space. One more result concerning the analogue of
{\eqref{coarea_e}} belongs to V.~Magnani. In 2000, he proved a
{\it coarea inequality} for mappings of Carnot groups
{\cite{M00}}. The equality was proved only for the case of a
mapping defined on a Heisenberg group with values in Euclidean
space~$\mathbb R^k$~{\cite{M05}}. Until now, the question about
the validity of coarea formula even for a model case of a mapping
of Carnot groups was open.

Main results  of Section~{\ref{coareacarnot}}  are formulated in
{\cite{vk3}}.

\section{Geometry of Carnot--Carath\'{e}odory Spaces}\label{geomcarnot}

\subsection{Preliminary Results}\label{prelim}

Recall the definition of a Carnot manifold.

\begin{defn}[compare with
{\cite{G,nsw}}]\label{carnotmanifold} Fix a connected Riemannian
$C^{\infty}$-mani\-fold~$\mathbb M$ of a dimension~$N$. The
manifold~$\mathbb M$ is called a {\it Carnot manifold} if, in the
tangent bundle $T\mathbb M$, there exists a tangent subbundle
$H\mathbb M$ with a finite collection of natural numbers $\dim
H_1<\ldots<\dim H_i<\ldots<\dim H_M=N$, $1<i<M$, and each point
$p\in\mathbb M$ possesses a neighborhood $U\subset\mathbb M$
 with a collection of $C^1$-smooth vector fields
$X_1,\dots,X_N$ on $U$ enjoying the following properties.

For each $v\in U$ we have

$(1)$ $X_1(v),\dots,X_N(v)$ constitutes a basis of $T_v\mathbb M$;

$(2)$ $H_i(v)=\operatorname{span}\{X_1(v),\dots,X_{\dim H_i}(v)\}$
is a subspace of $T_v\mathbb M$ of a dimension $\dim H_i$,
$i=1,\ldots,M$, where $H_1(v)=H_v\mathbb M$;

$(3)$
\begin{equation}\label{tcomm}[X_i,X_j](v)=\sum\limits_{\operatorname{deg}
X_k\leq \operatorname{deg} X_i+\operatorname{deg}
X_j}c_{ijk}(v)X_k(v)
\end{equation}
where the {\it degree} $\deg X_k$ equals $\min\{m\mid X_k\in
H_m\}$;

(4) a quotient mapping $[\,\cdot ,\cdot\, ]_0:H_1\times
H_j/H_{j-1}\mapsto H_{j+1}/H_{j}$ induced by Lie brackets, is an
epimorphism for all $1\leq j<M$.
\end{defn}

\begin{rem}
Note {\cite{nsw}} that  Condition~4 is necessary only for
obtaining results of Subsections~{\ref{applications}} and
~{\ref{diffsmmap}} and its corollaries, and
Section~\ref{coareacarnot}. The point is that in
 statements of these subsections, we use the fact that any two points of a local
Carnot group (see the definition below) can be joined by a
horizontal (with respect to the local Carnot group) curve that
consists of at most $L$ segments of integral lines of horizontal
(with respect to the local Carnot group) vector fields. The latter
is impossible without Condition~4.
\end{rem}

\begin{rem}\label{cijkcont} Consider a $C^2$-smooth local
diffeomorphism $\eta:U\to\mathbb R^N$, $U\subset\mathbb M$. Then
$\eta_*X_i=D\eta\langle X_i\rangle$ are also $C^1$-vector fields,
$i=1,\ldots, N$. We have the following relations instead of
{\eqref{tcomm}}:
$$
\eta_*[X_i,X_j](w)=[\eta_*X_i,\eta_*X_j](w)=\sum\limits_{\operatorname{deg}
X_k\leq \operatorname{deg} X_i+\operatorname{deg}
X_j}c_{ijk}(\eta^{-1}(w))\eta_*X_k(w).
$$
Denote by $X(w)$ the matrix, the $i$th column of which consists of
the coordinates of $\eta_*X_i(w)$ in the standard basis
$\{\partial_j\}_{j=1}^N$. Then the entries of $X(w)$ are
$C^1$-functions. Note that
$$
\eta_*[X_i,X_j](w)=X(w)(c_{ij1}(\eta^{-1}(w)),\ldots,c_{ijN}(\eta^{-1}(w)))^T.
$$
Consequently,
$$
(c_{ij1}(\eta^{-1}(w)),\ldots,c_{ijN}(\eta^{-1}(w)))^T=(X(w))^{-1}\cdot\eta_*[X_i,X_j](w).
$$
From here it follows that $c_{ijk}\circ\eta^{-1}$ are continuous,
$k=1,\ldots, N$. Since $\eta$ is continuous, then we have that
each $c_{ijk}=(c_{ijk}\circ\eta^{-1})\circ\eta$ is also
continuous, $k=1,\ldots,N$.
\end{rem}

\begin{ex} [A Carnot Manifold with $C^1$-Smooth Vector Fields] Consider the $C^M$ vector fields $X_1,\ldots,X_n\in H$.
Choose a basis in \linebreak$H_2=\operatorname{span}\{H,[H,H]\}$
by the following way:
$$
X_k(v)=\sum\limits_{i,j}a^k_{ij}(v)[X_i,X_j](v)+\sum\limits_lb^k_l(v)X_l,
$$
where $a^k_{ij}(v), b^k_l(v)\in C^1$, $i,j,l=1,\ldots, N$,
$k=n+1,\ldots,\dim H_2$. Similarly, we choose the following basis
in $H_{m+1}=\operatorname{span}\{H_{m},[H,H_{m}]\}$,
$m=2,\ldots,M-1$:
$$
X_k(v)=\sum\limits_{i,j}a^k_{ij}(v)[X_i,X_j](v)+\sum\limits_lb^k_l(v)X_l,
$$
where $a^k_{ij}(v), b^k_l(v)\in C^1$, $i=1,\ldots, N$, $j,l=\dim
H_{m-1}+1,\ldots,\dim H_m$, $k=\dim H_{m}+1,\ldots,\dim H_{m+1}$.
\end{ex}

\begin{assump}
\label{vectorfields} Throughout the paper, we assume that all the
basis vector fields $X_1,\dots,X_N$ are $C^{1,\alpha}$-smooth,
and, consequently, their commutators are $H^{\alpha}$-continuous,
$\alpha\in[0,1]$.

In some parts of this paper, we consider cases when the
derivatives of the basis vector fields are $H^{\alpha}$-continuous
with respect to some nonnegative symmetric function $\mathfrak
d:U\times U\to\mathbb R$, $U\Subset\mathbb M$, such that
$\mathfrak d\geq C\rho$, $0<C<\infty$, where $C$ depends only on
$U$, and $\rho$ is Riemannian distance. Some additional properties
of $\mathfrak d$ are described below when it is necessary.
\end{assump}

\begin{notat}\label{zero}
In the paper:

\begin{enumerate}
\item The symbol $X\in C^{1,0}$ means that $X\in C^1$, and the
symbol $X\in C^0$ means $X\in C$.

\item $0$-H\"{o}lder continuity means the ordinary continuity. We
denote a modulus of continuity of a mapping $f$ by
$\omega_f(\delta)$.

\item The Riemannian distance is denoted by the symbol $\rho$.
\end{enumerate}
\end{notat}

\begin{thm}\label{Liealg} The coefficients $\bar{c}_{ijk}=c_{ijk}(u)$ of \eqref{tcomm}
with $\operatorname{deg}
X_i+\operatorname{deg}X_j=\operatorname{deg}X_k$ define a graded
nilpotent Lie algebra.
\end{thm}
\begin{proof}
Fix an arbitrary point $u\in\mathbb M$ and show that the
collection $\{c_{ijk}(u)\}$  with ${\operatorname{deg}
X_k=\operatorname{deg} X_i+\operatorname{deg} X_j}$, enjoy the
Jacobi identity and, thus, define the structure of a Lie algebra.
The property $\bar{c}_{ijk}=-\bar{c}_{jik}$ is evident. Prove that
the collection $\{\bar{c}_{ijk}\}$ under consideration enjoys
Jacobi identity.

$1^{\text{\sc st}}$ {\sc Step.} We may assume without loss of
generality that $X_1,\ldots, X_N$ are the vector fields on an open
set of  $\mathbb R^N$ (otherwise, consider the local
$C^2$-diffeomorphism $\eta$ similarly to Remark {\ref{cijkcont}}).

For a vector field
$X_i(x)=\sum\limits_{j=1}^N\eta_{ij}(x)\partial_j$, consider the
mollification
$(X_i)_h(x)=\sum\limits_{j=1}^N(\eta_{ij}*\omega_h)(x)\partial_j$,
$i=1,\ldots, N$, where the function $\omega\in
C_0^{\infty}(B(0,1))$ is such that
$\int\limits_{B(0,1)}\omega(x)\,dx=1$, and
$\omega_h(x)=\frac{1}{h^N}\omega\bigl(\frac{x}{h}\bigr)$. By the
properties of mollification $\eta_{ij}*\omega_h$,
$i,j=1,\ldots,N$, we have
$(X_i)_h\underset{h\to0}{\overset{C^1}\longrightarrow}X_i$ locally
in some neighborhood of~$u$. Note that the vector fields
$(X_i)_h(v)$, $i=1,\ldots, N$, meet the Jacobi identity, and are a
basis of $T_v\mathbb M$ for $v$ belonging to some neighborhood of
$u$, if the parameter $h$ is small enough. Consequently, setting
$[(X_i)_h,(X_j)_h]=\sum\limits_{k=1}^N(c_{ijk})_h(X_k)_h$, we have
\begin{multline*}
\sum\limits_k\sum\limits_l{(c_{ijk})_h(c_{kml})_h}(X_l)_h
+\sum\limits_k\sum\limits_l{(c_{mik})_h(c_{kjl})_h}(X_l)_h\\
+\sum\limits_k\sum\limits_l{(c_{jmk})_h(c_{kil})_h}(X_l)_h
-\sum\limits_l[(X_m)_h(c_{ijl})_h](X_l)_h\\
-\sum\limits_l[(X_j)_h(c_{mil})_h](X_l)_h
-\sum\limits_l[(X_i)_h(c_{jml})_h](X_l)_h=0.
\end{multline*}
Note that, since
$(X_i)_h\underset{h\to0}{\overset{C^1}\longrightarrow}X_i$
locally, and the vector fields $\{(X_i)_h\}_{i=1}^N$ are linearly
independent for all $h\geq0$ small enough, we have $(c_{ijk})_h\to
c_{ijk}$ as $h\to0$.

Now, fix $1\leq l\leq N$. Since the vector fields
$\{(X_i)_h\}_{i=1}^N$  are linearly independent for $h>0$ small
enough, we have
\begin{multline}\label{jacid1}
\sum\limits_{k}(c_{ijk})_h(c_{kml})_h+\sum\limits_{k}{(c_{mik})_h(c_{kjl})_h}
+\sum\limits_{k}{(c_{jmk})_h(c_{kil})_h}\\
-[(X_m)_h(c_{ijl})_h] -[(X_j)_h(c_{mil})_h]-[(X_i)_h(c_{jml})_h]=0
\end{multline}
for each fixed $l$ in some neighbourhood of $u$. Fix $i,j,m$ and
$l$ such that $\operatorname{deg} X_l=\operatorname{deg}
X_i+\operatorname{deg} X_j+\operatorname{deg} X_m$, and  consider
a test function $\varphi\in C_0^{\infty}(U)$ on some small compact
neighborhood $U\ni u$, $U\Subset \mathbb M$. We multiply both
sides of \eqref{jacid1} on $\varphi$ and integrate the result
over $U$. For $h>0$ small enough, we have
\begin{multline*}
0=\int\limits_U\Bigl[\sum\limits_{k}(c_{ijk})_h(v)(c_{kml})_h(v)
+\sum\limits_{k}{(c_{mik})_h(v)(c_{kjl})_h(v)}\\
+\sum\limits_{k}{(c_{jmk})_h(v)(c_{kil})_h(v)}\Bigr]\cdot\varphi(v)\,dv
-\int\limits_U[(X_m)_h(c_{ijl})_h](v)\cdot\varphi(v)\,dv\\
-\int\limits_U[(X_j)_h(c_{mil})_h](v)\cdot\varphi(v)\,dv
-\int\limits_U[(X_i)_h(c_{jml})_h](v)\cdot\varphi(v)\,dv.
\end{multline*}
Show that, among the last three integrals, the first one tends to
zero as $h\to0$. Indeed,
$$
\int\limits_U[(X_m)_h(c_{ijl})_h](v)\cdot\varphi(v)\,dv=
-\int\limits_U[(X_m)^*_h\varphi](v) \cdot(c_{ijl})_h(v)\,dv,
$$
where $(X_i)^*_h$ is an adjoint operator to $(X_i)_h$. The
right-hand part integral tends to zero as $h\to0$ since the value
$[(X_m)^*_h\varphi](v)$ is uniformly bounded in $U$ as $h\to0$,
and $(c_{ijl})_h(v)\to0$ as $h\to0$ in view of the choice of~$l$.
The similar conclusion is true regarding the last two integrals.

Consequently, taking into account the facts that $(c_{ijk})_h\to
c_{ijk}$ locally, and $c_{ijk}=0$ for $\operatorname{deg}
X_k>\operatorname{deg} X_i+\operatorname{deg} X_j$, and using
du~Bois--Reymond Lemma for $h\to0$ we infer
\begin{multline}\label{jacid2}
\sum\limits_{k:\ \operatorname{deg} X_k\leq\operatorname{deg}
X_i+\operatorname{deg} X_j}c_{ijk}(v)c_{kml}(v)+\sum\limits_{k:\
\operatorname{deg} X_k\leq\operatorname{deg}
X_m+\operatorname{deg}
X_i}{c_{mik}(v)c_{kjl}(v)}\\+\sum\limits_{k:\ \operatorname{deg}
X_k\leq\operatorname{deg} X_j+\operatorname{deg}
X_m}{c_{jmk}(v)c_{kil}(v)}=0
\end{multline}
for all $v\in\mathbb M$ close enough to $u$.

$2^{\text{\sc nd}}$ {\sc Step.} For fixed $l$, such that
$\operatorname{deg} X_l=\operatorname{deg} X_i+\operatorname{deg}
X_j+\operatorname{deg} X_m$, investigate the properties of the
index $k$. Consider the first sum. Since $\operatorname{deg}
X_l\leq \operatorname{deg} X_k+\operatorname{deg} X_m$, we have
$\operatorname{deg} X_k\geq \operatorname{deg}
X_l-\operatorname{deg} X_m=\operatorname{deg}
X_i+\operatorname{deg} X_j$. By {\eqref{tcomm}},
$\operatorname{deg} X_k\leq \operatorname{deg}
X_i+\operatorname{deg} X_j$, and, consequently,
$\operatorname{deg} X_k=\operatorname{deg} X_i+\operatorname{deg}
X_j$. The other two cases are considered similarly. Thus, the sum
{\eqref{jacid2}} with $\operatorname{deg} X_l=\operatorname{deg}
X_i+\operatorname{deg} X_j+\operatorname{deg} X_m$ and $v=u$ is
\begin{multline*}
\sum\limits_{\operatorname{deg} X_k=\operatorname{deg}
X_i+\operatorname{deg}
X_j}c_{ijk}(u)c_{kml}(u)+\sum\limits_{\operatorname{deg}
X_k=\operatorname{deg} X_m+\operatorname{deg}
X_i}{c_{mik}(u)c_{kjl}(u)}\\+\sum\limits_{\operatorname{deg}
X_k=\operatorname{deg} X_j+\operatorname{deg}
X_m}{c_{jmk}(u)c_{kil}(u)}=0.
\end{multline*}
The coefficients $\{\bar{c}_{ijk}=c_{ijk}(u)\}_{\operatorname{deg}
X_k=\operatorname{deg} X_i+\operatorname{deg} X_j}$ enjoy the
Jacobi identity, and, thus, they define the Lie algebra. The
theorem follows.
\end{proof}

We construct the Lie algebra $\mathfrak{g}^u$ from Theorem
{\ref{Liealg}}  as a graded nilpotent Lie algebra of vector fields
$\{(\widehat{X}_i^u)^{\prime}\}_{i=1}^N$ on $\mathbb
R^N$~{\cite{post}}. Thus the relation
$$
[(\widehat{X}_i^u)^{\prime},(\widehat{X}_j^u)^{\prime}]=\sum\limits_{\operatorname{deg}
X_k=\operatorname{deg} X_i+\operatorname{deg}
X_j}c_{ijk}(u)(\widehat{X}_k^u)^{\prime}
$$
holds for the vector fields
$\{(\widehat{X}_i^u)^{\prime}\}_{i=1}^N$ everywhere on $\mathbb
R^N$.

\begin{notat}
We use the following standard notations: for each $N$-dimensional
multi-index $\mu=(\mu_1,\ldots,\mu_N)$, its {\it homogeneous norm}
equals $|\mu|_h=\sum\limits_{i=1}^N\mu_i\deg X_i$.
\end{notat}

\begin{defn}
\label{groupoperator} The Carnot group ${\mathbb G}_u\mathbb M$
corresponding to the Lie algebra $\mathfrak{g}^u$, is called the
{\it nilpotent tangent cone} of  $\mathbb M$ at $u\in\mathbb M$.
 We construct ${\mathbb G}_u\mathbb M$ in $\mathbb R^N$ as a groupalgebra
{\cite{post}}.
By Campbell--Hausdorff formula, the group operation is defined for
the  basis vector fields~$(\widehat{X}_i^u)^{\prime}$ on $\mathbb
R^N$, $i=1,\ldots,N$, to be left-invariant~{\cite{post}}: if
$$x=\exp\Bigl(\sum\limits_{i=1}^Nx_i(\widehat{X}_i^u)^{\prime}\Bigr),\ y=\exp\Bigl(\sum\limits_{i=1}^Ny_i(\widehat{X}_i^u)^{\prime}\Bigr)\ \text{then} \
x\cdot
y=z=\exp\Bigl(\sum\limits_{i=1}^Nz_i(\widehat{X}_i^u)^{\prime}\Bigr),
$$
where
\begin{align}
z_i&=x_i+y_i,\quad \deg X_i=1,\notag
\\
z_i&=x_i+y_i+\sum\limits_{\substack{|e_l+e_j|_h=2,\notag\\
l<j}}{F}^i_{e_l,e_j}(u)(x_ly_j-y_lx_j), \quad\deg X_i=2,
\\
z_i&=x_i+y_i+ \sum\limits_{\substack{|\mu+\beta|_h=k,\notag\\
\mu>0,\,\beta>0}}{F}^i_{\mu,\beta}(u) x^\mu\cdot y^\beta
\\\label{group}
&=x_i+y_i+\sum\limits_{\substack{|\mu+e_l+\beta+e_j|_h=k,\\
l<j}} {G}^i_{\mu,\beta,l,j}(u)
x^{\mu}y^{\beta}(x_ly_j-y_lx_j),\quad \deg X_i=k.
\end{align}
\end{defn}

Theorem {\ref{Liealg}} implies

\begin{thm}
[{\cite{fs}}]\label{propfs} If $\{\frac{\partial}{\partial
x_l}\}_{l=1}^N$ is the standard basis in  $\mathbb R^N$ then the
$j$-th coordinate of a vector field
$(\widehat{X}_i^u)^{\prime}(x)=\sum\limits_{j=1}^Nz^{j}_i(u,x)
\frac{\partial}{\partial x_j}$ can be written as
\begin{equation}\label{mat_tc}
z^{j}_i(u,x)=\begin{cases}\delta_{ij}\quad&\text{if}\quad j\leq
\dim H_{\deg X_i},
\\
\sum\limits_{\substack{|\mu+e_i|_h=\deg X_j,\\
\mu>0}}{F}_{\mu,e_i}^j(u)x^\mu \quad&\text{if}\quad j> \dim
H_{\deg X_i}.
\end{cases}
\end{equation}
\end{thm}

\begin{defn}
Suppose that $u\in\mathbb M$ and $(v_1,\ldots, v_N)\in B_E(0, r)$
where $B_E(0,r)$ is an Euclidean ball in $\mathbb R^N$. Define a
mapping $\theta_u(v_1,\ldots, v_N):B_E(0,r)\to\mathbb M$ as
follows:
$$
\theta_u(v_1,\ldots,
v_N)=\exp\biggl(\sum\limits_{i=1}^Nv_iX_i\biggr)(u).
$$
It is known, that $\theta_u$ is a $C^1$-diffeomorphism if $0<r\leq
r_u$ for some $r_u>0$. The collection $\{v_i\}_{i=1}^N$ is called
{\it the normal coordinates} or {\it the coordinates of the
$1^{\text{st}}$ kind  $($with respect to $u\in\mathbb M)$} of the
point $v=\theta_u(v_1,\ldots, v_N)$.
\end{defn}

\begin{assump}
\label{thetav} The compactly embedded neighborhood
$U\subset\mathbb M$ under consideration is such that
$\theta_u(B_E(0,r_u))\supset U$ for all $u\in U$.
\end{assump}

By means of the exponential map we can push-forward the vector
fields $(\widehat{X}_i^u)^{\prime}$ onto $U$ for obtaining the
vector fields
$\widehat{X}_i^u=(\theta_u)_*(\widehat{X}_i^u)^{\prime}$ where
$$(\theta_u)_*\langle Y\rangle(\theta_u(x))=D\theta_u(x)\langle
Y\rangle,$$ $Y\in T_x\mathbb R^N$. Note that
$\widehat{X}_i^u(u)=X_i(u)$. Indeed, on the one hand, by the
definition, we have $(\theta_u)^{-1}_*{X}_i(0)=e_i$. On the other
hand, Theorem {\ref{propfs}} implies
$(\widehat{X}_i^u)^{\prime}(0)=e_i$. Thus
$\widehat{X}_i^u(u)=X_i(u)$.

\begin{thm}
\label{tcone} The vector fields $\widehat{X}_i^u$, $i=1,\ldots,
N$, are locally $H^{\alpha}$-continu\-ous on~$u$.
\end{thm}

The proofs of this theorem and of many other assertions concerning
smoothness use often the following lemma (see its proof in
Section~{\ref{proofdiff}}).

\begin{thm}
\label{ODE} Consider the ODE
\begin{equation}\label{diffequation}
\begin{cases}
\frac{dy}{dt}=f(y,v,u),\\
y(0)=0
\end{cases}
\end{equation}
where $t\in[0,1]$, $y,v,u\in W\subset\mathbb R^N$ and
$\operatorname{Lip}_y(f)=L<1$.

\begin{enumerate}
\item If the mapping $f(y,v,u)=f(y,u)\in C^1(y)\cap C^{\alpha}(u)$
then the solution $y(t,u)\in C^{\alpha}(u)$ locally.

\item If $f(y,v,u)\in C^{1,\alpha}(y,u)\cap C^{1}(v)$ and
$\frac{\partial f}{\partial v}\in C^{1,\alpha}(y,u)$ then
$\frac{dy(t,v,u)}{dv}\in C^{\alpha}(u)$ locally.
\end{enumerate}
\end{thm}

\begin{rem}\label{remODE}
The following statements are proved similarly to Theorem
{\ref{ODE}}.
\begin{enumerate}
\item If the mapping $f(y,v,u)$ from {\eqref{diffequation}} does
not depend on $v$, and it is $C^1$-smooth in $y$ and it is locally
$\alpha$-H\"{o}lder with respect to nonnegative symmetric function
$\mathfrak d$ defined on $U\times U$, $U\Subset\mathbb M$, such
that $\mathfrak d\geq C\rho$, $0<C<\infty$, where $C$ depends only
on $U$, then the solution $y(t,u)$ is also locally
$\alpha$-H\"{o}lder with respect to $\mathfrak d$.

\item If $f(y,v,u)\in C^{1}(y,u)\cap C^{1}(v)$, its derivatives in
$y$ and in $u$ are locally $\alpha$-H\"{o}lder with respect to
$\mathfrak d$, $\frac{\partial f}{\partial v}\in C^{1}(y,u)$, and
the derivatives of $\frac{\partial f}{\partial v}$ in $y$ and in
$u$ are locally $\alpha$-H\"{o}lder with respect to $\mathfrak d$
then $\frac{dy(t,v,u)}{dv}$ is locally $\alpha$-H\"{o}lder with
respect to $\mathfrak d$.
\end{enumerate}
\end{rem}

\begin{rem}
One of particular cases of $\mathfrak d$ is $d_{\infty}$.
\end{rem}

\begin{notat}
Hereinafter, we denote a nonnegative symmetric function defined on
$U\times U$, $U\Subset\mathbb M$, possessing properties from
item~1 of Remark~{\ref{remODE}}, by $\mathfrak d$.
\end{notat}

\begin{proof}[Proof of Proposition {\ref{tcone}}]
$1^{\text{\sc st}}$ {\sc Step.} Taking into account Assumption
{\ref{vectorfields}}, we have the table
$$
[(\widehat{X}^u_i)^{\prime},(\widehat{X}^u_j)^{\prime}](v)=\sum\limits_{\operatorname{deg}
X_k=\operatorname{deg} X_i+\operatorname{deg}
X_j}c_{ijk}(u)(\widehat{X}^u_k)^{\prime}(v).
$$
 By means of Assumption
{\ref{vectorfields}} and Definition {\ref{carnotmanifold}}, the
functions $c_{ijk}(u)$ from {\eqref{tcomm}} are
$H^{\alpha}$-continuous.

If  $X=\sum\limits_{i=1}^Nx_i(\widehat{X}^u_i)^{\prime}$ and
$Y=\sum\limits_{i=1}^Ny_i(\widehat{X}^u_i)^{\prime}$ then by
Campbell--Hausdorff formula we have $\exp tY\circ\exp tX(g)=\exp
Z(t)(g)$ where $Z(t)=tZ_1+t^2Z_2+\ldots+t^MZ_M$ and $Z_1$,
$Z_2,\dots$ are some vector fields independent of $t$. Dynkin
formula (see, for instance, {\cite{post}}) for calculating
$Z_l(t)$, $1\leq l\leq M$, gives
\begin{multline*}
Z_l= \frac1{n}
\sum\limits_{k=1}^l\frac{(-1)^{k-1}}{k}\sum\limits_{(p)(q)}
\frac{(\operatorname{ad} Y)^{q_k}(\operatorname{ad} X)^{p_k}\dots
(\operatorname{ad} Y)^{q_1}(\operatorname{ad}
X)^{p_1-1}X}{p_1!q_1!\dots p_k!q_k!}\\
= \sum\limits_{(p)(q)}C_{(p)(q)} (\operatorname{ad}
Y)^{q_k}(\operatorname{ad} X)^{p_k}\dots (\operatorname{ad}
Y)^{q_1}(\operatorname{ad} X)^{p_1-1}X,
\end{multline*}
where $C_{(p)(q)}=\operatorname{const}$, $(p)=(p_1,\dots,p_k)$,
$(q)=(q_1,\dots,q_k)$. We sum over all natural $p_1$, $q_1$,
$\dots p_k$, $q_k$, such that $p_i+q_i>0$,
$p_1+q_1+\dots+p_k+q_k=l$, and $(\operatorname{ad} A)B=[A,B]$,
$(\operatorname{ad} A)^0B=B$. The each summand can be represented
as a sum
$$
Z_l(v)=\sum\limits_{i=1}^Nd_{j,l}(u,
x,y)(\widehat{X}^u_j)^{\prime}(v),
$$
where $d_{j,l}(u,x,y)$ are polynomial functions of
$x=(x_1,\ldots,x_N)$, $y=(y_1,\ldots,y_N)$ coefficients of which
are polynomial functions of $\{c_{lmk}(u)\}$ and, consequently,
are H\"{o}lder in~$u$. More exactly,
$$
\sum\limits_{l=2}^M Z_{l}= \sum\limits_{l=2}^M
\sum\limits_{j=1}^Nd_{j,l}(u,
x,y)(\widehat{X}^u_j)^{\prime}=\sum\limits_{j=1}^N\biggl[\sum\limits_{l=2}^M\sum\limits_{\begin{array}{c}\scriptstyle
|\mu+\beta|_h=l,\vspace{-2mm}\\ \scriptstyle \mu>0,\beta>0
\end{array}}
F^j_{\mu,\beta}(u) x^{\mu}\cdot
y^{\beta}\biggr](\widehat{X}^u_j)^{\prime}.
$$
Consequently,
$$
d_{j,l}(u,x,y)=\sum\limits_{l=2}^M\sum\limits_{\begin{array}{c}\scriptstyle
|\mu+\beta|_h=l,\vspace{-2mm}\\ \scriptstyle \mu>0,\beta>0
\end{array}}
F^j_{\mu,\beta}(u) x^{\mu}\cdot y^{\beta}.
$$ Hence,
$F^j_{\mu,\beta}(u)$ are $H^{\alpha}$-continuous in $u$, and
$(\widehat{X}_i^u)^{\prime}$ are also $H^{\alpha}$-continuous on
$u$ (see {\eqref{mat_tc}}).

$2^{\text{\sc nd}}$ {\sc Step.} Consider the following Cauchy
problem:
\begin{equation}\label{problem}
\begin{cases}
\frac{d\Phi(t,u,\xi)}{d t}=\sum\limits_{i=1}^N\xi_i
X_i(\Phi),\\
\Phi(0,u,\xi)=u,
\end{cases}
\end{equation}
where $\xi=(\xi_1,\ldots,\xi_N)$. Note that
$\Phi(t,u,\xi)=\exp\Bigl(\sum\limits_{i=1}^Nt\xi_iX_i\Bigr)(u)$.
We can assume without loss of generality, that $\mathbb M=\mathbb
R^N$. If Assumption {\ref{vectorfields}} holds then the mapping
$f(\xi,\Phi)=\sum\limits_{i=1}^N\xi_i X_i(\Phi)$ is $C^{1,
\alpha}$-smooth in $\xi$ and $\Phi$. From the definition, it
follows, that $\theta_u(\xi)=\Phi(1,u,\xi)$.

By theorem {\ref{ODE}} on smooth dependence of ODE solution on
parameters (see Section~{\ref{proofdiff}}
 for details), it is easy to see, that the differential
$D\theta_u(y)$ is $H^{\alpha}$-continuous in $u$.

Since
$\widehat{X}_i^u(x)=D\theta_u(y)(\widehat{X}_i^u)^{\prime}(y)$,
$x=\theta_u(y)$, the proposition follows from results of the
$1^{\text{\sc st}}$ and $2^{\text{\sc nd}}$  steps.
\end{proof}

\begin{rem}\label{nilphold}
If the derivatives of $X_i$, $i=1,\ldots, N$, are locally
H\"{o}lder with respect to $\mathfrak d$, then $\widehat{X}^u_i$,
$i=1,\ldots, N$, are locally H\"{o}lder on $u$ with respect to
$\mathfrak d$.
\end{rem}

\begin{defn}\label{locCargr}
The local Lie group corresponding to the Lie algebra
$\{\widehat{X}_i^u\}_{i=1}^N$, is called the {\it local Carnot
group} ${\mathcal G}^u\mathbb M$ at $u\in\mathbb M$. Define it in
such a way that the mapping $\theta_u$ is a {\it group
isomorphism} of some neighborhood of the unity of the group
$\mathbb G_u\mathbb M$ and ${\mathcal G}^u\mathbb M$. The
canonical Riemannian structure is defined by scalar product at the
unit of ${\mathcal G}^u\mathbb M$ coinciding with those in
$T_u\mathbb M$.
\end{defn}

\begin{rem}\label{locCarop}
Recall that the vector fields $\widehat{X}_i^u$, $i=1,\ldots, N$,
are locally $H^{\alpha}$-continuous on $\mathbb M$,
$\alpha\in[0,1]$. The exponential mapping
$\exp\Bigl(\sum\limits_{i=1}^Na_i\widehat{X}_i^u\Bigr)(g)$ is not
defined  correctly for such fields. Therefore, in view of
smoothness of $(\theta^{-1}_u)_*(\widehat X_i^u)$, $i=1,\ldots,
N$, we define the point
$$
a=\exp\biggl(\sum\limits_{i=1}^Na_i\widehat{X}_i^u\biggr)(g)
$$
according to Definition \ref{locCargr}: first, we obtain a point
$$
a_{u}=\exp\biggl(\sum\limits_{i=1}^Na_i\cdot
(\theta^{-1}_u)_*(\widehat{X}_i^u)\biggr)(\theta_u^{-1}(g)),
$$
and then we define $a=\theta_u(a_{u})$. Moreover, we similarly
define the whole curve corresponding to this exponential mapping.
Suppose that
$$
\begin{cases}
\dot{\gamma_u}(t)=\sum\limits_{i=1}^Na_i\cdot
(\theta^{-1}_u)_*(\widehat{X}_i^u)(\gamma_u(t))\\
\gamma_u(0)=\theta_u^{-1}(g).
\end{cases}
$$
Then, for the curve $\gamma(t)=\theta_u(\gamma_u(t))$, we have
$$
\begin{cases}
\dot{\gamma}(t)=\sum\limits_{i=1}^Na_i\widehat{X}_i^u(\gamma(t))\\
\gamma(0)=g.
\end{cases}
$$

In particular, we have:

\begin{enumerate}
\item The exponential mapping $
\widehat{\theta}_u(v_1,\ldots,v_n)=\exp\Bigl(\sum\limits_{i=1}^Nv_i\widehat{X}^u_i\Bigr)(u)
$ is defined as
$$
\theta_u\biggl[\exp\biggl(\sum\limits_{i=1}^Nv_i(\widehat{X}_i^u)^{\prime}\biggr)(0)\biggr];
$$
and the mapping $
\widehat{\theta}_u^w(v_1,\ldots,v_n)=\exp\Bigl(\sum\limits_{i=1}^Nv_i\widehat{X}^u_i\Bigr)(w)
$ is defined as
$$
\theta_u\biggl[\exp\biggl(\sum\limits_{i=1}^Nv_i(\widehat{X}_i^u)^{\prime}\biggr)(\theta_u^{-1}(w))\biggr].
$$

\item The inverse mapping $\exp^{-1}$ is also defined by the
unique way for vector fields $\{\widehat{X}_i^u\}_{i=1}^N$ since
it is defined by the unique way for
$\{(\widehat{X}_i^u)^{\prime}\}_{i=1}^N$.

\item The group operation is defined by the following way: if
$x=\exp\Bigl(\sum\limits_{i=1}^Nx_i\widehat{X}^u_i\Bigr)$,
$y=\exp\Bigl(\sum\limits_{i=1}^Ny_i\widehat{X}^u_i\Bigr)$ then
$x\cdot
y=\exp\Bigl(\sum\limits_{i=1}^Ny_i\widehat{X}^u_i\Bigr)\circ
\exp\Bigl(\sum\limits_{i=1}^Nx_i\widehat{X}^u_i\Bigr)
=\exp\Bigl(\sum\limits_{i=1}^Nz_i\widehat{X}^u_i\Bigr)$ where
$z_i$ are taken from Definition {\ref{groupoperator}}.

\item Using the normal coordinates $\widehat\theta_u^{-1}$, define
the action of the {\it dilation group} $\delta^u_\varepsilon$ in
the local Carnot group $\mathcal G^u\mathbb M$: to an element
$x=\exp\Bigl(\sum\limits_{i=1}^Nx_i\widehat X_i^u\Bigr)(u)$,
assign $ \delta^u_\varepsilon
x=\exp\Bigl(\sum\limits_{i=1}^Nx_i\varepsilon^{\deg X_i} \widehat
X_i^u\Bigr)(u) $ in the cases where the right-hand side makes
sense.
\end{enumerate}
\end{rem}

\begin{property}\label{proplocCargr}
For each vector field $\widehat{X}_i^u$, $i=1,\ldots, N$, we have
$(\delta^u_\varepsilon)_*\widehat{X}_i^u(g)=\varepsilon^{\operatorname{deg}X_i}\widehat{X}_i^u(\delta^u_\varepsilon
g)$.
\end{property}

This property comes from those on the ``canonical''  Carnot group
$T_u\mathbb M$ \cite{fs}.

\begin{lem}[{\cite{vod1}}]
\label{property_sum}Suppose that $u\in U$. The equality
$$
\sum\limits_{i=1}^{j}
\sum\limits_{\substack{ |\mu+e_i|_h=\deg X_j,\\
|\mu+e_i|=l,\,\mu>0}} x_i F^j_{\mu,e_i}(u)x^\mu=0,\quad
x=(x_1,\ldots,x_N)\in\mathbb R^N,
$$
holds for all $\deg X_j\geq 2$, $l=2,\ldots,\deg X_j$.
\end{lem}

\begin{proof} Consider a vector field $X=\sum\limits_{i=1}^{N}x_i
(\widehat{X}_i{}^u)'$. It is known that $\exp r sX\circ\exp rt
X(g)=\exp r(s+t)X(g)$. Therefore, by \eqref{group}, we have
$$
\sum\limits_{\substack{ |\mu+\beta|_h=\deg X_j,\\
\mu>0,\,\beta>0}} r^{|\mu+\beta|} F^j_{\mu,\beta}(g)s^{|\mu|}
x^\mu\cdot t^{|\beta|}x^\beta =0
$$
for all fixed $s$ and $t$, $\deg X_j\geq2$. It follows that the
coefficients at all powers of $r$ vanish. In particular, if
$|\mu+\beta|=l\geq2$ then
$$
\sum\limits_{\substack{|\mu+\beta|_h=\deg X_j,\\
\mu>0,\,\beta>0,\,|\mu+\beta|=l}} F^j_{\mu,\beta}(g)s^{|\mu|}
x^\mu\cdot t^{|\beta|}x^\beta =0.
$$
Consequently, if $|\beta|=1$ then we infer
$$
P(s)=\sum\limits_{l=2}^{\deg X_j} s^{l-1}\sum\limits_{i=1}^{j}
\sum\limits_{\substack{|\mu+e_i|_h= \deg
X_j,\\|\mu+e_i|=l,\,\mu>0}}x_i F^j_{\mu,e_i}(g)
 x^\mu\equiv0,
$$
where $s$ is an arbitrarily small parameter.  Therefore, all
coefficients of the polynomial $P(s)$ at the powers of $s$ vanish.
The lemma follows.
\end{proof}

\begin{lem}[{\cite{vod1}}]
\label{coincide} Let  $u\in U$ be an arbitrary point. Then
$$
a=\exp\biggl(\sum\limits_{i=1}^Na_iX_i\biggr)(u)=\exp\biggl(\sum\limits_{i=1}^Na_i\widehat{X}^u_i\biggr)(u)
$$
for all $|a_i|<r_u$, $i=1,\ldots, N$.
\end{lem}

\begin{proof}
Lemma {\ref{property_sum}}  implies that the line $\mathbb R\ni
t\mapsto t(a_1,\ldots,a_N)$ is the integral line of the vector
field $\sum\limits_{i=1}^Na_i(\widehat{X}_i^u)^{\prime}$ starting
at 0 as $t=0$.
 By the definition of the exponential map,
 we infer  $\mathbb R^N\ni(a_1,\ldots,a_N)=
\sum\limits_{i=1}^Na_i(\widehat{X}_i^u)^{\prime}(a_1,\ldots,a_N)=
\exp\Bigl(\sum\limits_{i=1}^Na_i(\widehat{X}_i^u)^{\prime}\Bigr)$,
i.~e. the exponential map equals the identity. From this, it
 follows immediately that
\begin{multline*}
a=\theta_u(a_1,\ldots,a_N)
=\theta_u\biggl(\sum\limits_{i=1}^Na_i(\widehat{X}_i^u)^{\prime}\biggr)\\
=\theta_u\biggl(\exp\biggl(\sum\limits_{i=1}^Na_i(\widehat{X}_i^u)^{\prime}\biggr)\biggr)
=\exp\biggl(\sum\limits_{i=1}^Na_i\widehat{X}_i^u\biggr)
\end{multline*}
according to Remark~\ref{locCarop}.
\end{proof}

\begin{defn}
Suppose that $\mathbb M$ is a Carnot manifold, and $u\in\mathbb
M$. For $a,p\in{\mathcal G}^u\mathbb M$, where
$$
a=\exp\biggl(\sum\limits_{i=1}^Na_i\widehat{X}^u_i\biggr)(p),
$$
we define the quasimetric
$d_{\infty}^u(a,p)=\max\limits_{i}\{|a_i|^{\frac{1}
{\operatorname{deg} X_i}}\}$ on ${\mathcal G}^u\mathbb M$.
\end{defn}

The following  properties comes from those on the ``canonical''
Carnot group $T_u\mathbb M$ \cite{fs}.

\begin{property}\label{propmetrlCg}
It is easy to see that $d_{\infty}^u(x,y)$ is a quasimetric on
${\mathcal G}^u\mathbb M$ meeting  the following properties:
\begin{enumerate}
\item $d_\infty^u(x,y)\geq0$, $d_\infty^u(x,y)=0$
 if and only if $x=y$;

\item $d_\infty^u(u,v)=d_\infty^u(v,u)$;

\item the quasimetric $d_\infty^u(x,y)$ is continuous with respect
to each of its variables;

\item there exists a constant $Q_\triangle=C_\triangle(U)$ such
that the inequality
$$
d_\infty^u(x,y)\leq Q_\triangle(d_\infty^u(x,z)+d_\infty^u(z,y))
$$
holds for every triple of points $x$, $y$, $z\in U$.
\end{enumerate}
\end{property}

\begin{property}
Let
$$
w_{\varepsilon}=\exp\biggl(\sum\limits_{i=1}^N\varepsilon^{\operatorname{deg}
X_i}w_i\widehat{X}^u_i\biggr)(v)\text{ and
}g_{\varepsilon}=\exp\biggl(\sum\limits_{i=1}^N\varepsilon^{\operatorname{deg}
X_i}g_i\widehat{X}^u_i\biggr)(v).
$$
Then $d_{\infty}^u(w_{\varepsilon},g_{\varepsilon})=\varepsilon
d_{\infty}^u(w_1,g_1)$.
\end{property}

 By $\operatorname{Box}^u(x,r)$ we denote a set
$\{y\in\mathbb M:\,d_{\infty}^u(x,y)<r\}$.

\begin{property}
We have
$\delta^u_\varepsilon(\operatorname{Box}^u(u,r))=\operatorname{Box}^u(u,\varepsilon
r)$.
\end{property}

\subsection{Gromov's Theorem on the Nilpotentization of Vector
Fields~and Estimate of the Diameter of a Box}\label{GrTh}

\begin{defn}
Suppose that $\mathbb M$ is a Carnot manifold, and let
$U\subset\mathbb M$ be as in Assumption \ref{thetav}. Given
$$
v=\exp\biggl(\sum\limits_{i=1}^Nv_iX_i\biggr)(u)
$$
$u,v\in  U$, define the quasimetric
$d_{\infty}(u,v)=\max\limits_{i}\{|v_i|^{\frac{1}{\operatorname{deg}
X_i}}\}$. By $\operatorname{Box}(x,r)$ we denote a set
$\{y\in\mathbb M:\,d_{\infty}(x,y)< r\}$, $r\leq r_x$.
\end{defn}

\begin{defn}
Using the normal coordinates $\theta_u^{-1}$, define the action of
the {\it dilation group} $\Delta^u_\varepsilon$ in a neighborhood
of a point $u\in\mathbb M$: to an element
$x=\exp\Bigl(\sum\limits_{i=1}^Nx_i X_i\Bigr)(u)$, assign $
\Delta^u_\varepsilon
x=\exp\Bigl(\sum\limits_{i=1}^Nx_i\varepsilon^{\deg X_i}
X_i\Bigr)(u) $ in the cases where the right-hand side makes sense.
\end{defn}

\begin{property}\label{coindil}
By Lemma \ref{coincide} we have $\Delta^u_\varepsilon
x=\delta^u_\varepsilon x$.
\end{property}

\begin{property}
By Lemma \ref{coincide} we have
$\operatorname{Box}^u(u,r)=\operatorname{Box}(u,r)$.
\end{property}

\begin{property}
We have $\Delta^u_\varepsilon
(\operatorname{Box}(u,r))=\operatorname{Box}(u,\varepsilon r)$,
$r\in(0,r_u]$.
\end{property}

\begin{property}\label{propmetr}
The quasimetric $d_{\infty}$ has the following properties:
\begin{enumerate}
\item $d_\infty(u,v)\geq0$, $d_\infty(u,v)=0$
 if and only if $u=v$;

\item $d_\infty(u,v)=d_\infty(v,u)$;

\item the quasimetric $d_\infty(u,v)$ is continuous with respect
to each of its variables;

\item there exists a constant $Q=Q(U)$ such that the inequality
$$
d_\infty(u,v)\leq Q(d_\infty(u,w)+d_\infty(w,v))
$$
holds for every triple of points $u$, $w$, $v\in U$.
\end{enumerate}
\end{property}

\begin{proof}
The proof of properties 1--3 is based on known properties of
solutions to ODE's. We prove the generalized triangle inequality
at the end of current subsection (see Corollary \ref{trin}).
\end{proof}

\begin{thm}\label{th_expr} Let $X_j\in C^{1}$. Fix $u\in\mathbb M$. If $d_{\infty}(u,w)=C\varepsilon$, then
$$
\widehat{X}_j^u(w)=\sum\limits_{k:\,\operatorname{deg}X_k\leq\operatorname{deg}X_j}[\delta_{kj}+O(\varepsilon)]X_k
+\sum\limits_{k:\,\operatorname{deg}X_k>\operatorname{deg}X_j}
o(\varepsilon^{{\operatorname{deg}X_k-\operatorname{deg}X_j}})X_k(w),
$$ $j=1,\ldots, N$. All $o(\cdot)$ are uniform  in $u$
belonging to some compact subset of $U$.
\end{thm}

\begin{proof}
$1^{\text{\sc st}}$ {\sc Step.} Applying the mapping
$\theta_u^{-1}$ to each vector field $\widehat{X}^u_j$,
$j=1,\ldots, N$, we deduce
$$
D\theta_u^{-1}\widehat{X}^u_j(s)=\sum\limits_{k=1}^N z^k_j(s)e_k,
$$
where $\{e_k\}_{k=1}^N$ is the collection of the vectors of the
standard basis in $\mathbb R^N$, and by \eqref{mat_tc}
$$
z^k_j(s)=\delta_{kj}
+\sum\limits_{|\mu|_h=\operatorname{deg}X_k-\operatorname{deg}X_j>0}F^k_{\mu,e_j}(u)s^{\mu}.
$$
Note that, here
$|s^{\mu}|=O(\varepsilon^{\operatorname{deg}X_k-\operatorname{deg}X_j})$,
since
$$d_{\infty}(0,s)=d_{\infty}(\theta_u^{-1}(u),s)= d^u_{\infty}(u,\theta_u(s))=O(\varepsilon).
$$
Then
$$
\widehat{X}^u_j(\theta_u(s))=\sum\limits_{k=1}^N
z^k_j(s)D\theta_u(s)e_k=\sum\limits_{k=1}^N
z^k_j(s)\biggl(X_k(\theta_u(s))+\frac12\biggl[X_k,\sum\limits_{l=1}^Ns_lX_l\biggr](\theta_u(s))\biggr),
$$
since
$D\theta_u(s)e_k=X_k(\theta_u(s))+\frac12\Bigl[X_k,\sum\limits_{l=1}^Ns_lX_l\Bigr](\theta_u(s))$,
where $s=(s_1,\ldots, s_N)$.

To understand the latter, it is enough to consider the following
equalities
\begin{multline*}
\theta_u(s+re_k)=\exp\biggl(\sum\limits_{l=1}^Ns_lX_l+rX_k\biggr)(u)\\=\exp\biggl(\sum\limits_{l=1}^Ns_lX_l+rX_k\biggr)\circ\exp\biggl(-\sum\limits_{l=1}^Ns_lX_l\biggr)\circ\exp\biggl(\sum\limits_{l=1}^Ns_lX_l\biggr)(u)\\=\exp\biggl(rX_k+\frac
r2\biggl[X_k,\sum\limits_{l=1}^Ns_lX_l\biggr]+o(r)\biggr)(\theta_u(s)),
\end{multline*}
and note (see justification of this calculation for $C^1$-vector
fields in \cite{agrgam})  that
$$
D\theta_u(s)e_k=\frac{\partial}{\partial
r}\theta_u(s+re_k)\Bigl|_{r=0}=X_k(\theta_u(s))+\frac12\biggl[X_k,\sum\limits_{l=1}^Ns_lX_l\biggr](\theta_u(s)).
$$

In view of the properties of the point $s$, we get
$|s_l|=O(\varepsilon^{\operatorname{deg}X_l})$, $l=1,\ldots, N$.
Moreover, taking into account the definition of a Carnot manifold,
we have
$$
\biggl[X_k,\sum\limits_{l=1}^Ns_lX_l\biggr](\theta_u(s))=
\sum\limits_{l=1}^N
\sum\limits_{m:\operatorname{deg}X_m\leq\operatorname{deg}X_k+\operatorname{deg}X_l}
c_{klm}(\theta_u(s))X_m(\theta_u(s)).
$$
Consequently,
\begin{multline*}
\widehat{X}^u_j(\theta_u(s))=\sum\limits_{k=1}^N
z^k_j(s)X_k(\theta_u(s))\\+\frac12\sum\limits_{k=1}^N\sum\limits_{l=1}^N\sum\limits_{\operatorname{deg}X_m\leq\operatorname{deg}X_k+\operatorname{deg}X_l}z^k_j(s)s_lc_{klm}(\theta_u(s))X_m(\theta_u(s))\\
=\sum\limits_{k=1}^N
\Bigl[z^k_j(s)+\frac12\sum\limits_{m,l:\operatorname{deg}X_k
\leq\operatorname{deg}X_m+\operatorname{deg}X_l}z^m_j(s)s_lc_{mlk}
(\theta_u(s))\Bigr]X_k(\theta_u(s)).
\end{multline*}
where
$\bigl|z^m_j(s)\big|=O(\varepsilon^{\operatorname{deg}X_m-\operatorname{deg}X_j})$
and
\begin{equation}\label{calc1}
\bigl|z^m_j(s)s_l\bigr|=O(\varepsilon^{\operatorname{deg}X_m+\operatorname{deg}X_l-\operatorname{deg}X_j}).
\end{equation}
Represent the last sum as
\begin{multline}\label{sum_rep}
\sum\limits_{k:\,\operatorname{deg}X_k<\operatorname{deg}X_j}
\Bigl[z^k_j(s)+\frac12\sum\limits_{m,l:\operatorname{deg}X_k
\leq\operatorname{deg}X_m+\operatorname{deg}X_l}z^m_j(s)s_lc_{mlk}
(\theta_u(s))\Bigr]X_k(\theta_u(s))\\
+\sum\limits_{k:\,\operatorname{deg}X_k=\operatorname{deg}X_j}
\Bigl[z^k_j(s)+\frac12\sum\limits_{m,l:\operatorname{deg}X_j
\leq\operatorname{deg}X_m+\operatorname{deg}X_l}z^m_j(s)s_lc_{mlj}
(\theta_u(s))\Bigr]X_j(\theta_u(s))\\
+\sum\limits_{k:\,\operatorname{deg}X_k>\operatorname{deg}X_j}
\Bigl[z^k_j(s)+\frac12\sum\limits_{m,l:\operatorname{deg}X_k
\leq\operatorname{deg}X_m+\operatorname{deg}X_l}z^m_j(s)s_lc_{mlk}
(\theta_u(s))\Bigr]X_k(\theta_u(s)).
\end{multline}
Note that, we have $z^k_j(s)=0$ if $k<j$. Next, if $k<j$ and
$\operatorname{deg}X_k
=\operatorname{deg}X_m+\operatorname{deg}X_l$, we have $m<j$ and
$z^m_j(s)=0$. Thus, for the first sum equals
$$
\sum\limits_{k:\,\operatorname{deg}X_k<\operatorname{deg}X_j}
\Bigl[\frac12\sum\limits_{m,l:\operatorname{deg}X_k
<\operatorname{deg}X_m+\operatorname{deg}X_l}z^m_j(s)s_lc_{mlk}
(\theta_u(s))\Bigr]X_k(\theta_u(s)).
$$
Similarly, for the second sum we have $z^k_j(s)=\delta_{kj}$, and
if $\operatorname{deg}X_j
=\operatorname{deg}X_m+\operatorname{deg}X_l$ then $z^m_j(s)=0$
since this relation implies $m<j$. Thus, we obtain
$$
\sum\limits_{k:\,\operatorname{deg}X_k=\operatorname{deg}X_j}\Bigl[\delta_{kj}+\frac12\sum\limits_{m,l:\operatorname{deg}X_j
<\operatorname{deg}X_m+\operatorname{deg}X_l}z^m_j(s)s_lc_{mlj}
(\theta_u(s))\Bigr]X_j(\theta_u(s)).
$$
In the third sum, the functions $z^k_j(s)$ and $z^m_j(s)$ can take
any possible values.

$2^{\text{\sc nd}}$ {\sc Step.} Now, we calculate  more exact
estimates of \eqref{calc1}.

\begin{itemize}

\item Let $\operatorname{deg}X_k>\operatorname{deg}X_j$ and
$\operatorname{deg}X_k=\operatorname{deg}X_m+\operatorname{deg}X_l$.
From the above estimate we infer
$$
\bigl|z^m_j(s)s_l\bigr|
=O(\varepsilon^{\operatorname{deg}X_k-\operatorname{deg}X_j}).
$$

Next, suppose that $\operatorname{deg}X_k>\operatorname{deg}X_j$
and
$\operatorname{deg}X_k<\operatorname{deg}X_m+\operatorname{deg}X_l$.
Then all the situations
$\operatorname{deg}X_m>\operatorname{deg}X_j$,
$\operatorname{deg}X_m=\operatorname{deg}X_j$ and
\linebreak$\operatorname{deg}X_m<\operatorname{deg}X_j$ are
possible. Here we have
$$
\bigl|z^m_j(s)s_l\bigr|=
\begin{cases}
\varepsilon O(\varepsilon^{\operatorname{deg}X_l})\leq\varepsilon
O(\varepsilon^{\operatorname{deg}X_k-\operatorname{deg}X_j})&\text{
if }\operatorname{deg}X_m>\operatorname{deg}X_j,\\
O(\varepsilon^{\operatorname{deg}X_l})\leq \varepsilon
O(\varepsilon^{\operatorname{deg}X_k-\operatorname{deg}X_j})&\text{
if }\operatorname{deg}X_m=\operatorname{deg}X_j,\\
0&\text{ if }\operatorname{deg}X_m<\operatorname{deg}X_j.
\end{cases}
$$

\item Let now $\operatorname{deg}X_k=\operatorname{deg}X_j$ and
$\operatorname{deg}X_k<\operatorname{deg}X_m+\operatorname{deg}X_l$.
We again have to consider the situations
$\operatorname{deg}X_m>\operatorname{deg}X_j$,
$\operatorname{deg}X_m=\operatorname{deg}X_j$ and
$\operatorname{deg}X_m<\operatorname{deg}X_j$. It follows
$$
\bigl|z^m_j(s)s_l\bigr|=
\begin{cases}
\varepsilon O(\varepsilon^{\operatorname{deg}X_l})\leq\varepsilon
O(1)&\text{
if }\operatorname{deg}X_m>\operatorname{deg}X_j,\\
O(\varepsilon^{\operatorname{deg}X_l})\leq \varepsilon O(1)&\text{
if }\operatorname{deg}X_m=\operatorname{deg}X_j,\\
0&\text{ if }\operatorname{deg}X_m<\operatorname{deg}X_j.
\end{cases}
$$

\item Finally, let $\operatorname{deg}X_k<\operatorname{deg}X_j$
and
$\operatorname{deg}X_k<\operatorname{deg}X_m+\operatorname{deg}X_l$.
In three situations $\operatorname{deg}X_m>\operatorname{deg}X_j$,
$\operatorname{deg}X_m=\operatorname{deg}X_j$ and
$\operatorname{deg}X_m<\operatorname{deg}X_j$, we obtain the same
result as in the previous case:
$$
\bigl|z^m_j(s)s_l\bigr|=
\begin{cases}
\varepsilon O(\varepsilon^{\operatorname{deg}X_l})\leq\varepsilon
O(1)&\text{
if }\operatorname{deg}X_m>\operatorname{deg}X_j,\\
O(\varepsilon^{\operatorname{deg}X_l})\leq \varepsilon O(1)&\text{
if }\operatorname{deg}X_m=\operatorname{deg}X_j,\\
0&\text{ if }\operatorname{deg}X_m<\operatorname{deg}X_j.
\end{cases}
$$
\end{itemize}

Thus, in the first sum of {\eqref{sum_rep}}, the coefficients at
$X_k$ equal $O(\varepsilon)$, and in the second sum the
coefficient at $X_j$ equals $1+O(\varepsilon)$, and the
coefficients at $X_k$ for $k\neq j$ equal~$O(\varepsilon)$.

$3^{\text{\sc rd}}$ {\sc Step.} Consider the last sum (where
$\operatorname{deg}X_k>\operatorname{deg}X_j$). Note that,
\begin{equation}\label{cmlkinc}
c_{mlk}(\theta_u(s))=c_{mlk}(u)+o(1).
\end{equation}
Then, taking into account \eqref{calc1} and the results of the
$2^{\text{nd}}$ step, we deduce
\begin{multline}\label{sumcmlk}
\sum\limits_{m,l:\operatorname{deg}X_k\leq\operatorname{deg}X_m+\operatorname{deg}X_l}z^m_j(s)s_lc_{mlk}(\theta_u(s))\\=
\sum\limits_{m,l:\operatorname{deg}X_k=\operatorname{deg}X_m+\operatorname{deg}X_l}z^m_j(s)s_lc_{mlk}(\theta_u(s))\\
+\sum\limits_{m,l:\operatorname{deg}X_k<\operatorname{deg}X_m+\operatorname{deg}X_l}z^m_j(s)s_lc_{mlk}(\theta_u(s))\\
=\sum\limits_{m,l:\operatorname{deg}X_k=\operatorname{deg}X_m+\operatorname{deg}X_l}z^m_j(s)s_lc_{mlk}(u)
+o(1)\cdot\varepsilon^{\operatorname{deg}X_k-\operatorname{deg}X_j}\\+\varepsilon\cdot
O(\varepsilon^{\operatorname{deg}X_k-\operatorname{deg}X_j})\\=
\sum\limits_{m,l:\operatorname{deg}X_k=\operatorname{deg}X_m+\operatorname{deg}X_l}z^m_j(s)s_lc_{mlk}(u)+o(1)\cdot
\varepsilon^{\operatorname{deg}X_k-\operatorname{deg}X_j}.
\end{multline}
Consequently,
\begin{multline*}
\widehat{X}^u_j(\theta_u(s))=\sum\limits_{k:\,\operatorname{deg}X_k\leq\operatorname{deg}X_j}[\delta_{kj}+O(\varepsilon)]X_k\\
+\sum\limits_{k:\,\operatorname{deg}X_k>\operatorname{deg}X_j}
\Bigl[z^k_j(s)+\frac12\sum\limits_{m,l}z^m_j(s)s_lc_{mlk}(u)
+o(\varepsilon^{\operatorname{deg}X_k-\operatorname{deg}X_j})\Bigr]X_k(\theta_u(s)),
\end{multline*}
where $m,l$ in the last sum are such that
$\operatorname{deg}X_k=\operatorname{deg}X_m+\operatorname{deg}X_l$.

$4^{\text{\sc th}}$ {\sc Step.} It only remains to show that
\begin{equation}\label{100}
z^k_j(s)+\frac12\sum\limits_{m,l:\operatorname{deg}X_k=\operatorname{deg}X_m+\operatorname{deg}X_l}z^m_j(s)s_lc_{mlk}(u)=\delta_{kj}.
\end{equation}
For obtaining this, consider the mapping
$\widehat{\theta}_u(x)=\exp\Bigl(\sum\limits_{i=1}^Nx_i\widehat{X}_i^u\Bigr)(u)=\theta_u(x)$,
and apply the arguments of the $1^{\text{st}}$ step with the
following difference: it is known, that the vector fields
$\widehat{X}_i^u$, $i=1,\ldots N$, are continuous, but they may
not be differentiable, and formally, we cannot consider
commutators of such vector fields. Therefore we modify previous
arguments. For doing this, we consider the following
representation of the identical mapping:
$$
\widehat{\theta}_0(s)=\exp\biggl(\sum\limits_{i=1}^Ns_iD\widehat{\theta}_u^{-1}(\widehat{X}_i^u)\biggr)(0)=s,
$$
and represent $e_k=D\widehat{\theta}_0(s)(e_k)$ as before we
represented $D{\theta}_u(s)(e_k)$. It is possible, since the
vector fields $D\widehat{\theta}_u^{-1}(\widehat{X}_i^u)$,
$i=1,\ldots, N$, are smooth. Similarly to the $1^{\text{st}}$
step, we infer
$$
D\widehat{\theta}_0(s)(e_k)=D\widehat{\theta}_u^{-1}(\widehat{X}^u_k)(\widehat{\theta}_0(s))+\frac12\biggl[D\widehat{\theta}_u^{-1}(\widehat{X}^u_k),\sum\limits_{l=1}^Ns_lD\widehat{\theta}_u^{-1}(\widehat{X}^u_l)\biggr](\widehat{\theta}_0(s)).
$$
Since $\widehat\theta_0(s)=s$ and in view of properties of the
vector fields $D\widehat{\theta}_u^{-1}(\widehat{X}_i^u)$,
$i=1,\ldots, N$, we deduce
\begin{multline*}
e_k=D\widehat{\theta}_u^{-1}(\widehat{X}^u_k)(s)+\frac12\biggl[D\widehat{\theta}_u^{-1}(\widehat{X}^u_k),\sum\limits_{l=1}^Ns_lD\widehat{\theta}_u^{-1}(\widehat{X}^u_l)\biggr](s)\\
=D\widehat{\theta}_u^{-1}(\widehat{X}^u_k)(s)+\frac12\sum\limits_{l=1}^Ns_l\sum\limits_{\operatorname{deg}X_m=\operatorname{deg}X_k+\operatorname{deg}X_l}c_{klm}(u)D\widehat{\theta}_u^{-1}(\widehat{X}^u_m)(s).
\end{multline*}
It follows
$$
D\widehat{\theta}_u(s)e_k=\widehat{X}^u_k(\theta_u(s))+\frac12\sum\limits_{l=1}^Ns_l\sum\limits_{\operatorname{deg}X_m=\operatorname{deg}X_k+\operatorname{deg}X_l}c_{klm}(u)\widehat{X}^u_m(\theta_u(s)).
$$
Applying further the arguments of the $1^{\text{st}}$ step, we
have
$$
\widehat{X}^u_j(\theta_u(s))=\sum\limits_{k=1}^N
\Bigl[z^k_j(s)+\frac12\sum\limits_{m,l:\operatorname{deg}X_k=\operatorname{deg}X_m+\operatorname{deg}X_l}z^m_j(s)s_lc_{mlk}(u)\Bigr]\widehat{X}^u_k(\theta_u(s)),
$$
and thus \eqref{100} is proved.

Taking into account the result of the $3^{\text{rd}}$ step, we
obtain
$$
\widehat{X}_j^u(w)=\sum\limits_{k:\,\operatorname{deg}X_k\leq\operatorname{deg}X_j}[\delta_{kj}+O(\varepsilon)]X_k
+\sum\limits_{k:\,\operatorname{deg}X_k>\operatorname{deg}X_j}
o(\varepsilon^{{\operatorname{deg}X_k-\operatorname{deg}X_j}})X_k(w),
$$
$j=1,\ldots, N$. The theorem follows.
\end{proof}

\begin{rem}\label{remcoef}
{\bf 1.} If the vector fields $X_i$, $i=1,\ldots, N$, belong to
the class $C^{1,\alpha}$, $\alpha\in(0,1]$, then in
{\eqref{cmlkinc}} and, consequently, in {\eqref{sumcmlk}}, we
obtain $o(1)=O(\rho(u,\theta_u(s))^{\alpha})$. In this case, we
have
\begin{multline*}
\widehat{X}_j^u(w)=\sum\limits_{k:\,\operatorname{deg}X_k\leq\operatorname{deg}X_j}[\delta_{kj}+O(\varepsilon)]X_k\\
+\sum\limits_{k:\,\operatorname{deg}X_k>\operatorname{deg}X_j}
\rho(u,\theta_u(s))^{\alpha}\cdot
o(\varepsilon^{{\operatorname{deg}X_k-\operatorname{deg}X_j}})X_k(w).
\end{multline*}

{\bf 2.} If the derivatives of the basis vector fields are
H\"{o}lder with respect to $d_{\infty}$, we obtain
$o(1)=O(d_{\infty}(u,\theta_u(s))^{\alpha})=O(\varepsilon^{\alpha})$,
and
$$
\widehat{X}_j^u(w)=\sum\limits_{k:\,\operatorname{deg}X_k\leq\operatorname{deg}X_j}[\delta_{kj}+O(\varepsilon)]X_k
+\sum\limits_{k:\,\operatorname{deg}X_k>\operatorname{deg}X_j}
O(\varepsilon^{{\operatorname{deg}X_k-\operatorname{deg}X_j}+\alpha})X_k(w).
$$

{\bf 3.} If the derivatives of the basis vector fields are
H\"{o}lder with respect to $\mathfrak d$, we have
\begin{multline*}
\widehat{X}_j^u(w)=\sum\limits_{k:\,\operatorname{deg}X_k\leq\operatorname{deg}X_j}[\delta_{kj}+O(\varepsilon)]X_k\\
+\sum\limits_{k:\,\operatorname{deg}X_k>\operatorname{deg}X_j}
\mathfrak d(u,\theta_u(s))^{\alpha}\cdot
o(\varepsilon^{{\operatorname{deg}X_k-\operatorname{deg}X_j}})X_k(w).
\end{multline*}
\end{rem}

\begin{cor}\label{lem_matrix}
For $x\in\operatorname{Box}(u,\varepsilon)$, the coefficients
$\{a_{j,k}(x)\}_{j,k=1}^N$ from the equality
\begin{equation}\label{expan}
X_j(x)=\sum\limits_{k=1}^Na_{j,k}(x)\widehat{X}^u_k(x)
\end{equation}
enjoy the following property$:$
\begin{equation}\label{estimate}
a_{j,k}(x)=
\begin{cases}
O(\varepsilon) &\text{ if }\ \operatorname{deg}X_j>\operatorname{deg}X_k,\\
\delta_{kj}+O(\varepsilon) &\text{ if }\ \operatorname{deg}X_j=\operatorname{deg}X_k,\\
o(\varepsilon^{\operatorname{deg}X_k-\operatorname{deg}X_j})
&\text{ if }\ \operatorname{deg}X_j<\operatorname{deg}X_k,
\end{cases}
\end{equation}
$j=1,\ldots, N$. All ``o'' are uniform  in $u$ belonging to some
compact  subset of  $U$.
\end{cor}

\begin{proof}
According to  Theorem~{\ref{th_expr}}, the coefficients
$b_{j,k}(x)$ from the relation
$$
\widehat{X}_j^u(x)=\sum\limits_{k=1}^Nb_{j,k}(x)X_k(x),
$$
$j=1,\ldots, N$,  have the same properties. Put
$A(x)=(a_{j,k}(x))_{j,k=1}^N$ and $B(x)=(b_{j,k}(x))_{j,k=1}^N$.
Then $A(x)=B(x)^{-1}$.

We use the well-known formula of calculation of the entries of the
inverse matrix to estimate all $a_{j,k}(x)$, $j,k=1\ldots, N$. We
estimate the value $|a_{j,k}(x)|=\frac{|\det B_{j,k}(x)|}{|\det
B(x)|}$, where the $(N-1)\times(N-1)$-matrix $B_{j,k}$ is
constructed from the matrix $B(x)$ by deleting its $j$th column
and $k$th line.

It is easy to see that $|\det B(x)|=1+O(\varepsilon)$, where
$O(\varepsilon)$ is uniform for $x$ belonging to some compact
neighborhood $U\subset\mathbb M$.

Next, we estimate $|\det B_{j,k}(x)|$. Obviously, $|\det
B_{j,j}(x)|=1+O(\varepsilon)$, where $O(\varepsilon)$ is uniform
for $x$ belonging to some compact neighborhood $U\subset\mathbb
M$, $j=1,\ldots, N$.

Let now $k>j$. By construction, the diagonal elements with numbers
$(i,i)$, $j\leq i<k$, equal
$o(\varepsilon^{\operatorname{deg}X_{i+1}-\operatorname{deg}X_i})$,
and the elements under these ones equal $1+O(\varepsilon)$. Note
that, $\det B_{j,k}(x)$ up to a multiple $(1+O(\varepsilon))$
equals the product of determinants of the following three
matrices: the first $P(x)=p_{i,l}(x)$ is a
$(j-1)\times(j-1)$-matrix with $p_{i,l}(x)=b_{i,l}(x)$, the second
$Q(x)=q_{i,l}(x)$ is a $(k-j)\times(k-j)$-matrix with
$q_{i,l}(x)=b_{i+j-1,l+j}(x)$, and the third $R(x)=r_{i,l}(x)$
with $r_{i,l}(x)=b_{i+k-1,l+k-1}(x)$.

For the matrices $P(x)$ and $R(x)$ we have $|\det
P(x)|=1+O(\varepsilon)$ and $|\det R(x)|=1+O(\varepsilon)$. By
construction,
$q_{i,i}(x)=o(\varepsilon^{\operatorname{deg}X_{i+1}-\operatorname{deg}X_i})$
and $q_{i+1,i}(x)=1+O(\varepsilon)$. We have that the product of
the diagonal elements of $Q(x)$ equals
$$
\prod\limits_{i=j}^{k-1}o(\varepsilon^{\operatorname{deg}X_{i+1}-\operatorname{deg}X_i})
=o(\varepsilon^{\operatorname{deg}X_{k}-\operatorname{deg}X_j}).
$$
It is easy to see that, for all other summands constituting $\det
Q(x)$, we have the same estimate.

Similarly, we show that for $k<j$ we have $|\det
B_{jk}(x)|=O(\varepsilon)$. Here $O(\varepsilon)$ is uniform for
$x$ belonging to some compact neighborhood $U\subset\mathbb M$.
The lemma follows.
\end{proof}

\begin{rem}\label{rem2coef}
Similarly to Remark {\ref{remcoef}}:
\begin{itemize}

\item if $X_i\in C^{1,\alpha}$ then
$$
a_{j,k}(x)=
\begin{cases}
O(\varepsilon) &\text{ if }\ \operatorname{deg}X_j>\operatorname{deg}X_k,\\
\delta_{kj}+O(\varepsilon) &\text{ if }\ \operatorname{deg}X_j=\operatorname{deg}X_k,\\
\rho(u,x)^{\alpha}\cdot
o(\varepsilon^{\operatorname{deg}X_k-\operatorname{deg}X_j})
&\text{ if }\ \operatorname{deg}X_j<\operatorname{deg}X_k,
\end{cases}
$$

\item if the derivatives of the basis vector fields are H\"{o}lder
with respect to $d_{\infty}$ then
$$
a_{j,k}(x)=
\begin{cases}
O(\varepsilon) &\text{ if }\ \operatorname{deg}X_j>\operatorname{deg}X_k,\\
\delta_{kj}+O(\varepsilon) &\text{ if }\ \operatorname{deg}X_j=\operatorname{deg}X_k,\\
O(\varepsilon^{\operatorname{deg}X_k-\operatorname{deg}X_j+\alpha})
&\text{ if }\ \operatorname{deg}X_j<\operatorname{deg}X_k,
\end{cases}
$$

\item if the derivatives of the basis vector fields are H\"{o}lder
with respect to $\mathfrak d$ then
$$
a_{j,k}(x)=
\begin{cases}
O(\varepsilon) &\text{ if }\ \operatorname{deg}X_j>\operatorname{deg}X_k,\\
\delta_{kj}+O(\varepsilon) &\text{ if }\ \operatorname{deg}X_j=\operatorname{deg}X_k,\\
\mathfrak d(u,x)^{\alpha}\cdot
O(\varepsilon^{\operatorname{deg}X_k-\operatorname{deg}X_j})
&\text{ if }\ \operatorname{deg}X_j<\operatorname{deg}X_k,
\end{cases}
$$
\end{itemize}
$j=1,\ldots, N$.
\end{rem}

Corollary {\ref{lem_matrix}} imply instantly Gromov's Theorem on
the nilpotentization of vector fields~\cite{G}.


\begin{cor}[Gromov's Theorem \cite{G}]\label{GTheorem}
Let $X_j\in C^{1}$. On  $\operatorname{Box}(g,\varepsilon r_g)$,
consider the vector fields
$\{^{\varepsilon}X_i\}{=}\{\varepsilon^{\deg X_i}X_i\}$,
$i=1,\dots,N$. Then the uniform convergence
$$
X_i^{\varepsilon}=\bigl(\Delta^g_{\varepsilon^{-1}}\bigr)_*{^{\varepsilon}X_i}
\to \widehat{X}^g_i\quad \text{as $\varepsilon\to 0$},\quad
i=1,\dots,N,
$$
holds at the points of the box~$\operatorname{Box}(g,r_g)$ and
this convergence is uniform in  $g$ belonging to some compact set.
\end{cor}

\begin{proof} Really, by \eqref{expan}, \eqref{estimate} and in view of
Corollary \ref{lem_matrix} and  Property~{\ref{proplocCargr}}, we
infer
\begin{multline*}
X_i^{\varepsilon}(x)=\bigl(\bigl(\Delta^g_{\varepsilon^{-1}}\bigr)_*{^{\varepsilon}X_i}\bigr)(x)
=\varepsilon^{\deg
X_i}\sum\limits_{k=1}^Na_{i,k}\bigl(\Delta^g_{\varepsilon}(x)\bigr)
\bigl(\bigl(\Delta^g_{\varepsilon^{-1}}\bigr)_*\widehat X^g_k\bigr)(x)\\
=\sum\limits_{k=1}^N\varepsilon^{\deg X_i-\deg X_k}a_{i,k}\bigl(\Delta^g_{\varepsilon}(x)\bigr)\widehat X^g_k(x)\\
= \sum\limits_{k:\,\operatorname{deg}X_k\leq
\operatorname{deg}X_i}\varepsilon^{\deg X_i-\deg
X_k}(\delta_{ik}+O(\varepsilon)) \widehat X^g_k(x)+
\sum\limits_{k:\,\operatorname{deg}X_k>\operatorname{deg}X_i}o(1)\widehat
X^g_k(x)
\end{multline*}
as $\varepsilon\to0$. It follows the uniform convergence
$X_i^{\varepsilon}=\bigl(\Delta^g_{\varepsilon^{-1}}\bigr)_*{^{\varepsilon}X_i}
\to \widehat{X}^g_i$ as $\varepsilon\to 0$, $i=1,\dots,N$, at the
points of the box~$\operatorname{Box}(g,r_g)$ and this convergence
is uniform in  $g$ belonging to some compact set.
\end{proof}

\begin{rem}\label{Ber}
For $C^\infty$-vector fields, the above corollary is formulated in
\cite{Me,rs} in another way: $\widehat{X}^g_i$ is an homogeneous
part of ${X}_i$, $1=1,\ldots,N$. This statement implies  Corollary
\ref{GTheorem}. It is shown in \cite{greshn1} that, applying
similar arguments, the smoothness of vector fields can be reduced
to be $2M+1$.

Estimates \eqref{estimate} were written in the proof of
\cite[Thereom 3.1]{V3} as a consequence of the Gromov's Theorem
which can be proved  by method of \cite{rs} under an additional
smoothness of vector fields: $X_j\in C^{2M-\deg X_j}$. Corollary
\ref{GTheorem} shows that estimates \eqref{estimate}  are not only
necessary but also sufficient for the validity of the  Gromov's
Theorem. In our paper estimates \eqref{estimate}  are obtained
under minimal assumption on the smoothness of vector fields:
$X_j\in C^{1}$,
 $j=1,\ldots,N$. Thus, taking into account
 the footnote in \cite[p.~253]{V3},  all results of papers \cite{vod1,vod2,vod3,V3,vk1}
 are valid under the same assumptions on the smoothness of basis vector fields.

Recall that Gromov \cite[p.~130]{G} has formulated the theorem
under assumption  $X_j\in C^{1}$.  Valeri\v{\i} Berestovski\v{\i}
sent us the following example confirming that arguments of
Gromov's proof have to be corrected.

\begin{ex}\label{cex}
Let $X = \frac{\partial}{\partial x}$, $Y = xy
\frac{\partial}{\partial x} + \frac{\partial}{\partial y} + x
\frac{\partial}{\partial z}$. Then $Z := [X, Y ] =
y\frac{\partial}{\partial x} + \frac{\partial}{\partial z}$,
$[X,Z] = 0$, $[Y,Z] = \frac{\partial}{\partial x} - y\bigl(y
\frac{\partial}{\partial x} + \frac{\partial}{\partial z} \bigr) =
(1 - y^2) \frac{\partial}{\partial x} - y \frac{\partial}{\partial
z}$. One can easily see that $X, Y, Z$ constitutes a global frame
of smooth vector fields over the ring of smooth functions in
$\mathbb R^3$. Also for corresponding one-parameter subgroups
$X(x)$, $Y (y)$, $Z(z)$, we have $(X(x)\circ Y (y)\circ Z(z))(0,
0, 0) = (x, y, z)$. Under this $X = \frac{\partial}{\partial x}$
on $\mathbb R^3$, $Y = \frac{\partial}{\partial y}$ on $x = 0$, $Z
= \frac{\partial}{\partial z}$ on $z$-line (even on $y = 0$). On
the other hand, $\frac{\partial}{\partial y}Z = X \ne [Y,Z]$ (see
above) on $x = 0$. This contradicts to the Gromov's statement that
(A) of \cite[p.~131]{G} implies (B) of  \cite[p.~132]{G} in
general case.
\end{ex}
\end{rem}

\begin{cor}[Estimate of the Diameter of a Box]\label{estdiam}
In a compact neighborhood $U\subset\mathbb M$, for each point
$u\in U$ and each $\varepsilon>0$ small enough, we have
$\operatorname{diam}(\operatorname{Box}(u,\varepsilon))\leq
L\varepsilon$, where $L$ depends only on $U$.
\end{cor}

\begin{proof}
Assume the contrary: there exist sequences
$\{\varepsilon_k\}_{k\in\mathbb N}$,
 $\{u_k\}_{k\in\mathbb
N}$, $\{v_k\}_{k\in\mathbb N}$   and $\{w_k\}_{k\in\mathbb N}$
  such that
$\varepsilon_k\to0$ as $k\to\infty$,
$d_{\infty}(u_k,v_k)=\varepsilon_k$ and
$d_{\infty}(u_k,w_k)\leq\varepsilon_k$ but
$d_{\infty}(v_k,w_k)>k\varepsilon_k$. Since $U\subset\mathbb M$ is
compact, we may assume without loss of generality that $u_k\to
u_0$ as $k\to\infty$. Then $v_k\to u_0$ and $w_k\to u_0$ as
$k\to\infty$.

Assume without loss of generality that
$\varepsilon^{\operatorname{deg}X_i}D\Delta^{u_k}_{\varepsilon^{-1}}X_i(x)\to\widehat{X}^{u_k}_i(x)$
as $\varepsilon\to0$ for $x\in\operatorname{Box}(u_0,Kr_0)$
uniformly in $u_k$, $i=1,\ldots, N$, where $K=\max\{5, 5c^4\}$,
$c$ is such that $d_{\infty}^{u_k}(v,w)\leq
c(d_{\infty}^{u_k}(u,v)+d_{\infty}^{u_k}(u,w))$ for all
$k\in\mathbb N$ big enough, and $k\in\mathbb N$ is big enough (see
Corollary~{\ref{GTheorem}}). Note that, $c<\infty$ since $c(u_k)$
continuously depends on values of
$\{F^j_{\mu,\beta}(u_k)\}_{j,\mu,\beta}$, consequently, it depends
continuously on~$u_k$. Moreover, the choice of $K$ implies the
following:
\begin{enumerate}
\item For $k$ big enough, we have that an integral curve of a
vector field with constant coefficients connecting
$\Delta^{u_k}_{r_0\varepsilon_k^{-1}}(w_k)$ and
$\Delta^{u_k}_{r_0\varepsilon_k^{-1}}(v_k)$ in the local Carnot
group $\mathcal G^{u_k}\mathbb M$ lies in $\operatorname{Box}(u_k,
Kr_0)$;

\item We may choose $k$ by the following way:
$d_{\infty}(u_0,u_k)<r_0$ and the Riemannian distance between the
integral curves corresponding to the collections
$\{\widehat{X}^{u_k}_i\}_{i=1}^N$ and
$\{({r_0}^{-1}\varepsilon_k)^{\operatorname{deg}X_i}D\Delta^{u_k}_{r_0\varepsilon_k^{-1}}X_i\}_{i=1}^N$
(with constant coefficients) that connect points
$\Delta^{u_k}_{r_0\varepsilon_k^{-1}}(w_k)$ and
$\Delta^{u_k}_{r_0\varepsilon_k^{-1}}(v_k)$, is less than $r_0$.
\end{enumerate}
Fix $k\in\mathbb N$. Then
$v_k=\exp\Bigl(\sum\limits_{i=1}^N\xi_i\varepsilon_k^{\operatorname{deg}X_i}X_i\Bigr)(u_k)$,
$w_k=\exp\Bigl(\sum\limits_{i=1}^N\eta_i\varepsilon_k^{\operatorname{deg}X_i}X_i\Bigr)(u_k)$,
and
$w_k=\exp\Bigl(\sum\limits_{i=1}^N\zeta_i(\varepsilon_k)\varepsilon_k^{\operatorname{deg}X_i}X_i\Bigr)(v_k)$.
Apply the mapping $\Delta^{u_k}_{r_0\varepsilon_k^{-1}}$ to $v_k$
and $w_k$. We have
$$
\Delta^{u_k}_{r_0\varepsilon_k^{-1}}(w_k)=
\exp\biggl(\sum\limits_{i=1}^N\zeta_i(\varepsilon)
\varepsilon_k^{\operatorname{deg}X_i}D\Delta^{u_k}_{r_0\varepsilon_k^{-1}}X_i\biggr)
\bigl(\Delta^{u_k}_{r_0\varepsilon_k^{-1}}(v_k)\bigr).
$$
Note that,
$d_{\infty}\bigl(u_k,\Delta^{u_k}_{r_0\varepsilon_k^{-1}}(v_k)\bigr)=r_0$
and
$d_{\infty}\bigl(u_k,\Delta^{u_k}_{r_0\varepsilon_k^{-1}}(w_k)\bigr)\leq
r_0$. In view of Corollary~{\ref{GTheorem}}, the vector fields
$({r_0}^{-1}\varepsilon_k)^{\operatorname{deg}X_i}D\Delta^{u_k}_{r_0\varepsilon_k^{-1}}X_i(x)=
\widehat{X}_i^{u_k}(x)+o(1)$, $i=1,\ldots, N$, where $o(1)$ is
uniform in $x$ and in $u_k$. Consequently, since
\linebreak$\dim\operatorname{span}\{\widehat{X}^{u_k}_i(x)\}_{i=1}^N=N$
at each $x\in\operatorname{Box}(u_0,r_0)$, the Riemannian distance
between $\Delta^{u_k}_{r_0\varepsilon_k^{-1}}(w_k)$ and
$\Delta^{u_k}_{r_0\varepsilon_k^{-1}}(v_k)$ is bounded from above
for all $k\in\mathbb N$ big enough. Therefore, the coefficients
$\zeta_i(\varepsilon_k)$, $i=1,\ldots, N$, are bounded from above
for all $k\in\mathbb N$ big enough. The assumption
$d_{\infty}(v_k,w_k)>k\varepsilon_k$ contradicts this conclusion.

Thus there exists a constant $L=L(U)$ such that
$\operatorname{diam}(\operatorname{Box}(u,\varepsilon))\leq
L\varepsilon$ for $u\in U$. The statement follows.
\end{proof}

From the previous statement we come immediately to the following

\begin{cor}[Triangle inequality]
\label{trin}
 The quasimetric $d_\infty(x,y)$
meets locally the generalized triangle inequality $($see
Property~\ref{propmetr}$)$.
\end{cor}

\begin{cor}[Decomposition of the basis vector fields]
Fix a point $\theta_u(s)\in\operatorname{Box}(u,O(\varepsilon))$.
Remarks~{\ref{remcoef}} and~{\ref{rem2coef}} imply the following
decomposition of $D\theta_u^{-1}X_i$, $i=1,\ldots, N$:
$$
[D\theta_u^{-1}X_i(s)]_j=[(\widehat{X}_i^u)^{\prime}(s)]_j+\sum\limits_{k=1}^N
(a_{i,k}(\theta_u(s))-\delta_{ik})[(\widehat{X}_k^u)^{\prime}(s)]_j.
$$
If $d_{\infty}(u,\theta_u(s))=O(\varepsilon)$, we have
\begin{multline*}
[D\theta_u^{-1}X_i(s)]_j=z^j_i(u,s)+\sum\limits_{k:\,\operatorname{deg}X_k\leq\operatorname{deg}X_i}
O(\varepsilon)z^j_k(u,s)\\
+\sum\limits_{k:\,\operatorname{deg}X_k>\operatorname{deg}X_i}
a_{i,k}(\theta_u(s))z^j_k(u,s).
\end{multline*}
If $\operatorname{deg}X_j\leq\operatorname{deg}X_i$ then
$[D\theta_u^{-1}X_i(s)]_j=\delta_{ij}+O(\varepsilon)$. For
$\operatorname{deg}X_j>\operatorname{deg}X_i$ we have:
\begin{itemize}
\item If the basis vector fields are $C^1$-smooth then we deduce
$[D\theta_u^{-1}X_i(s)]_j=z^j_i(u,s)+O(\varepsilon^{\operatorname{deg}X_j-\operatorname{deg}X_i+1})
+o(1)\cdot
\varepsilon^{\operatorname{deg}X_j-\operatorname{deg}X_i}$, and
therefore
$$
[D\theta_u^{-1}X_i(s)]_j=z^j_i(u,s)+o(\varepsilon^{\operatorname{deg}X_j-\operatorname{deg}X_i}).
$$
\item If the derivatives of the basis vector fields are
$H^{\alpha}$-continuous with respect to $\mathfrak d$, then if
$\operatorname{deg}X_j>\operatorname{deg}X_i$ we have
$$
[D\theta_u^{-1}X_i(s)]_j=z^j_i(u,s)+\mathfrak
d(u,\theta_u(s))^{\alpha}\cdot
O(\varepsilon^{\operatorname{deg}X_j-\operatorname{deg}X_i}).
$$
\end{itemize}
In particular, for $\alpha=1$ and $\mathfrak d=d_{\infty}$ or
$\mathfrak d=d^z_{\infty}$, where
$d_{\infty}(u,z)=O(\varepsilon)$, we have
$$
[D\theta_u^{-1}X_i(s)]_j=z^j_i(u,s)+
O(\varepsilon^{\operatorname{deg}X_j-\operatorname{deg}X_i+1}).
$$
\end{cor}

\subsection{Comparison of Geometries of Tangent Cones}\label{geomconeconeglob}

The goal of Subsections {\ref{geomconeconeglob}},
{\ref{geomconecone1}} and {\ref{geomconecone}} is to compare the
geometries of two local Carnot groups. The main result of
Section~{\ref{geomcarnot}} is the following

\begin{thm}
\label{est_chain} Let $u,u^{\prime}\in U$ be such that
$d_{\infty}(u,u^{\prime})=C\varepsilon$. For a fixed $Q\in\mathbb
N$, consider points $w_0$, $d_{\infty}(u,w_0)={\mathcal
C}\varepsilon$, and
$$
w_j^{\varepsilon}=\exp\biggl(\sum\limits_{i=1}^{N}w_{i,j}{\varepsilon}^{\operatorname{deg}
X_i}\widehat{X}_i^u\biggr)(w^{\varepsilon}_{j-1}),\quad
{w_j^{\varepsilon}}^{\prime}=\exp\biggl(\sum\limits_{i=1}^{N}w_{i,j}{\varepsilon}^{\operatorname{deg}
X_i}\widehat{X}_i^{u^{\prime}}\biggr)({w^{\varepsilon\prime}_{j-1}}),
$$
$w^{\varepsilon\prime}_0=w^{\varepsilon}_0=w_0^{\prime}=w_0$,
$j=1,\ldots,Q$. {$($Here $Q\in\mathbb N$ is such that all these
points belong to the neighborhood $U\subset\mathbb M$, for all
$\varepsilon>0$.$)$} Then for $\alpha>0$,
\begin{equation}\label{est_l}
\max\{d_{\infty}^u(w^{\varepsilon}_Q,w^{\varepsilon\prime}_Q),
d_{\infty}^{u^{\prime}}(w^{\varepsilon}_Q,w^{\varepsilon\prime}_Q)\}=\varepsilon
\cdot [\Theta(C,{\mathcal
C},Q,\{F^j_{\alpha,\beta}\}_{j,\alpha,\beta})]\rho(u,u^{\prime})^{\frac{\alpha}{M}}.
\end{equation}
In the case of $\alpha=0$, we have
$$
\max\{d_{\infty}^u(w^{\varepsilon}_Q,w^{\varepsilon\prime}_Q),
d_{\infty}^{u^{\prime}}(w^{\varepsilon}_Q,w^{\varepsilon\prime}_Q)\}=\varepsilon\cdot
[\Theta(C,{\mathcal
C},Q,\{F^j_{\alpha,\beta}\}_{j,\alpha,\beta})][\omega(\rho(u,u^{\prime}))]^{\frac{1}{M}}
$$ where $\omega\to0$ is a modulus of continuity. $($Here ${\Theta}$ is
uniform in $u, u^{\prime}, w_0\in U$ and $\{w_{i,j}\}$,
$i=1,\ldots, N$, $j=1,\ldots,Q$, belonging to some compact
neighborhood of $0$, { and it depends on $Q$ and
$\{F^j_{\mu,\beta}\}_{j,\mu,\beta}$.}$)$
\end{thm}

\begin{rem}
If the derivatives of $X_i$, $i=1,\ldots, N$, are locally
H\"{o}lder with respect to $\mathfrak d$, then we have $\mathfrak
d(u,u^{\prime})^{\frac{\alpha}{M}}$ instead of
$\rho(u,u^{\prime})^{\frac{\alpha}{M}}$ in {\eqref{est_l}}.
\end{rem}

In the current subsection we prove the ''base`` of the main
result, i.~e., we obtain it for $Q=1$ and $\varepsilon=1$. The
full proof is written in Subsection~{\ref{geomconecone}}.

Fix points $u,u^{\prime}\in U$, where $U$ is such that Assumption
{\ref{vectorfields}} holds. Recall that the collections of vector
fields $\{\widehat{X}_i^u\}_{i=1}^{N}$ and
$\{\widehat{X}_i^{u^{\prime}}\}_{i=1}^{N}$ are frames in
${\mathcal G}^u\mathbb M$ and in ${\mathcal G}^{u^{\prime}}\mathbb
M$ respectively.

\begin{defn}
By $\widehat{X}^p(q)$, we denote the matrix, such that its $i$th
column consists of the coordinates of the vector
$\widehat{X}^p_i(q)$, $i=1,\ldots, N$, $p\in\mathbb M$,
$q\in{\mathcal G}^p\mathbb M$, in the frame $\{\widehat
X_j\}_{j=1}^N$.
\end{defn}

\begin{lem}
\label{smoothxi} Suppose that Assumption {\ref{vectorfields}}
holds. Let $\Xi(u,u^{\prime},q)$, $q\in\mathbb M$, be the matrix
such that
\begin{equation}\label{trans_mat}
\widehat{X}^{u^{\prime}}(q)=\widehat{X}^{u}(q)\Xi(u,u^{\prime},q).
\end{equation}
Then the entries of $\Xi(u,u^{\prime},q)$ are $($locally$)$
$H^{\alpha}$-continuous in $u$ and $u^{\prime}$.
\end{lem}

\begin{proof}
The proof of this statement follows from Theorem {\ref{tcone}}.
Indeed, it implies that the vector fields
$\{\widehat{X}_i^u\}_{i=1}^N$ are locally $H^{\alpha}$-continuous
in $u$. Since we prove a local property, and $\mathbb M$ is a
Riemannian manifold, then, instead of $\mathbb M$, we may consider
without loss of generality some neighborhood $U\subset\mathbb R^N$
containing $u$ and $u^{\prime}$. Then it is easy to see that the
entries of the matrices $\widehat{X}^u$ and
$\widehat{X}^{u^{\prime}}$ are (locally) $H^{\alpha}$-continuous
on $U\times U$. Since both matrices are non-degenerate in
 $U\subset\mathbb M$, we have that
$\Xi(u,u^{\prime},q)=\widehat{X}^{u}(q)^{-1}\widehat{X}^{u^{\prime}}(q)$
is also non-degenerate, and its entries $\Xi_{ij}(u,u^{\prime},q)$
belong locally to $C^{\alpha}(U\times U)$, $i,j=1,\ldots,N$.
\end{proof}

\begin{rem}
If the derivatives of $X_i$, $i=1,\ldots, N$, are locally
H\"{o}lder with respect to $\mathfrak d$, then the entries of
$\Xi$ are also locally H\"{o}lder with respect to $\mathfrak d$
(see Remark~{\ref{nilphold}}).
\end{rem}

\begin{rem}
\label{omega} Suppose that Assumption {\ref{vectorfields}} holds.
Since $\Xi(u,u^{\prime},q)$ equals the unit matrix if
$u=u^{\prime}$ then $\Xi_{ij}=\delta_{ij}+\Theta
\rho(u,u^{\prime})^{\alpha}$ where $\Theta=\Theta(u,u^{\prime},q)$
is a bounded measurable function: $|\Theta|\leq C$, and the
constant $C\geq 0$ depends only on the neighborhood
$U\subset\mathbb M$.
\end{rem}

\begin{proof}
Note that $\Xi(u,u,q)$ equals the unit matrix. Then the
$\alpha$-H\"{o}lder continuity of all vector fields implies
$|\Xi_{ij}(u,u^{\prime},q)-\delta_{ij}|\leq
C(\rho(u,u^{\prime})^{\alpha})$, where
$$
C=\sup\limits_{u,u^{\prime},q\in
U}\frac{|\Xi_{ij}(u,u^{\prime},q)-\delta_{ij}|}{\rho(u,u^{\prime})^{\alpha}}<\infty
$$
depends only on the neighborhood $U\subset\mathbb M$.
\end{proof}

\begin{rem}\label{xihold}
If the derivatives of $X_i$, $i=1,\ldots, N$, are locally
H\"{o}lder with respect to $\mathfrak d$, then
$\Xi_{ij}(u,u^{\prime},q)=\delta_{ij}+\Theta\mathfrak
d(u,u^{\prime})^{\alpha}$.
\end{rem}

\begin{notat}
\label{func_theta} Throughout the paper, by the symbol $\Theta$,
we denote some bounded function absolute values of which do not
exceed some $0\leq C<\infty$, where $C$ depends only on the
neighborhood where $\Theta$ is defined (i.~e., it does not depend
on {\it points} of this neighborhood).
\end{notat}

\begin{thm}\label{difference}  Let
$$
w=\exp\Bigl(\sum\limits_{i=1}^{N}w_i\widehat{X}_i^u\Bigr)(v),\quad
w^{\prime}=\exp\Bigl(\sum\limits_{i=1}^{N}w_i\widehat{X}_i^{u^{\prime}}\Bigr)(v).
$$
Then, for $\alpha>0$, we have
\begin{equation}\label{eqdifference}
\max\{d_{\infty}^u(w,w^{\prime}),
d_{\infty}^{u^{\prime}}(w,w^{\prime})\}=\Theta[\rho(u,u^{\prime})^{\alpha}
\rho(v,w)]^{\frac{1}{M}},
\end{equation}
where $u,u^{\prime},v\in U$, $\{w_i\}_{i=1}^N\in
U(0)\subset\mathbb R^N$.

In the case of $\alpha=0$,
$$
\max\{d_{\infty}^u(w,w^{\prime}),
d_{\infty}^{u^{\prime}}(w,w^{\prime})\}=\Theta[\omega_{\Xi}(\rho(u,u^{\prime}))
\rho(v,w)]^{\frac{1}{M}}.
$$
\end{thm}

\begin{rem}
Here (see Notation {\ref{func_theta}}), the value
$$
\sup|\Theta(u,u^{\prime},v,\{w_i\}_{i=1}^N)|<\infty
$$
depends only on $U\subset\mathbb M$ and $U(0)\subset\mathbb R^N$.
\end{rem}

\begin{proof}[Proof of Theorem {\ref{difference}}]
$1^{\text{\sc st}}$ {\sc Step.} Fix $q\in\mathbb M$. Notice that
both collections of vectors $\{\widehat{X}^u_i(q)\}_{i=1}^{N}$ and
$\{\widehat{X}^{u^{\prime}}_i(q)\}_{i=1}^{N}$ are frames of
$T_q\mathbb M$. Consequently, there exists the transition
$(N\times N)$-matrix
$\Xi(u,u^{\prime},q)=(\Xi(u,u^{\prime},q))_{i,k}$ such that
\begin{equation}\label{xiik}
\widehat{X}^{u^{\prime}}_i(q)=\sum\limits_{k=1}^{N}(\Xi(u,u^{\prime},q))_{i,k}\widehat{X}^u_k(q).
\end{equation}
Remark {\ref{omega}} implies that
\begin{equation}\label{xiest}
\Xi(u,u^{\prime},q)_{i,j}=
\begin{cases}1+{\Theta_{i,j}} \rho(u,u^{\prime})^{\alpha} &\text{ if }i=j,\\
{\Theta_{i,j}}\rho(u,u^{\prime})^{\alpha} &\text{ if }i\neq j.
\end{cases}
\end{equation}
Thus
$\widehat{X}^{u^{\prime}}_i(q)=\widehat{X}^u_i(q)+[\Xi(u,u^{\prime},q)-I]\widehat{X}^u_i(q)$
where $|[\Xi(u,u^{\prime},q)-I]|_{k,j}=\Theta_{k,j}
\rho(u,u^{\prime})^{\alpha}$ for all $k,j=1,\ldots, N$.

$2^{\text{\sc nd}}$ {\sc Step.} Consider the integral line
$\gamma(t)$ of the vector field
$\sum\limits_{i=1}^{N}w_i\widehat{X}_i^{u^{\prime}}$ starting at
$v$ with the endpoint $w^{\prime}$. Rewrite the tangent vector to
$\gamma(t)$ in the frame $\{\widehat{X}_i^u\}_{i=1}^{N_1}$ as
$\dot{\gamma}(t)=\sum\limits_{i=1}^{N}w_i^u(\gamma(t))\widehat{X}_i^u(\gamma(t))$.
From {\eqref{xiik}} it follows that
$$
w_i^u(q)=\sum\limits_{k=1}^Nw_k(\Xi(u,u^{\prime},q))_{k,i}.
$$
From {\eqref{xiest}} we can estimate the coefficient $w_i^u$ at
$\widehat{X}_i^u$:
\begin{equation}\label{w_i^u}
w_i^u=w_i+\sum\limits_{k=1}^{N}[w_k{\Theta}_{k,i}
\rho(u,u^{\prime})^{\alpha}],\quad i=1,\ldots,N.
\end{equation}

$3^{\text{\sc rd}}$ {\sc Step.} Next, we estimate the Riemannian
distance between $w$ and $w^{\prime}$. By $\kappa(t)$ denote the
integral line of the vector field
$\sum\limits_{i=1}^{N}w_i\widehat{X}_i^u$ connecting $v$ and $w$,
i.~e., a line such that $\kappa(0)=v$ and
$$
\dot{\kappa}(t)=\sum\limits_{i=1}^{N}w_i\widehat{X}_i^u(\kappa(t)).
$$

By means of the mapping $\theta_u^{-1}$ we transport $\kappa(t)$
and $\gamma(t)$ to $\mathbb R^N$. Let
$\kappa_u(t)=\theta_u^{-1}(\kappa(t))$ and
$\gamma_u(t)=\theta_u^{-1}(\gamma(t))$. Then
$$
\dot{\kappa}_u(t)=(\theta_u^{-1})_*(\kappa(t))\dot{\kappa}(t)
=\sum\limits_{i=1}^Nw_i(\widehat{X}^u_i)^{\prime}(\kappa_u(t))
$$
and similarly
$$
\dot{\gamma}_u(t)=\sum\limits_{i=1}^{N}w_i(\theta_u^{-1})_*\widehat{X}_i^{u^{\prime}}
=\sum\limits_{i=1}^Nw_i^u(t)(\widehat{X}^u_i)^{\prime}(\gamma_u(t))
$$
since
$(\theta_u^{-1})_*\widehat{X}^{u^{\prime}}_i(q)=\sum\limits_{k=1}^{N}(\Xi(u,u^{\prime},q))_{i,k}(\widehat{X}^u_i)^{\prime}(q)$
(see {\eqref{trans_mat}}). Using formula {\eqref{mat_tc}} rewrite
the tangent vectors in Cartesian coordinates:
$$
\dot{\kappa_u}(t)=\sum\limits_{i=1}^Nw_i\sum\limits_{j=1}^Nz^j_i(u,\kappa_u(t))\frac{\partial}{\partial
x_j}=\sum\limits_{j=1}^NW_j(u,\kappa_u(t))\frac{\partial}{\partial
x_j}
$$
where
$$
W_j(u,\kappa_u(t))=\sum\limits_{i=1}^Nw_iz^j_i(u,\kappa_u(t))
=w_j+\sum\limits_{i=1}^{j-1}w_iz^j_i(u,\kappa_u(t)).
$$
Similarly
$$
\dot{\gamma_u}(t)=\sum\limits_{j=1}^NW_j^u(u,\gamma_u(t))\frac{\partial}{\partial
x_j}
$$
where
$$
W_j(u,\gamma_u(t))=w_j^u(t)+\sum\limits_{i=1}^{j-1}w_i^u(t)z^j_i(u,\gamma_u(t)).
$$
Now we estimate the length of the curve
$\lambda_u(t)=\gamma_u(t)-\kappa_u(t)+\theta_u^{-1}(w)$ with
endpoints $\theta_u^{-1}(w)$ and $\theta_u^{-1}(w^{\prime})$. The
tangent vector to $\lambda_u(t)$ equals
\begin{multline}\label{betau}
\dot{\lambda_u}(t)=\dot{\gamma_u}(t)-\dot{\kappa_u}(t)
=\sum\limits_{j=1}^N[W_j^u(u,\gamma_u(t))-W_j(u,\kappa_u(t))]\frac{\partial}{\partial
x_j}\\
=\sum\limits_{j=1}^N\Bigl[(w_j^u(t)-w_j)
+\sum\limits_{i<j}w_i(z^j_i(u,\gamma_u(t))-z^j_i(u,\kappa_u(t)))\Bigr]\\
+\sum\limits_{i,j=1}^N(w_i^u(t)-w_i)z^j_i(u,\gamma_u(t)).
\end{multline}
Notice that for the last sum we have
$$
\sum\limits_{i,j=1}^N(w_i^u(t)-w_i)z^j_i(u,\gamma_u(t))={\Theta}\rho(u,u^{\prime})^{\alpha}
\rho(v,w)
$$
since $w_i^u(t)=w_i+{\Theta}\rho(u,u^{\prime})^{\alpha} \rho(v,w)$
by {\eqref{w_i^u}}. By properties of $z^j_i$,
$$
z^j_i(u,\gamma_u(t))-z^j_i(u,\kappa_u(t))
={\Theta}\Bigl[\sum\limits_{|\mu|=1}{F}^j_{\mu,e_i}(\gamma_u^{\mu}(t)-\kappa_u^{\mu}(t))\Bigr].
$$
Notice that
$$
|\gamma_u(t)-\kappa_u(t)|\leq
\int\limits_{0}^{t}|\dot{\gamma_u}(\tau)-\dot{\kappa_u}(\tau)|\,d\tau.
$$
Consequently
$$
\max\limits_t|{\gamma_u}(t)-{\kappa_u}(t)|
\leq\max\limits_t|\dot{\gamma_u}(t)-\dot{\kappa_u}(t)|=\max\limits_t|\dot{\lambda_u}(t)|.
$$
Applying these estimates to {\eqref{betau}} we obtain
$$
\max\limits_t|\dot{\lambda_u}(t)|={\Theta}\rho(u,u^{\prime})^{\alpha}
\rho(v,w)+{\Theta}\rho(v,w)\max\limits_t|\dot{\lambda_u}(t)|.
$$
From here it follows
$$
\max\limits_t|\dot{\lambda_u}(t)|=\frac{{\Theta}\rho(u,u^{\prime})^{\alpha}
\rho(v,w)}{1-{\Theta}\rho(v,w)}\leq{\Theta}\rho(u,u^{\prime})^{\alpha}
\rho(v,w)
$$
if ${\Theta}\rho(v,w)\leq\frac{1}{2}$. Thus
$$
\rho(\theta_u^{-1}(w),
\theta_u^{-1}(w^{\prime}))\leq\int\limits_0^1|\dot{\lambda_u}(t)|\,dt
\leq\max\limits_t|\dot{\lambda_u}(t)|={\Theta}\rho(u,u^{\prime})^{\alpha}
\rho(v,w),
$$
and $\rho(w,w^{\prime})\leq {\Theta}\rho(u,u^{\prime})^{\alpha}
\rho(v,w)$.

$4^{\text{\sc th}}$ {\sc Step.} By the inequality
$d_{\infty}^u(p,q)\leq C\rho(p,q)^{\frac{1}{M}}$, we obtain the
estimate of $d_{\infty}^u(w,w^{\prime})$:
$$
d_{\infty}^u(w,w^{\prime})={\Theta}[\rho(u,u^{\prime})^{\alpha}
 \rho(v,w)]^{\frac{1}{M}}
$$
in some compact neighborhood of $g$. The same estimate is true for
$d_{\infty}^{u^{\prime}}(w,w^{\prime})$. The theorem follows.
\end{proof}

\begin{rem}\label{geomhold}
If the derivatives of $X_i$, $i=1,\ldots, N$, are locally
H\"{o}lder with respect to $\mathfrak d$, then we have $\mathfrak
d(u,u^{\prime})^{\alpha}$ instead of $\rho(u,u^{\prime})^{\alpha}$
in {\eqref{eqdifference}} (the proof is similar, see
Remarks~{\ref{nilphold}} and~{\ref{xihold}}).
\end{rem}

\subsection{Comparison of Local Geometries of Tangent Cones}\label{geomconecone1}

Consider points
$$
w_{\varepsilon}=\exp\Bigl(\sum\limits_{i=1}^Nw_i\varepsilon^{\operatorname{deg}
X_i}\widehat{X}^{u}_i\Bigr)(v) \text{ and }
w_{\varepsilon}^{\prime}=\exp\Bigl(\sum\limits_{i=1}^Nw_i\varepsilon^{\operatorname{deg}
X_i}\widehat{X}^{u^{\prime}}_i\Bigr)(v).
$$

\begin{thm}
\label{est_eps} Suppose that
$d_{\infty}(u,u^{\prime})=C\varepsilon$ and $d_{\infty}(u,
v)={\mathcal C}\varepsilon$ for some $C,{\mathcal C}<\infty$.
Then, for $\alpha>0$, we have
\begin{equation}\label{locgeomest}
\max\{d_{\infty}^u(w_{\varepsilon},w^{\prime}_{\varepsilon}),
d_{\infty}^{u^{\prime}}(w_{\varepsilon},w^{\prime}_{\varepsilon})\}=\varepsilon
[\Theta(C,{\mathcal C})]\rho(u,u^{\prime})^{\frac{\alpha}{M}}.
\end{equation}
In the case of $\alpha=0$, we have
\begin{multline*}
\max\{d_{\infty}^u(w_{\varepsilon},w^{\prime}_{\varepsilon}),
d_{\infty}^{u^{\prime}}(w_{\varepsilon},w^{\prime}_{\varepsilon})\}\\
=\varepsilon [\Theta(C,{\mathcal C})]
\max\{\omega_{\Xi}(\rho(u,u^{\prime})),
\omega_{\Delta^u_{\varepsilon^{-1},v}\circ\Delta^{u^{\prime}}_{\varepsilon,v}}(\rho(u,u^{\prime}))\}^{\frac{\alpha}{M}},
\end{multline*}
where $\Delta^u_{\varepsilon^{-1},v}$ is defined below in
{\eqref{delta_epsilon}} and {\eqref{d_e_exact}}. $($Here
${\Theta}$ is uniform in $u, u^{\prime}, v\in U\subset\mathbb M$,
and in $\{w_i\}_{i=1}^N$ belonging to some compact neighborhood of
$0$ $($see Notation {\ref{func_theta}}$)$.$)$
\end{thm}

\begin{rem}
If the derivatives of $X_i$, $i=1,\ldots, N$, are locally
$\alpha$-H\"{o}lder with respect to $\mathfrak d$ (instead of
$\rho$), then we have $\mathfrak
d(u,u^{\prime})^{\frac{\alpha}{M}}$ instead of
$\rho(u,u^{\prime})^{\frac{\alpha}{M}}$ in {\eqref{locgeomest}}
(the proof is similar, see Remark~{\ref{geomhold}}).
\end{rem}

\begin{proof}[Proof of Theorem {\ref{est_eps}}]
$1^{\text{\sc st}}$ {\sc Step.} Let $w=w_1$ and
$w^{\prime}=w^{\prime}_1$ as it was earlier. In the frame
$\{\widehat{X}^u_i\}_{i=1}^N$ we have
$$
w^{\prime}=\exp\Bigl(\sum\limits_{i=1}^Nw_i^{\prime}\widehat{X}^{u}_i\Bigr)(v).
$$
Consider the point
$$
\omega_{\varepsilon}=\exp\Bigl(\sum\limits_{i=1}^Nw_i^{\prime}\varepsilon^{\operatorname{deg}
X_i}\widehat{X}^{u}_i\Bigr)(v).
$$
Note that $\omega_1=w^{\prime}$. In view of the generalized
triangle inequality,
$d_{\infty}^u(w_{\varepsilon},w^{\prime}_{\varepsilon})\leq
c(d_{\infty}^u(w_{\varepsilon},\omega_{\varepsilon})
+d_{\infty}^u(\omega_{\varepsilon},w^{\prime}_{\varepsilon}))$. By
the above estimate
\begin{equation}\label{est1}
d_{\infty}^u(\omega_{\varepsilon},w_{\varepsilon})=\varepsilon
d_{\infty}^u(w,w^{\prime})=\varepsilon
{\Theta}(\rho(u,u^{\prime})^{\alpha}
d_{\infty}^u(v,w))^{\frac{1}{M}}.
\end{equation}
Note that, if $\alpha=0$, then we obtain here
$\omega_{\Xi}(\rho(u,u^{\prime}))$.

Now we estimate the distance
$d_{\infty}^u(\omega_{\varepsilon},w^{\prime}_{\varepsilon})$.
Represent $w^{\prime}_{\varepsilon}$ in the frame
$\{\widehat{X}^u_i\}_{i=1}^N$:
\begin{equation}\label{est_w}
w^{\prime}_{\varepsilon}=\exp\Bigl(\sum\limits_{i=1}^N\alpha_i(\varepsilon)\varepsilon^{\operatorname{deg}
X_i}\widehat{X}^{u}_i\Bigr)(v),
\end{equation}
and consider the point
$$
\omega^{\prime}=\exp\Bigl(\sum\limits_{i=1}^N\alpha_i(\varepsilon)\widehat{X}^{u}_i\Bigr)(v).
$$
Here the coefficients $\alpha_i(\varepsilon)$, $i=1,\ldots, N$,
depend on $u$ and $\{w_i\}_{i=1}^N$.

$2^{\text{\sc nd}}$ {\sc Step.} Next, we show that the
coefficients $\alpha_i(\varepsilon)$, $i=1,\ldots, N$, are
uniformly bounded for all $\varepsilon>0$ uniformly on $u$ and
$\{w_i\}_{i=1}^N$. By another words, there exists $S<\infty$ such
that $d_{\infty}^{u}(v,w_{\varepsilon}^{\prime})\leq S\varepsilon$
for all $\varepsilon>0$ small enough and all $u$ and
$\{w_i\}_{i=1}^N$. Indeed, by the generalized triangle inequality
for Carnot groups, we have
$$
d_{\infty}^u(v,w_{\varepsilon}^{\prime})\leq
c(d_{\infty}^u(u,v)+d_{\infty}^u(u,w_{\varepsilon}^{\prime})).
$$
Next,
$d_{\infty}^u(u,w_{\varepsilon}^{\prime})=d_{\infty}(u,w_{\varepsilon}^{\prime})$.
Since $d_{\infty}(u, v)={\mathcal C}\varepsilon$, it is enough to
show that $d_{\infty}(u,w_{\varepsilon}^{\prime})\leq
K\varepsilon$. To do this, we estimate the value
$d_{\infty}(u^{\prime}, w_{\varepsilon}^{\prime})$. Since
$d_{\infty}(u^{\prime},
w_{\varepsilon}^{\prime})=d_{\infty}^{u^{\prime}}(u^{\prime},
w_{\varepsilon}^{\prime})$, then in view of the generalized
triangle inequality for Carnot groups, we have
$$
d_{\infty}^{u^{\prime}}(u^{\prime}, w_{\varepsilon}^{\prime})\leq
c(d_{\infty}^{u^{\prime}}(u^{\prime},
v)+d_{\infty}^{u^{\prime}}(v, w_{\varepsilon}^{\prime})).
$$
The conditions $d_{\infty}(u,u^{\prime})=C\varepsilon$,
$d_{\infty}(u, v)={\mathcal C}\varepsilon$ and Theorem
{\ref{estdiam}} imply
$$
d_{\infty}^{u^{\prime}}(u^{\prime}, v)=d_{\infty}(u^{\prime},
v)\leq L\max\{C,\mathcal C\}\varepsilon.
$$
Applying Theorem {\ref{estdiam}} again, we infer
$$
d_{\infty}(u,w_{\varepsilon}^{\prime})\leq K\varepsilon.
$$
From here and from the fact that $d_{\infty}(u, v)={\mathcal
C}\varepsilon$, we have
$$
d_{\infty}^u(v,w_{\varepsilon}^{\prime})\leq S\varepsilon
$$
for all $\varepsilon>0$ small enough and all $u$ and
$\{w_i\}_{i=1}^N$ belonging to some compact neighborhoods.

From here, we have that all $\alpha_i(\varepsilon)$, $i=1,\ldots,
N$, are bounded uniformly in $\varepsilon>0$.

$3^{\text{\sc rd}}$ {\sc Step.} Note that
$d_{\infty}^u(\omega_{\varepsilon},w^{\prime}_{\varepsilon})=\varepsilon
d_{\infty}^u(\omega^{\prime},w^{\prime})$. Consider the mapping
\begin{equation}\label{delta_epsilon}
\Delta_{\varepsilon,v}^u(x)=\exp\Bigl(\sum\limits_{i=1}^Nx_i\varepsilon^{\operatorname{deg}
X_i}\widehat{X}^{u}_i\Bigr)(v).
\end{equation}
More exactly,
\begin{multline}\label{d_e_exact}
\mathbb M\ni x\mapsto\{x_1,\ldots, x_N\}\text{ by such a way that
}x=\exp\Bigl(\sum\limits_{i=1}^Nx_i\widehat{X}^{u}_i\Bigr)(v)\\
\overset{\Delta_{\varepsilon,v}^u}\longmapsto
\exp\Bigl(\sum\limits_{i=1}^Nx_i\varepsilon^{\operatorname{deg}
X_i}\widehat{X}^{u}_i\Bigr)(v).
\end{multline}
Show that the coordinate functions are $H^{\alpha}$-continuous
 in $u\in\mathbb M$ uniformly on $\varepsilon>0$.

{\bf 1. The case of ${\alpha}>0$.} Indeed, the mapping
$$
\theta_{v,u}(x_1,\ldots,x_N)=\exp\Bigl(\sum\limits_{i=1}^Nx_i
\widehat{X}^{u}_i\Bigr)(v),
$$
where $(x_1,\ldots,x_N)\in\operatorname{Box}(0,T\varepsilon)$, is
$H^{\alpha}$-continuous in $u\in\mathbb M$ as a solution to an
equation with $H^{\alpha}$-continuous right-hand part (see
Section~\ref{proofdiff}), and its H\"{o}lder constant does not
depend on $v$ belonging to some compact set. This mapping is also
quasi-isometric on
$(x_1\varepsilon^{-\operatorname{deg}X_1},\ldots,x_N\varepsilon^{-\operatorname{deg}X_N})\in\mathbb
R^N$ with respect to the Riemannian metric. Consider now the
inverse mapping, which assigns to a given point $x\in\mathbb M$,
$d^u_{\infty}(v,x)\leq T\varepsilon$, the ``coordinates''
$x_1(u,x)\varepsilon^{-\operatorname{deg}X_1},\ldots,x_N(u,x)\varepsilon^{-\operatorname{deg}X_N}$
such that
$$
x=\exp\Bigl(\sum\limits_{i=1}^Nx_i(u,x)\widehat{X}^{u}_i\Bigr)(v).
$$

Note that the quasi-isometric coefficients of the mapping
$\theta_{v,u}$ are independent from $(x_1,\ldots, x_N)$, $u$ and
$v$ belonging to some compact set (here we suppose that
$d^u_{\infty}(v,x)\leq T\varepsilon$). Show that the functions
$x_1(u,x)\varepsilon^{-\operatorname{deg}X_1},\ldots,$
$x_N(u,x)\varepsilon^{-\operatorname{deg}X_N}$ are
$H^{\alpha}$-continuous in $u\in U$ for a fixed $x\in\mathbb M$,
and their H\"{o}lder constants are bounded locally uniformly in
$x$, $v$ and in $\varepsilon>0$. (Here, to guarantee the uniform
boundedness of
$x_1(u,x)\varepsilon^{-\operatorname{deg}X_1},\ldots,x_N(u,x)\varepsilon^{-\operatorname{deg}X_N}$,
we assume that
\begin{itemize}
\item both values $d_{\infty}(u,v)$ and $d^u_{\infty}(v,x)$ are
comparable to $\varepsilon$

\item the point $u$ can be changed only by a point $u^{\prime}$,
such that the distance $d_{\infty}(u,u^{\prime})$ is also
comparable to $\varepsilon$ (see $2^{\text{nd}}$ step).)
\end{itemize}

The latter statement follows from the fact, that
$\theta_{u,v}(x_1,\ldots,x_N)$ is locally H\"{o}lder in $u$, and
its H\"{o}lder constant is independent of $v$ belonging to some
compact set, and of $(x_1,\ldots,x_N)$ belonging to some compact
neighborhood $U(0)$ of zero. Since we prove a local property of a
mapping then we may assume that $u$, $u^{\prime}$, $x$ and $v$
meet our above condition on $d_{\infty}$-distances and they belong
to some compact neighborhood $U$ such that the mapping
$\theta_{u,v}$ is bi-Lipschitz on
$(x_1\varepsilon^{-\operatorname{deg}X_1},\ldots,x_N\varepsilon^{-\operatorname{deg}X_N})$
if $u\in U$; moreover, its bi-Lipschitz coefficients are
independent of $u$,
$(x_1\varepsilon^{-\operatorname{deg}X_1},\ldots,x_N\varepsilon^{-\operatorname{deg}X_N})$
and $v$ belonging to some compact set. Indeed, consider the
mapping $\theta_v(u,x_1,\ldots,x_N)=\theta_{u,v}(x_1,\ldots,x_N)$
and suppose that for any $L>0$ there exist $\varepsilon>0$, points
$v,x\in U$, a level set $\theta_{v}^{-1}(x)$, and points
$(u,x_1(u),\ldots,x_N(u))$ and $(u^{\prime},x_1(u^{\prime}),$
$\ldots,x_N(u^{\prime}))$ on it such that
\begin{multline}\label{nbl}
\Bigl|(x_1(u)\varepsilon^{-\operatorname{deg}X_1},\ldots,x_N(u)\varepsilon^{-\operatorname{deg}X_N})\\
-(x_1(u^{\prime})\varepsilon^{-\operatorname{deg}X_1},\ldots,x_N(u^{\prime})\varepsilon^{-\operatorname{deg}X_N})\Bigr|\geq
L|u-u^{\prime}|^{\alpha}
\end{multline}
for some $u$ and $u^{\prime}$. The assumption {\eqref{nbl}} leads
to the following contradiction:
\begin{multline}\label{cont_nbl}
0=\Bigl|\theta_v(u,x_1(u)\varepsilon^{-\operatorname{deg}X_1},\ldots,x_N(u)\varepsilon^{-\operatorname{deg}X_N})\\
-\theta_v(u^{\prime},x_1(u^{\prime})\varepsilon^{-\operatorname{deg}X_1},\ldots,x_N(u^{\prime})\varepsilon^{-\operatorname{deg}X_N})\Bigr|\\
\geq\Bigl|\theta_v(u,x_1(u)\varepsilon^{-\operatorname{deg}X_1},\ldots,x_N(u)\varepsilon^{-\operatorname{deg}X_N})\\
-\theta_v(u,x_1(u^{\prime})\varepsilon^{-\operatorname{deg}X_1},\ldots,x_N(u^{\prime})\varepsilon^{-\operatorname{deg}X_N})\Bigr|\\
-\Bigl|\theta_v(u,x_1(u^{\prime})\varepsilon^{-\operatorname{deg}X_1},\ldots,x_N(u^{\prime})\varepsilon^{-\operatorname{deg}X_N})\\
-\theta_v(u^{\prime},x_1(u^{\prime})\varepsilon^{-\operatorname{deg}X_1},\ldots,x_N(u^{\prime})\varepsilon^{-\operatorname{deg}X_N})\Bigr|\\
\geq
C_x\Bigl|(x_1(u)\varepsilon^{-\operatorname{deg}X_1},\ldots,x_N(u)\varepsilon^{-\operatorname{deg}X_N})\\
-(x_1(u^{\prime})\varepsilon^{-\operatorname{deg}X_1},\ldots,x_N(u^{\prime})\varepsilon^{-\operatorname{deg}X_N})\Bigr|\\
-C_u|u-u^{\prime}|^{\alpha}\geq
(LC_x-C_u)|u-u^{\prime}|^{\alpha}>0
\end{multline}
if $L>\frac{C_u}{C_x}$.

Note that
$\omega^{\prime}=\Delta_{\varepsilon^{-1},v}^u(\Delta_{\varepsilon,v}^{u^{\prime}}(w^{\prime}))$,
and
$w^{\prime}=\Delta_{\varepsilon^{-1},v}^{u^{\prime}}(\Delta_{\varepsilon,v}^{u^{\prime}}(w^{\prime}))$.
Here, for the point
$w^{\prime}_{\varepsilon}=\Delta_{\varepsilon,v}^{u^{\prime}}(w^{\prime})$,
we have
$x_i(u,w^{\prime}_{\varepsilon})=\alpha_i(\varepsilon)\cdot\varepsilon^{\operatorname{deg}
X_i}$ on the one hand, and we have
$x_i(u^{\prime},w^{\prime}_{\varepsilon})=w_i\cdot\varepsilon^{\operatorname{deg}
X_i}$ on the other hand, $i=1,\ldots, N$. Since the points $u,
u^{\prime}, v$ and $w^{\prime}_{\varepsilon}$ meet our assumption
on points, we have that the H\"older constants of
$x_i(u,x)\varepsilon^{-\operatorname{deg}X_i}$ are bounded
uniformly in $\{w_j\}_{j=1}^N$ belonging to some neighborhood of
zero. Hence,
$\rho(\omega^{\prime},w^{\prime})={\Theta}\rho(u,u^{\prime})^{\alpha}$,
and
\begin{equation}\label{rhoalpha}
d_{\infty}^u(\omega^{\prime},w^{\prime})={\Theta}\rho(u,u^{\prime})^{\frac{\alpha}{M}}.
\end{equation}

{\bf 2. The case of ${\alpha}=0$} is proved similarly to the
previous case. We prove that the functions
$x_1(u,x)\varepsilon^{-\operatorname{deg}X_1},\ldots,x_N(u,x)\varepsilon^{-\operatorname{deg}X_N}$
are uniformly continuous in $u\in U$ for a fixed $x\in\mathbb M$,
and this continuity is uniform in $x,v$ and $\varepsilon>0$. The
points under consideration meet the above condition.

To prove our result, we assume the contrary that there exists
$\sigma>0$ such that for any $\delta>0$ there exist
$\varepsilon>0$, points $v,x\in U$, a level set
$\theta_{v}^{-1}(x)$, and points $(u,x_1(u),\ldots,x_N(u))$ and
$(u^{\prime},x_1(u^{\prime}),\ldots,x_N(u^{\prime}))$ on it such
that $|u-u^{\prime}|<\delta$, and in the right-hand part of
{\eqref{nbl}} instead of $L|u-u^{\prime}|^{\alpha}$, we obtain
$\sigma$.

Repeating further the scheme of the proof almost verbatim and
replacing $(LC_x-C_u)|u-u^{\prime}|^{\alpha}$ by $\sigma
C_x-\omega_{\theta_v}(u)$ in the right-hand part of
{\eqref{cont_nbl}}, we deduce
\begin{equation}\label{rhoalpha0}
\rho(\omega^{\prime},w^{\prime})=
\omega_{\Delta^u_{\varepsilon^{-1},v}\circ\Delta^{u^{\prime}}_{\varepsilon,v}}(\rho(u,u^{\prime})).
\end{equation}
We may assume without loss of generality, that
$\omega_{\Delta^u_{\varepsilon^{-1},v}\circ\Delta^{u^{\prime}}_{\varepsilon,v}}$
does not depend on $x$ and $v$ (see {\eqref{nbl}} and
{\eqref{cont_nbl}}).

$4^{\text{\sc th}}$ {\sc Step.} Taking {\eqref{est1}},
{\eqref{rhoalpha}} and {\eqref{rhoalpha0}} into account we obtain
$$
d_{\infty}^u(w_{\varepsilon},w^{\prime}_{\varepsilon})=\varepsilon
[\Theta(C,{\mathcal C})]\rho(u,u^{\prime})^{\frac{\alpha}{M}}
$$
for $\alpha>0$. Similarly, we obtain the theorem for $\alpha=0$.
The theorem follows.
\end{proof}

\begin{cor}
\label{rem_eps} {\bf 1.} Note that
$d_{\infty}(u,u^{\prime})=C\varepsilon$ implies
$\rho(u,u^{\prime})<C\varepsilon$. Then, for $\alpha>0$, we have
$$
d_{\infty}^u(w_{\varepsilon},
w_{\varepsilon}^{\prime})=O(\varepsilon^{1+\frac{\alpha}{M}})\text
{ as }\varepsilon\to0
$$
where $O$ is uniform in $u, u^{\prime}, v\in U\subset\mathbb M$,
and in $\{w_i\}_{i=1}^N$ belonging to some compact neighborhood of
$0$, and depends on $C$ and ${\mathcal C}$.

{\bf 2.} If  $\alpha=0$ then
$$
d_{\infty}^u(w_{\varepsilon},
w_{\varepsilon}^{\prime})=o(\varepsilon)\text { as
}\varepsilon\to0
$$
where $o$ is uniform in $u, u^{\prime}, v\in U\subset\mathbb M$,
and in $\{w_i\}_{i=1}^N$ belonging to some compact neighborhood of
$0$, and depends on $C$ and ${\mathcal C}$.
\end{cor}

\begin{rem}
The estimate $O(\varepsilon^{1+\frac{\alpha}{M}})$ is also true
for the case of vector fields $X_i$, $i=1,\ldots, N$, which are
H\"{o}lder with respect to such $\mathfrak d$ that
$d_{\infty}(u,u^{\prime})=C\varepsilon$ implies $\mathfrak
d(u,u^{\prime})=K\varepsilon$, where $K$ is bounded for
$u,u^{\prime}\in U$.

A particular case is $\mathfrak d=d^z_{\infty}$, where
$d_{\infty}(z,u)\leq Q\varepsilon$ (see Local Approximation
Theorem~{\ref{apth_d_1}}, case $\alpha=0$, below).
\end{rem}

\begin{rem}
The estimate $O(\varepsilon^{1+\frac{\alpha}{M}})$ is also true
for the case of vector fields $X_i$, $i=1,\ldots, N$, which are
H\"{o}lder with respect to such $\mathfrak d$ that
$d_{\infty}(u,u^{\prime})=C\varepsilon$ implies $\mathfrak
d(u,u^{\prime})=K\varepsilon$, where $K$ is bounded for
$u,u^{\prime}\in U$.

A particular case is $\mathfrak d=d^z_{\infty}$, where
$d_{\infty}(z,u)\leq Q\varepsilon$ (see Local Approximation
Theorem~{\ref{apth_d_1}}, case $\alpha=0$, below).
\end{rem}

\subsection{The Approximation Theorems}\label{apprths}

In this subsection, we prove two Approximation Theorems. Their
proofs use the following geometric property.

\begin{prop}
\label{prop21vg} For a neighborhood $U$, there exist positive
constants $C>0$ and $r_0>0$ depending on $U$, $M$, and $N$, such
that for any points $u$ and $v$ from a neighborhood $U$ the
following inclusion is valid:
$$
\bigcup\limits_{x\in\operatorname{Box}^u (v,r)}
\operatorname{Box}^u (x,\xi)\subseteq \operatorname{Box}^u
(v,r+C\xi),\quad 0<\xi,\,  r\leq r_0.
$$
\end{prop}

\begin{proof}
Let $x=\exp\Bigl(\sum\limits_{i=1}^Nx_i\widehat{X}^u_i\Bigr)(v)$,
$d_{\infty}^u(v,x)\leq r$, and
$z=\exp\Bigl(\sum\limits_{i=1}^Nz_i\widehat{X}^u_i\Bigr)(x)$,
$d_{\infty}^u(x,z)\leq\xi$. We estimate the distance
$d_{\infty}^u(v,z)$ applying {\eqref{group}} to points $x$ and
$z$. Let
$z=\exp\Bigl(\sum\limits_{i=1}^N\zeta_i\widehat{X}^u_i\Bigr)(v)$.

{\sc Case} of $\operatorname{deg} X_i=1$. Then $|\zeta_i|\leq
|x_i|+|z_i|\leq (r+\xi)^{\operatorname{deg} X_i}$.

{\sc Case} of $\operatorname{deg} X_i=2$. Then
\begin{multline*}
|\zeta_i|\leq |x_i|+|z_i|
+\sum\limits_{\substack{|e_l+e_j|_h=2,\\l<j}}|{F}^i_{e_l,e_j}(u)||x_lz_j-z_lx_j|\\
\leq r^2+\xi^2+c_i(u)r\xi\leq
r^2+2r\frac{c_i(u)}{2}\xi+\Bigl(\frac{c_i(u)}{2}\xi\Bigr)^2\\
=\Bigl(r+\frac{c_i(u)}{2}\xi\Bigr)^{\operatorname{deg}
X_i}=(r+C_i(u)\xi)^{\operatorname{deg} X_i}.
\end{multline*}
Here we assume that $C_i(u)\geq 1$.

{\sc Case} of $\operatorname{deg} X_i=k>2$. Then we obtain
analogously to the previous case
\begin{multline*}
|\zeta_i|\leq |x_i|+|z_i|+\sum\limits_{|\mu+\beta|_h=k,
\mu>0,\beta>0}|{F}^i_{\mu,\beta}(u)|x^\mu\cdot
z^\beta\\
\leq
r^k+\xi^k+\sum\limits_{|\mu+\beta|_h=k}c_i^{\mu\beta}(u)r^{|\mu|_h}\xi^{|\beta|_h}\leq
(r+C_i(u)\xi)^{\operatorname{deg} X_i}.
\end{multline*}
Here we assume that $C_i(u), c_i(u)\geq 1$. Denote by
$C(u)=\max\limits_iC_i(u)$. From above estimates we obtain
$$
d_{\infty}^u(v,x)=\max\limits_i\{|\zeta_i|^{\operatorname{deg}
X_i}\}\leq \max\limits_i\{(r+C_i(u)\xi)^{\frac{\operatorname{deg}
X_i}{\operatorname{deg} X_i}}\}\leq r+C(u)\xi.
$$
Since all the $C_i(u)$'s are continuous on $u$ then we may choose
$C<\infty$ such that $C(u)\leq C$ for all $u$ belonging to a
compact neighborhood. The lemma follows.
\end{proof}

\begin{thm}[Approximation Theorem]
\label{lat} Let $u,u^{\prime},v,w\in U$. Then the following
estimate is valid:
\begin{equation}\label{apthglob}
|d_{\infty}^u(v,w)-d_{\infty}^{u^{\prime}}(v,w)|={\Theta}[\rho(u,u^{\prime})^{\alpha}
 \rho(v,w)]^{\frac{1}{M}}.
\end{equation}
\end{thm}

\begin{proof}
Let $p=\exp\Bigl(\sum\limits_{i=1}^{N}p_i\widehat{X}_i^u\Bigr)(v)$
and
$p^{\prime}=\exp\Bigl(\sum\limits_{i=1}^{N}p_i\widehat{X}_i^{u^{\prime}}\Bigr)(v)$.
Notice that if $z\in\operatorname{Box}^u(v, d_{\infty}^u(v,w))$
then $z^{\prime}\in\operatorname{Box}^{u^{\prime}}(v,
d_{\infty}^u(v,w))$ and
$z\in\operatorname{Box}^{u^{\prime}}(z^{\prime},
R(u,u^{\prime}))$, where
$$
R(u,u^{\prime})=\sup\limits_{p^{\prime}\in\operatorname{Box}^{u^{\prime}}(v,
d_{\infty}^u(v,w))}d_{\infty}^{u^{\prime}}(p,p^{\prime}).
$$
Using Proposition {\ref{prop21vg}} we have that
\begin{multline*}
\operatorname{Box}^u(v,
d_{\infty}^u(v,w))\subset\bigcup\limits_{x\in\operatorname{Box}^{u^{\prime}}(v,
d_{\infty}^u(v,w))}\operatorname{Box}^{u^{\prime}}(x,R(u,u^{\prime}))\\
\subset \operatorname{Box}^{u^{\prime}}(v,
d_{\infty}^u(v,w)+CR(u,u^{\prime}))
\end{multline*}
for some $C>0$. Consequently, in view of Theorem
{\ref{difference}} we can write
\begin{multline*}
\operatorname{Box}^u(v, d_{\infty}^u(v,w))\subset
\operatorname{Box}^{u^{\prime}}(v, d_{\infty}^u(v,w)+CR(u,u^{\prime}))\subset\\
\operatorname{Box}^{u^{\prime}}(v,
d_{\infty}^u(v,w)+{\Theta}[\rho(u,u^{\prime})^{\alpha}
\rho(v,w)]^{\frac{1}{M}}).
\end{multline*}

If $d_{\infty}^u(v,w)\leq
{\Theta}[\rho(u,u^{\prime})^{\alpha}\rho(v,w)]^{\frac{1}{M}}$ then
the theorem follows:
$$
|d_{\infty}^u(v,w)-d_{\infty}^{u^{\prime}}(v,w)|\leq
d_{\infty}^u(v,w)+d_{\infty}^{u^{\prime}}(v,w)\\
={\Theta}[\rho(u,u^{\prime})^{\alpha}
 \rho(v,w)]^{\frac{1}{M}}.
$$

If $d_{\infty}^u(v,w)>{\Theta}[\rho(u,u^{\prime})^{\alpha}
 \rho(v,w)]^{\frac{1}{M}}$
then applying again Proposition {\ref{prop21vg}} we obtain
$$
\operatorname{Box}^{u^{\prime}}(v,
d_{\infty}^u(v,w)-{\Theta}[\rho(u,u^{\prime})^{\alpha}
\rho(v,w)]^{\frac{1}{M}})\\
\subset \operatorname{Box}^u(v, d_{\infty}^u(v,w)).
$$
From the latter relation it follows that
\begin{multline*}
d_{\infty}^u(v,w)-{\Theta}[\rho(u,u^{\prime})^{\alpha}
 \rho(v,w)]^{\frac{1}{M}}\leq
d_{\infty}^{u^{\prime}}(v,w)\\
\leq d_{\infty}^u(v,w)+{\Theta}[\rho(u,u^{\prime})^{\alpha}
 \rho(v,w)]^{\frac{1}{M}},
\end{multline*}
and the theorem follows.
\end{proof}

\begin{rem}
If the derivatives of $X_i$, $i=1,\ldots, N$, are locally
H\"{o}lder with respect to $\mathfrak d$, then we have $\mathfrak
d(u,u^{\prime})^{\alpha}$ instead of $\rho(u,u^{\prime})^{\alpha}$
in {\eqref{apthglob}} (the proof is similar).
\end{rem}

Approximation Theorem and local estimates (see Theorem
{\ref{est_eps}}) imply Local Approximation Theorem.

\begin{thm}[Local Approximation Theorem]
\label{apth_d_1} Assume that
$d_{\infty}(u,u^{\prime})=C\varepsilon$, $d_{\infty}(u,
v)={\mathcal C}\varepsilon$ and $d_{\infty}(u, w)={\mathbb
C}\varepsilon$ for some $C,{\mathcal C},\mathbb C<\infty$.

{\bf 1.} If $\alpha>0$, then
\begin{equation}\label{apthloc}
|d_{\infty}^u(v,w)-d_{\infty}^{u^{\prime}}(v,w)|=\varepsilon
{\Theta}[\rho(u,u^{\prime})]^{\frac{\alpha}{M}}\Theta(d_{\infty}^u(v,w)+o(1)).
\end{equation}
Moreover, if $u^{\prime}=v$ and $\alpha>0$, then
\begin{equation*}
|d_{\infty}^u(v,w)-d_{\infty}(v,w)|=\varepsilon
{\Theta}[\rho(u,v)]^{\frac{\alpha}{M}}\Theta(d_{\infty}^u(v,w)+o(1)).
\end{equation*}

{\bf 2.} If  $\alpha=0$, then
$$
|d_{\infty}^u(v,w)-d_{\infty}^{u^{\prime}}(v,w)|=\varepsilon
o(1)=o(\varepsilon)
$$
where $o$ is uniform in $u, u^{\prime}, v, w\in U\subset\mathbb
M$. Moreover, if $u^{\prime}=v$ and $\alpha=0$, then
$$
|d_{\infty}^u(v,w)-d_{\infty}(v,w)|=o(\varepsilon)
$$
where $o$ is uniform in $u, v, w\in U\subset\mathbb M$.
\end{thm}

\noindent {\it Proof} follows the scheme as the proof of
Approximation Theorem~{\ref{lat}} with
$R(u,u^{\prime})=\varepsilon [\Theta(C,\mathcal C,\mathbb
C)]\rho(u,u^{\prime})^{\frac{\alpha}{M}}$. The latter equality is
valid by the uniformity assertion of Theorem {\ref{est_eps}}.

\begin{rem}
If the derivatives of $X_i$, $i=1,\ldots, N$, are locally
H\"{o}lder with respect to $\mathfrak d$, then we have $\mathfrak
d(u,u^{\prime})^{\frac{\alpha}{M}}$ instead of
$\rho(u,u^{\prime})^{\frac{\alpha}{M}}$ in {\eqref{apthloc}} (the
proof is similar).
\end{rem}

\subsection{Comparison of Local Geometries of Two Local Carnot Groups}\label{geomconecone}

\begin{proof}[Proof of Theorem {\ref{est_chain}}]
$1^{\text{\sc st}}$ {\sc Step.} Consider the case of $\alpha>0$.
The case of $Q=1$ is proved in Theorem {\ref{est_eps}}.

$2^{\text{\sc nd}}$ {\sc Step.} Consider the case of $Q=2$. First,
for the points $w_2=w^1_2$ and $w_2^{\prime}={w_2^1}^{\prime}$, we
have
\begin{equation}\label{step2_1}
w_2=\exp\Bigl(\sum\limits_{i=1}^N\omega_{i,2}\widehat{X}^u_i\Bigr)(w_0)
\end{equation}
and
\begin{equation}\label{step2_2}
w_2^{\prime}=\exp\Bigl(\sum\limits_{i=1}^N{\omega_{i,2}^{\prime}}\widehat{X}^{u^{\prime}}_i\Bigr)(w_0).
\end{equation}
By the formulas of group operation, $\omega_{i,2}$ differs from
$\omega_{i,2}^{\prime}$ in the values of
$\{F^{j}_{\mu,\beta}(u)\}_{j,\mu,\beta}$. By Assumption
{\ref{vectorfields}},
$F^{j}_{\mu,\beta}(u^{\prime})=F^{j}_{\mu,\beta}(u)+\Theta\rho(u,u^{\prime})^{\alpha}$.

Consider the auxiliary points
$$
w_2^{\prime\prime}=\exp\Bigl(\sum\limits_{i=1}^N\omega_{i,2}\widehat{X}^{u^{\prime}}_i\Bigr)(w_0)\text{
and
}{w_2^{\prime\prime}}^{\varepsilon}=\exp\Bigl(\sum\limits_{i=1}^N\omega_{i,2}\varepsilon^{\operatorname{deg}
X_i}\widehat{X}^{u^{\prime}}_i\Bigr)(w_0)
$$
and estimate the value
$d_{\infty}^{u^{\prime}}(w_2^{\prime\prime},w_2^{\prime})$. For
doing this, we use the group operation in the local Carnot group
$\mathcal G^{u^{\prime}}\mathbb M$ and Approximation Theorem
{\ref{lat}}. Note that,
$|\omega_{i,2}-{\omega_{i,2}^{\prime}}|=\Theta\rho(u,u^{\prime})^{\alpha}$.
Next, note that while applying the group operation, all summands
look like $\omega_{i,2}-{\omega_{i,2}^{\prime}}$ or
$\omega_{i,2}-{\omega_{i,2}^{\prime}}+\sum\Theta(\omega_{k,2}{\omega_{j,2}^{\prime}}-\omega_{j,2}{\omega_{k,2}^{\prime}})$.
By~{\eqref{group}}, we deduce
\begin{multline*}
\omega_{k,2}{\omega_{j,2}^{\prime}}-\omega_{j,2}{\omega_{k,2}^{\prime}}\\
=\omega_{k,2}(\omega_{j,2}+{\Theta}\rho(u,u^{\prime})^{\alpha})
-\omega_{j,2}(\omega_{k,2}+{\Theta}\rho(u,u^{\prime})^{\alpha})={\Theta}\rho(u,u^{\prime})^{\alpha},
\end{multline*}
$d_{\infty}^{u^{\prime}}(w_2^{\prime\prime},w_2^{\prime})={\Theta}(\rho(u,u^{\prime})^{\frac{\alpha}{M}})$.
Here $\Theta$ depends on $C$, ${\mathcal C}$, $Q=2$ and
$\{F^j_{\mu,\beta}(u^{\prime})\}_{j,\mu,\beta}$.

It follows from the formulas of group operation in $\mathcal
G^{u}\mathbb M$ and $\mathcal G^{u^{\prime}}\mathbb M$, that
$$
w_2^{\varepsilon}=\exp\Bigl(\sum\limits_{i=1}^N\omega_{i,2}\varepsilon^{\operatorname{deg}
X_i}\widehat{X}^u_i\Bigr)(w_0)
$$
and
$$
{w_2^{\varepsilon}}^{\prime}=\exp\Bigl(\sum\limits_{i=1}^N{\omega_{i,2}^{\prime}}\varepsilon^{\operatorname{deg}
X_i}\widehat{X}^{u^{\prime}}_i\Bigr)(w_0).
$$
By Theorem {\ref{est_eps}}, we have
$d_{\infty}^{u^{\prime}}({w_2^{\prime\prime}}^{\varepsilon},w_2^{\varepsilon})=\varepsilon
{\Theta}\rho(u,u^{\prime})^{\frac{\alpha}{M}}$. By the homogeneity
of the distance $d_{\infty}^{u^{\prime}}$ we have
$$
d_{\infty}^{u^{\prime}}({w_2^{\prime\prime}}^{\varepsilon},w_2^{\varepsilon\prime})=\varepsilon
d_{\infty}^{u^{\prime}}(w_2^{\prime\prime},w_2^{\prime})=\varepsilon
{\Theta}\rho(u,u^{\prime})^{\frac{\alpha}{M}},
$$
and from the generalized triangle inequality we deduce
$$
d_{\infty}^{u^{\prime}}(w_2^{\varepsilon},{w_2^{\varepsilon}}^{\prime})=\varepsilon
{\Theta}\rho(u,u^{\prime})^{\frac{\alpha}{M}}.
$$
In view of Local Approximation Theorem {\ref{apth_d_1}}, we derive
$$
d_{\infty}^u(w_2^{\varepsilon},{w_2^{\varepsilon}}^{\prime})=\varepsilon
{\Theta}\rho(u,u^{\prime})^{\frac{\alpha}{M}}.
$$

$3^{\text{\sc rd}}$ {\sc Step.} In the case of $Q=3$, it is easy
to see from the previous case and the group operation, that if
$$
w_3=\exp\Bigl(\sum\limits_{i=1}^N\omega_{i,3}\widehat{X}^u_i\Bigr)(w_0)
$$
and
$$
w_3^{\prime}=\exp\Bigl(\sum\limits_{i=1}^N{\omega_{i,3}^{\prime}}\widehat{X}^{u^{\prime}}_i\Bigr)(w_0),
$$
then again
$|\omega_{i,3}-{\omega_{i,3}^{\prime}}|={\Theta}\rho(u,u^{\prime})^{\alpha}$.
Here $\Theta$ depends on $C$, ${\mathcal C}$, $Q=3$ and
$\{F^j_{\mu,\beta}\}_{j,\mu,\beta}$. (It suffices to apply the
group operation in local Carnot groups $\mathcal G^u\mathbb M$ and
$\mathcal G^{u^{\prime}}\mathbb M$ to expressions
{\eqref{step2_1}} and {\eqref{step2_2}} and to points $w_3$ and
$w_3^{\prime}$, respectively.) From now on, for obtaining estimate
{\eqref{est_l}} at $Q=3$, we repeat the arguments of the
$2^{\text{nd}}$ {Step}.

$4^{\text{\sc th}}$ {\sc Step.} It is easy to see analogously to
the $3^{\text{rd}}$ {Step}, that the group operation and the
induction hypothesis
$|\omega_{i,l-1}-{\omega_{i,l-1}^{\prime}}|={\Theta}\rho(u,u^{\prime})^{\alpha}$,
$3<l<Q$, imply
$|\omega_{i,l}-{\omega_{i,l}^{\prime}}|={\Theta}\rho(u,u^{\prime})^{\alpha}$.
Indeed, it suffices to put $\omega_{i,l}$ and
${\omega_{i,l}^{\prime}}$ instead of $\omega_{i,3}$ and
${\omega_{i,3}^{\prime}}$, and $\omega_{i,l-1}$ and
${\omega_{i,l-1}^{\prime}}$ instead of $\omega_{i,2}$ and
${\omega_{i,2}^{\prime}}$ in the $3^{\text{rd}}$ Step, and apply
arguments from the $2^{\text{nd}}$ {Step}.

The case of $\alpha=0$ can be proved by applying the similar
arguments.

The theorem follows.
\end{proof}

\subsection{Comparison of Local Geometries of a Carnot Manifold and a Local Carnot Group}\label{geomcarnotcone}

In this subsection, we compare the local geometry of a Carnot
manifold with the one of a local Carnot group.

\begin{thm}
\label{mainresult} Fix $Q\in\mathbb N$. Consider points $w_0$, $u$
such that $d_{\infty}(u,w_0)={\mathcal C}\varepsilon$ for some
${\mathcal C}<\infty$, and
$$
\widehat{w}_j^{\varepsilon}=\exp\Bigl(\sum\limits_{i=1}^{N}w_{i,j}{\varepsilon}^{\operatorname{deg}
X_i}\widehat{X}_i^u\Bigr)(\widehat{w}^{\varepsilon}_{j-1}),\quad
{w_j^{\varepsilon}}=\exp\Bigl(\sum\limits_{i=1}^{N}w_{i,j}{\varepsilon}^{\operatorname{deg}
X_i}{X}_i\Bigr)({w^{\varepsilon}_{j-1}}),
$$
${w^{\varepsilon}_0}=\widehat{w}^{\varepsilon}_0=\widehat{w}_0=w_0$,
$j=1,\ldots,Q$. {$($Here $Q\in\mathbb N$ is such that all these
points belong to a neighborhood $U\subset\mathbb M$ small enough
for all $\varepsilon>0$.$)$} Then for $\alpha>0$
\begin{equation}\label{mainestimate}
\max\{d_{\infty}^u(\widehat{w}^{\varepsilon}_Q,{w^{\varepsilon}_Q}),
d_{\infty}(\widehat{w}^{\varepsilon}_Q,{w^{\varepsilon}_Q})\}=
{\sum\limits_{k=1}^Q \Theta({\mathcal
C},k,\{F^j_{\mu,\beta}\}_{j,\mu,\beta})}\cdot
\varepsilon^{1+\frac{\alpha}{M}}.
\end{equation}
In the case of $\alpha=0$ we have
$$
\{d_{\infty}^u(\widehat{w}^{\varepsilon}_Q,{w^{\varepsilon}_Q}),
d_{\infty}(\widehat{w}^{\varepsilon}_Q,{w^{\varepsilon}_Q})\}=\varepsilon\cdot
\Theta({\mathcal
C},Q,\{F^j_{\mu,\beta}\}_{j,\mu,\beta})[\omega(\varepsilon)]^{\frac{1}{M}}
$$
where $\omega(\varepsilon)\to0$ as $\varepsilon\to0$.
 Here $|w_{i,j}|$ are bounded, and ${\Theta}$ is uniformly bounded for $u, w_0\in U$
and $\{w_{i,j}\}$, $i=1,\ldots, N$, $j=1,\ldots,Q$, belonging to
some compact neighborhood of $0$, and it depends on $Q$ and
$\{F^j_{\mu,\beta}\}_{j,\mu,\beta}$.
\end{thm}

\begin{proof}
For simplifying the notation we denote the points
$\widehat{w}_i^1$ by $\widehat{w}_i$, and we denote ${w}_i^1$ by
${w}_i$ for $\varepsilon=1$. First, consider the points
$\widehat{w}_Q$ and $w_Q$. Now we construct a following sequence
of points.

Let
\begin{multline*}
\omega_{k,j}=\exp\Bigl(\sum\limits_{i=1}^{N}w_{i,j}\widehat{X}_i^{w_k}\Bigr)(\omega_{k,{j-1}}),\\
k=0,\ldots, Q-1,\ j=1,\ldots, Q-k,\ \omega_{k,0}=w_k.
\end{multline*}
Hence, $\omega_{Q-1,1}=w_Q$ and
$$
d_{\infty}^u(w_Q,w^{\prime}_Q)=O\Bigl(d_{\infty}^u(w_Q,
\omega_{0,Q})+\sum\limits_{k=1}^{Q-1}d_{\infty}^u(\omega_{k,{Q-k}},\omega_{{k-1},{Q-k+1}})\Bigr).
$$
If $\alpha>0$ then, by Theorem {\ref{est_chain}},
$$
d_{\infty}^u(\widehat{w}_Q, \omega_{0,Q})=\Theta({\mathcal C},
Q,\{F^j_{\alpha,\beta}\}_{j,\alpha,\beta})\rho(u,w_0)^{\frac{\alpha}{M}},
$$
and each of the summands
$$
d_{\infty}^u(\omega_{k,{Q-k}},\omega_{{k-1},{Q-k+1}})=\Theta({\mathcal
C}, Q-k,\{F^j_{\alpha,\beta}\}_{j,\alpha,\beta})\rho(w_k,
w_{k-1})^{\frac{\alpha}{M}}.
$$
By the same theorem, if we replace $w_{i,j}$ by
$w_{i,j}\varepsilon^{\operatorname{deg} X_i}$ then it is easy to
see using induction by $k$ that firstly
$d_{\infty}^u(w_k^{\varepsilon},
w_{k-1}^{\varepsilon})=O(\varepsilon)$, secondly
$d_{\infty}^u(u,w^{\varepsilon}_k)\sim\varepsilon$ and
$d_{\infty}^u(u,\omega^{\varepsilon}_{k,Q-k})\sim\varepsilon$ for
all $k$, and thirdly
$$
d_{\infty}^u(\widehat{w}^{\varepsilon}_Q,
\omega^{\varepsilon}_{0,Q})=\varepsilon \Theta({\mathcal C},
Q,\{F^j_{\alpha,\beta}\}_{j,\alpha,\beta})\rho(u,w_0)^{\frac{\alpha}{M}}
$$
and
$$
d_{\infty}^u(\omega^{\varepsilon}_{k,{Q-k}},\omega^{\varepsilon}_{{k-1},{Q-k+1}})=\varepsilon
\Theta({\mathcal C},
Q-k,\{F^j_{\alpha,\beta}\}_{j,\alpha,\beta})\rho(w_k^{\varepsilon},
w_{k-1}^{\varepsilon})^{\frac{\alpha}{M}}.
$$ Thus we obtain
$d_{\infty}^u(\widehat{w}^{\varepsilon}_Q,{w^{\varepsilon}_Q})=\sum\limits_{k=1}^Q
\Theta({\mathcal C},
k,\{F^j_{\alpha,\beta}\}_{j,\alpha,\beta})\cdot
\varepsilon^{1+\frac{\alpha}{M}}$.

Since
$d_{\infty}^u(\widehat{w}^{\varepsilon}_Q,{w^{\varepsilon}_Q})=O(\varepsilon)$
and $d_{\infty}^u(u,\widehat{w}^{\varepsilon}_Q)=O(\varepsilon)$
then, by Local Approximation Theorem  {\ref{apth_d_1}}, we have
$$
d_{\infty}(\widehat{w}^{\varepsilon}_Q,{w^{\varepsilon}_Q})=\sum\limits_{k=1}^Q
\Theta({\mathcal C},
k,\{F^j_{\alpha,\beta}\}_{j,\alpha,\beta})\cdot
\varepsilon^{1+\frac{\alpha}{M}}.
$$

If $\alpha=0$, then we repeat the above arguments replacing
$\rho(\cdot,\cdot)^{\frac{\alpha}{M}}$ by $o(1)$. The theorem
follows.
\end{proof}

\begin{rem}\label{mainrem}
If the derivatives of $X_i$, $i=1,\ldots, N$, are locally
H\"{o}lder with respect to $\mathfrak d$, such that
$d_{\infty}(x,y)\leq\varepsilon$ implies $\mathfrak d(x,y)\leq
K\varepsilon$, where K is bounded on $U$, then the same estimate
as in {\eqref{mainestimate}} is true (the proof is similar).

A particular case of such $\mathfrak d$ is $d_{\infty}^z$,
$d_{\infty}(u,z)\leq Q\varepsilon$ (see Local Approximation
Theorem~{\ref{apth_d_1}}).
\end{rem}

\subsection{Applications}\label{applications}

\subsubsection{Rashevski\v{\i}--Chow  Theorem}

\begin{defn}
An absolutely continuous curve $\gamma:[0,a]\to\mathbb M$ is said
to be {\it horizontal} if $\dot{\gamma}(t)\in H_{\gamma(t)}\mathbb
M$ for almost all $t\in [0,a]$. Its length $l(\gamma)$ equals
$\int\limits_{0}^a|\dot{\gamma}(t)|\,dt$, where the value
$|\dot{\gamma}(t)|$ is calculated using the Riemann tensor on
$\mathbb M$. Analogously, the canonical Riemann tensor on
${\mathcal G}^{u}\mathbb M$ defines a length $\widehat{l}$ of an
absolutely continuous curve $\widehat\gamma:[0,a]\to{\mathcal
G}^{u}\mathbb M$.
\end{defn}

\begin{defn}
The {\it Carnot--Carath\'eodory distance} between points
$x,y\in\mathbb M$ is defined as $d_{c}(x,y)=\inf\limits_\gamma
l(\gamma)$ where the infimum is taken over all horizontal curves
with endpoints $x$ and $y$.
\end{defn}

\begin{cor}[of Theorem {\ref{mainresult}}]
\label{curves} Suppose that Assumption {\ref{vectorfields}} holds
for $\alpha\in(0,1]$. Let $g\in\mathbb M$. Let also $\varepsilon$
be small enough, and $u, v,
\widehat{w}\in\operatorname{Box}(g,\varepsilon)$. The points $v,
\widehat{w}\in\operatorname{Box}(g,\varepsilon)$ can be joined in
the local Carnot group
 $({\mathcal G}^{u}\mathbb M,d^u_1)\supset \operatorname{Box}(g,\varepsilon)$ by a
horizontal curve $\widehat{\gamma}$ composed by at most $L$
segments of integral curves of horizontal fields $\widehat X_i^u$,
$i=1,\ldots,\dim H_1$. To the curve $\widehat{\gamma}$ it
corresponds a curve $\gamma$, horizontal with respect to the
initial horizontal distribution $H\mathbb M$, constituted by at
most $L$ segments of integral curves of the given horizontal
fields $X_i$, $i=1,\ldots,n$. Moreover,

\begin{enumerate}
\item the curve $\gamma$ has endpoints
$v,w\in\operatorname{Box}(g,O(\varepsilon))$;

\item
$|l(\gamma)-\widehat{l}(\widehat{\gamma})|=o(\widehat{l}(\widehat{\gamma}))$;

\item $\max\{d_{\infty}^u(\widehat{w},w),
d_{\infty}(\widehat{w},w)\}\leq
C\varepsilon^{1+{\frac{\alpha}{M}}}$ where $C$ is independent of
$g,u, v, \widehat{w}$ in some compact set.
\end{enumerate}
\end{cor}

\begin{proof}
The desired curve comes from those on any Carnot group \cite{fs}:
given a Carnot group $\mathbb G$ with the vector fields
$\widehat{X}_1,\ldots,\widehat{X}_N$, each point $x$ can be joint
with $0$ by a horizontal curve $\gamma$ constituted by at most $L$
segments $\gamma_{j}$, $j=1,\ldots, L$, of integral curves of the
given basic horizontal vector fields
$\widehat{X}_1,\widehat{X}_2,\ldots,\widehat{X}_{\dim H_1}$,
i.~e.,
$$
\begin{cases}
\dot{\gamma}_1(t)=a_{1}\widehat{X}_{i_1}(\gamma_1(t))\\
\gamma_1(0)=0,
\end{cases}
$$
$$
\begin{cases}
\dot{\gamma}_j(t)=a_{j}\widehat{X}_{i_j}(\gamma_j(t))\\
\gamma_j(0)=\gamma_{j-1}(1),
\end{cases}
$$
$j=2,\ldots, L$, and from here we have $x=\gamma_{i_L}(1)$. By
another words,
\begin{equation*}
x=\exp(a_{L}\widehat{X}_{i_L})\circ\dots\circ
\exp(a_{1}\widehat{X}_{i_1}), \quad i_j=1,\ldots,\dim H_1,
\end{equation*}
 where $|a_{j}|$ is controlled by the
distance $d_{c}(0,x)$, $j=1,\ldots,L$, and $L$ is independent of
$x$.

Now we carry over a construction described above to the local
Carnot group $({\mathcal G}^{u}\mathbb M,d^u_1)\supset
\operatorname{Box}(g,\varepsilon)$: the given points
$\widehat{w},v\in {\mathcal G}^{u}\mathbb M$ can be joint by a
horizontal curve $\widehat{\gamma}$:
\begin{equation}\label{curve2}
w=\exp(a_{L}{\widehat X}_{i_L})\circ\dots\circ \exp(a_{1}{\widehat
X}_{i_1})(v), \quad i_j=1,\ldots,\dim H_1,
\end{equation}
$j=1,\ldots,L$. Then the curve $\gamma$ defined as
\begin{equation}\label{curve3}
w=\exp(a_{L}{X}_{i_L})\circ\dots\circ \exp(a_{1}{ X}_{i_1})(v),
\quad i_j=1,\ldots,\dim H_1,
\end{equation}
is horizontal and its length equals
$\widehat{l}(\widehat{\gamma})(1+o(1))$. The estimate
$$\max\{d_{\infty}^u(w,w^{\prime}), d_{\infty}(w,w^{\prime})\}\leq
C\varepsilon^{1+{\frac{\alpha}{M}}}$$ follows immediately from
\eqref{mainestimate}.
\end{proof}

\begin{thm}
\label{locestcc} Suppose that Assumption {\ref{vectorfields}}
holds for some $\alpha\in(0,1]$. Let $g\in\mathbb M$. Given two
points $w,v\in B(g,\varepsilon)$ where $\varepsilon$ is small
enough, there exist a curve $\gamma$, horizontal with respect to
the initial horizontal distribution $H\mathbb M$,  with endpoints
$w$ and $v$, and a horizontal curve $\widehat{\gamma}$ in the
local Carnot group $({\mathcal G}^{g}\mathbb M,d^g_1)$ with the
same endpoints, such that
\begin{enumerate}
\item $\widehat{l}(\widehat{\gamma})$ is equivalent to
$d_{\infty}^g(w,v)$$;$

\item
$|l(\gamma)-\widehat{l}(\widehat{\gamma})|=o(\widehat{l}(\widehat{\gamma}))$$;$

\item if $v=g$ then the length $l(\gamma)$ is equivalent to
$d_{\infty}(g,w)$.
\end{enumerate}
\noindent All these estimates are uniform in $w$, $v$ and $g$ of
some compact neighborhood as $\varepsilon\to0$.
\end{thm}

\begin{proof}
We can choose $\varepsilon$ from the condition of the theorem by
requests
$C^{2+{\frac{\alpha}{M}}}\varepsilon^{{\frac{\alpha^2}{M^2}}}\leq1$
and $\varepsilon\leq\frac12$, where $C$ is the constant from
Corollary {\ref{curves}}.

 Apply Corollary {\ref{curves}} to the points $u=g$, $v$ and $w$.
It gives a  horizontal curve $\gamma_1$ ($\widehat{\gamma}$)  with
respect to the initial horizontal distribution $H\mathbb M$ (in
the local Carnot group $({\mathcal G}^{g}\mathbb M,d^g_1)$)
 with  endpoints $v$ and $w_1$ ($v$ and $w$) constituted by at most
$L$ segments of integral curves of given horizontal fields $X_i$
($\widehat X_i^g$), $i=1,\ldots,n$. The curves $\widehat{\gamma}$
and $\gamma_1$
 have lengths comparable with $d^g_1(v,w)$,
 and $ \max\{d_{\infty}^g(w_1,w),
d_{\infty}(w_1,w)\}\leq C\varepsilon^{1+{\frac{\alpha}{M}}}$.

Next, we apply again Corollary {\ref{curves}} to the points
$u=v=w_1$ and $w$. It gives a horizontal curve $\gamma_2$ with
respect to
 $H\mathbb M$  with  endpoints $w_1$ and $w_2$.
Its length is $O(\varepsilon^{1+{\frac{\alpha}{M}}})$ where $O$ is
uniform in $u,w\in\operatorname{Box}(g,\varepsilon)$,
 and
$d_{\infty}(w_2,w)\leq
C(C\varepsilon^{1+{\frac{\alpha}{M}}})^{1+{\frac{\alpha}{M}}}
 \leq \varepsilon^{1+{\frac{2\alpha}{M}}}$.

 Assume that
we have points $w_1, \ldots, w_k$ and horizontal curves
$\gamma_{l}$, $l=2,\ldots, k$, with respect to $H\mathbb M$ with
endpoints $w_{l-1}$ and $w_{l}$, such that $\gamma_{l}$ has a
length $O(\varepsilon^{1+{\frac{l-1}{M}}\alpha})$, and
$d_{\infty}(w_l,w)\leq \varepsilon^{1+{\frac{l\alpha}{M}}}$.

We continue,  by the induction, applying  Corollary {\ref{curves}}
to the points $u=v=w_k$ and $w$. It results a horizontal curve
$\gamma_{k+1}$ with endpoints $w_{k}$ and $w_{k+1}$, such that
$\gamma_{k+1}$ has a length
$O(\varepsilon^{1+{\frac{k\alpha}{M}}})$ and
$d_{\infty}(w_{k+1},w)\leq
C(C\varepsilon^{1+{\frac{k\alpha}{M}}})^{1+{\frac{\alpha}{M}}}\leq
\varepsilon^{1+{\frac{k+1}{M}\alpha}}$.

A curve $\Gamma_m=\gamma_1\cup\ldots\cup\gamma_m$ is horizontal,
has endpoints $v$ and $w_m$, its length does not exceed
$l(\gamma_1)+C\sum\limits_{l=1}^\infty
\varepsilon^{1+{\frac{l\alpha}{M}}}\leq l(\gamma_1)
+C\varepsilon^{1+{\frac{\alpha}{M}}}$ and $d_{\infty}(w_m,w)\to0$
as $m\to \infty$. Therefore the sequence $\Gamma_m$ converges to a
horizontal curve $\gamma$ as $m\to\infty$ with properties 1--2
mentioned in the theorem.

Under $v=g$ we can take $d_{\infty}(g,w)$ as $\varepsilon$ in
above estimates: it gives an evaluation $l(\gamma)\leq
Cd_{\infty}(g,w)$. The opposite inequality can be verified
directly {by means of the above obtained estimate: { indeed}, if
$d_{\infty}(g,w)=\varepsilon$ then $d_{\infty}(g,w)=d^g_1(g,w)\leq
C\widehat{l}(\widehat{\gamma})\leq
Cl(\gamma)+o(\widehat{l}(\widehat{\gamma}))$; it follows that
$d_{\infty}(g,w)-o\bigl(d_{\infty}(g,w)\bigr)\leq Cl(\gamma)$ and
the estimate $d_{\infty}(g,w)\leq C_1l(\gamma)$ holds with $C_1$
independent of $g$ from some compact neighborhood if $v$ is close
enough to~$g$}. Thus we have obtained the property~3.
\end{proof}

As an application of Theorem~{\ref{locestcc}} we obtain a  version
of Rashevski\v{\i}--Chow  type connectivity theorem.

\begin{thm}
\label{Chow} Suppose that Assumption {\ref{vectorfields}} holds
for $\alpha\in(0,1]$. Every two points $v$, $w$ of a connected
Carnot manifold can be joined by a rectifiable absolutely
continuous horizontal curve $\gamma$ composed by not more than
countably many segments of integral lines of given horizontal
fields.
\end{thm}

\subsubsection{Comparison of metrics, and Ball--Box Theorem}

\begin{cor}
\label{equiv} Suppose that Assumption {\ref{vectorfields}} holds
for $\alpha\in(0,1]$. In some compact neighborhood the distance
$d_c$  is equivalent to the quasimetric $d_{\infty}$.
\end{cor}

\begin{proof}
An estimate $d_{c}(x,y)\leq C_1d_{\infty}(x,y)$ for points $x,y$
from a compact part $\mathbb M$ follows from Theorem
{\ref{locestcc}}. Our next goal is to prove the converse estimate.
Fix a compact part $K\subset \mathbb M$ and assume the contrary:
for any $n\in \mathbb N$ there exist points $x_n,y_n\in K$ such
that $d_{\infty}(x_n,y_n)\geq nd_{c}(x_n,y_n)$. In this case we
have $d_{\infty}(x_n,y_n)\to0$ as $n\to\infty$ since otherwise we
have simultaneously $d_{c}(x_n,y_n)\to0$ as $n\to\infty$, and
$d_{\infty}(x_n,y_n)\geq \alpha>0$ for all  $n\in \mathbb N$ what
is impossible. We can assume also that $x_n\to x\in K$ as
$n\to\infty$ and $x_n\ne y_n$. Setting
$d_{\infty}(x_n,y_n)=\varepsilon_n$ we have
$d_{\infty}\bigl(x_n,\Delta^{x_n}_{r_0\varepsilon_n^{-1}}y_n\bigr)=r_0$,
and
$d^n_{c}\bigl(x_n,\Delta^{x_n}_{r_0\varepsilon_n^{-1}}y_n\bigr)\leq
r_0n^{-1}$, where the distance $d^n_{c}$ is measured with respect
to the frame $\bigl\{X_i^{\varepsilon_n}\bigr\}$ with
pushed-forward Riemannian tensor. As far as the length of vectors
$X_i^{\varepsilon_n}$, $i=1,\ldots,\dim H_1$, is closed to the
lengths of corresponding nilpotentized vector fields
$\widehat{X}_i$, $i=1,\ldots,\dim H_1$, by Corollary
\ref{GTheorem}, the Riemannian distance
$\rho\bigl(x_n,\Delta^{x_n}_{r_0\varepsilon_n^{-1}}y_n\bigr)\to0$
as $n\to \infty$. It is in a contradiction with
$d_{\infty}\bigl(x_n,\Delta^{x_n}_{r_0\varepsilon_n^{-1}}y_n\bigr)=r_0$
for all $n\in \mathbb N$ (see Proposition \ref{prop1.2} for a
comparison of metrics).
\end{proof}

\begin{rem}\label{equivlit}
Note that, for obtaining the estimate $d_{\infty}(x,y)\leq
C_2d_{c}(x,y)$, the value $\alpha$ need not to be strictly greater
than zero. Thus, the estimate $d_{\infty}(x,y)\leq C_2d_{c}(x,y)$
is valid also for $\alpha=0$.
\end{rem}

Another corollary is so called ball-box theorem proved for smooth
vector fields in \cite{nsw,G}.

\begin{thm}
[Ball--Box Theorem]\label{ballbox} Suppose that Assumption
{\ref{vectorfields}} holds for $\alpha\in(0,1]$. The shape of a
small ball $B(x,r)$ in the metric $d_c$ looks like a box: given
compact set $K\subset\mathbb M$ there are constants $0< C_1\leq
C_2<\infty$ and $r_0$ independent from $x\in K$ such that
\begin{equation}\label{embeddings}
 \operatorname{Box}(x,C_1r)\subset B(x,r)\subset \operatorname{Box}(x,C_2r)
\end{equation}
for all $r\in(0,r_0)$.
\end{thm}

Theorem \ref{ballbox} implies

\begin{cor}\label{dimhaus}
Suppose that Assumption {\ref{vectorfields}} holds for
$\alpha\in(0,1]$. The Hausdorff dimension of  $\mathbb M$ equals
\begin{equation*}
\nu=\sum\limits_{i=1}^Mi(\dim H_i-\dim H_{i-1})
\end{equation*}
where $\dim H_0=0$.
\end{cor}

This Corollary extends Mitchell Theorem \cite{Mit} to
Carnot--Carath\'eodory spaces with minimal smoothness of vector
fields.

\begin{rem}\label{generalization}
Let Assumption {\ref{vectorfields}} holds for $\alpha\in(0,1]$.
Applying Corollary {\ref{ballbox}}, we obtain

\begin{enumerate}
\item the generalization of Theorem {\ref{difference}} for points
$w$ and $w^{\prime}$ close enough:
$$
\max\{d^u_{c}(w,w^{\prime}),
d_{c}(w,w^{\prime})\}={\Theta}[\rho(u,v)
\rho(v,w)]^{\frac{1}{M}}\leq {\Theta}[d_{c}(u,v)
d_{c}(v,w)]^{\frac{1}{M}};
$$

\item the generalization of Theorem {\ref{est_eps}}:
$$
\max\{d_{c}^u(w_{\varepsilon},w^{\prime}_{\varepsilon}),
d_{c}(w_{\varepsilon},w^{\prime}_{\varepsilon})\}=\varepsilon
[\Theta(C,{\mathcal
C})]\rho(u,v)^{\frac{\alpha}{M}}(d_{c}^u(v,w)+o(1));
$$

\item the generalization of Theorem {\ref{mainresult}}:
$$
\max\{d_{c}^u(\widehat{w}^{\varepsilon}_Q,{w^{\varepsilon}_Q}),
d_{c}(\widehat{w}^{\varepsilon}_Q,{w^{\varepsilon}_Q})\}=
{\sum\limits_{k=1}^Q \Theta({\mathcal
C},k,\{F^j_{\mu,\beta}\}_{j,\mu,\beta})}\cdot
\varepsilon^{1+\frac{\alpha}{M}}.
$$
\end{enumerate}
\end{rem}

Corollary {\ref{equiv}} and Theorem 11.11 {\cite{hk}} imply the
following statement containing  a result~of~{\cite{greshn}} where
 only the first assertion is  obtained under
assumption of higher smoothness of vector fields.

\begin{prop}\label{vectfields}  Let $X$ and $Y$ be two families of vector
fields on $\mathbb M$ with the same horizontal distribution
$H\mathbb M$ for both of which  Assumption {\ref{vectorfields}}
holds with some $\alpha\in(0,1]$. Then in some compact
neighborhood the following assertions are equivalent$:$

$1)$ There exists a constant $C\geq 1$ such that
$C^{-1}d_{\infty}^X\leq d_{\infty}^Y\leq Cd_{\infty}^X$.

 $2)$ There exists a constant $C\geq 1$ such that
$C^{-1}|X_H\varphi|\leq |Y_H\varphi|\leq C|X_H\varphi|$ for all
$\varphi\in C^{\infty}(\mathbb M)$.
\end{prop}

Here $d_{\infty}^X$ and $d_{\infty}^Y$ are quasimetrics
constructed with respect to the bases $X$ and $Y$, and
$X_H\varphi$ and $Y_H\varphi$ are subgradients of $\varphi$.

Define the Riemannian quasimetric $d_{\operatorname{riem}}(u,v)$
between a point $u$ and a point $v=\exp
\Bigl(\sum\limits_{i=1}^Nx_iX_i\Bigr)(u)$ as
$d_{\operatorname{riem}}(u,v)=\max\{|x_i|\mid i=1,\dots,N\}$. The
well-known facts of differential geometry  imply that the metric
$d_{\operatorname{riem}}$ is equivalent to the Riemannian metric
$\rho$ on every compactly embedded domain $U\Subset\mathbb M$,
i.e., there exists a constant $c$ independent of the choice of the
points $u$, $v\in U$ and such that $c^{-1}\rho(u,v)\leq
d_{\operatorname{riem}}(u,v)\leq c\rho(u,v)$ for all $u,v\in U$
for which the quantities under consideration are defined. Hence,
we have

\begin{prop}\label{prop1.2} The relations
$$ c^{-1}\rho(u,v)\leq d_{\operatorname{riem}}(u,v)\leq
cd_\infty(u,v) \leq cd_{\operatorname{riem}}(u,v)^{\frac1{M}}
$$
hold for all $x,y\in U$.
 \end{prop}

\begin{rem}
If the derivatives of $X_i$, $i=1,\ldots, N$, are locally
H\"{o}lder with respect to $\mathfrak d$, where $\mathfrak d$
meets conditions of Remark~{\ref{mainrem}}, the statements of
Corollary~{\ref{curves}}, Theorem~{\ref{locestcc}},
Theorem~{\ref{Chow}}, Corollary~{\ref{equiv}},
Theorem~{\ref{ballbox}}, Corollary~{\ref{dimhaus}},
Remark~{\ref{generalization}} and Proposition~{\ref{vectfields}}
are also true.
\end{rem}

\section{Differentiability on a~Carnot
Manifold}\label{differentiabilitycarnot}

\subsection{Primitive calculus}

Further, we extend the dilations $\delta^{g}_{t}$ to negative $t$
by setting $\delta^g_tx=\delta^g_{|t|} (x^{-1})$ for $t<0$. The
convenience of this definition is seen from the comparison of
different kinds of differentiability.

\subsubsection{Definition}  Let $\mathbb M,{\mathbb N}$ be two Carnot
manifolds. We denote the vector fields on $\mathbb N$ by $Y_i$. We
label the remaining objects on $\mathbb N$ (the distance, the
tangent cone etc.) with the same symbols as on~$\mathbb M$ but
with a tilde $\,\,\tilde{}\,\,$ excluding the cases where the
objects under consideration are obvious: for example, for a given
mapping $f:E\to {\mathbb N}$, it is clear that  ${\mathcal G^g
\mathbb M}$ is the tangent cone at a~point $g\in\mathbb M$ and
$\mathcal G^{f(g)}\mathbb N$ is the tangent cone at the point
$f(g)\in\mathbb N$; $d_c^g$ is the metric in the cone  ${\mathcal
G^g} \mathbb M$, $d_c^{f(g)}$ is the metric in $\mathcal G^{f(g)}
\mathbb N$, etc.

Recall that a {\it horizontal homomorphism} of Carnot groups is a
continuous homomorphism $L:\mathbb G\to\widetilde{\mathbb G}$ of
Carnot groups such that

1)  $DL(0)(H\mathbb G)\subset  H\widetilde{\mathbb G}$.

\noindent The notion of a horizontal homomorphism
$L:\bigl({\mathcal G^g\mathbb M},d^g_c\bigr)\to\bigl({\mathcal
G}^q\mathbb N,\tilde d^q_c\bigr)$, $g\in\mathbb M$, $q\in\mathbb
N$, of local Carnot groups is different from this only in that the
inclusion $L({\mathcal G}^g\mathbb M\cap\exp H\mathcal G^g\mathbb
M)\subset {\mathcal G}^q\mathbb N\cap\exp\, H{\mathcal G}^q\mathbb
N$ holds only for $v\in{\mathcal G}^g\mathbb M\cap \exp H\mathcal
G^g\mathbb M$ such that $L(v)\in {\mathcal G}^q\mathbb N$.

Since a homomorphism of Lie groups is continuous, it can be proved
that a horizontal homomorphism $L:\mathbb G\to\widetilde{\mathbb
G}$ also has the property

2) $L(\delta_t v)=\tilde \delta_t L(v)$ for all $v\in\mathbb G$
and $t>0$ (in the case of a horizontal homomorphism
$L:\bigl({\mathcal G}^g\mathbb M,d^g_c\bigr)\to \bigl({\mathcal
G}^q\mathbb N,\tilde d^q_c\bigr)$ of local Carnot groups, the
equality $L(\delta_t v)=\tilde \delta_t L(v)$ is fulfilled only
for $v\in{\mathcal G}^g\mathbb M$ and $t>0$ such that  $\delta_t
v\in{\mathcal G}^g\mathbb M$ and $\tilde \delta_t L(v)\in
{\mathcal G}^q\mathbb N$).

\begin{defn}\label{defdiff}
Given two Carnot manifolds $\mathbb M$ and $\mathbb N$, and a~set
$E\subset\mathbb M$, a mapping $f:E\to {\mathbb N}$ is called {\it
$hc$-differentiable} at a point $g\in E$ if there exists a
horizontal homomorphism
 $L:\bigl(\mathcal G^g\mathbb M, d^g_c\bigr)\to \bigl(\mathcal G^{f(g)}\mathbb N,  d_c^{f(g)}\bigr)$
of the nilpotent tangent cones such that
\begin{equation}\label{estdiff}
  d_c^{f(g)}(f(v),L(v))=o\bigl(d^g_c(g,v)\bigr)
  \quad\text{as  $E\cap {\mathcal G^g\mathbb M}\ni v\to g$}.
\end{equation}
\end{defn}

A horizontal homomorphism $L:\bigl({\mathcal G}^g\mathbb
M,d^g_c\bigr)\to \bigl({\mathcal G}^{f(g)}\mathbb
N,{d}^{f(g)}_c\bigr)$ satisfying condition~\eqref{estdiff}, is
called a $hc$-{\it differential\/} of the mapping $f:E\to{\mathbb
N}$ at $g\in E$ on $E$ and is denoted by  $Df(g)$.
 It can be proved \cite{vod3} that if $g$ is a density point
 then the $hc$-differential is unique.

Moreover, it is easy to verify that the  $hc$-{differential\/}
commutes with the one-parameter dilation group:
\begin{equation}\label{commut}
\delta^{f(g)}_t\circ Df(g)=Df(g)\circ \delta^g_t.
\end{equation}

 \begin{prop}[\cite{vod3}]\label{Prop2.1}
 Definition~{\ref{defdiff}}  is equivalent to each
of the following assertion:

\begin{enumerate}
\item $ d_c^{f(g)}\bigl(\Delta^{f(g)}_{t^{-1}}
f\bigl(\delta^g_{t}(v)\bigr),L(v)\bigr)= o(1)$ as $t\to0$, where
$o(\cdot)$ is uniform in the~points $v$ of any compact part of
$\mathcal G^g\mathbb M;$

\item $\tilde d_{\infty} (f(v),L(v))=o\bigl(d^g_c(g,v)\bigr)$
  as  $E\cap {\mathcal G^g\mathbb M}\ni v\to g;$

\item $\tilde d_{\infty}(f(v),L(v))=o\bigl(d_{\infty}(g,v)\bigr)$
  as  $E\cap {\mathcal G^g\mathbb M}\ni v\to g;$

\item  $ \tilde d_{\infty}\bigl(\Delta^{f(g)}_{t^{-1}}
f\bigl(\delta^g_{t}(v)\bigr),L(v)\bigr)= o(1)$ as $t\to0$, where
$o(\cdot)$ is uniform in the points $v$ of any compact part of
$\mathcal G^g\mathbb M$.
\end{enumerate}
\end{prop}

\begin{proof}
 Consider a point $v$ of a compact part
of  $\mathcal G^g\mathbb M$ and a sequence $\varepsilon_i\to0$ as
$i\to0$ such that $\delta^g_{\varepsilon_i}v\in E$ for all  $i\in
\mathbb N$. From~\eqref{estdiff} we have $
d_c^{f(g)}\bigl(f\bigl(\delta^g_{\varepsilon_i}v\bigr),
 L\bigl(\delta^g_{\varepsilon_i}v\bigr)\bigr)=
 o\bigl(d^g_c\bigl(g,\delta^g_{\varepsilon_i}v\bigr)\bigr)
= o({\varepsilon_i})$. In view of~\eqref{commut}, we infer
 $$d_c^{f(g)}\bigl(\Delta^{f(g)}_{\varepsilon_i}
\bigl(\Delta^{f(g)}_{\varepsilon_i^{-1}}f
\bigl(\delta^g_{\varepsilon_i}v\bigr)\bigr),\delta^{f(g)}_{\varepsilon_i}L(v)
\bigr)=o(\varepsilon_i)\quad \text{uniformly in~$v$.}$$ From here
we obtain item~1. Obviously, the argument is reversible. Item~1 is
equivalent to item~4 since $\alpha\rho(x,y)\leq
d_{\infty}(x,y)\leq \beta\rho(x,y)^\frac{1}{M}$
($\alpha\rho(x,y)\leq d^g_{\infty}(x,y)\leq
\beta\rho(x,y)^\frac{1}{M}$) on any compact part of $\mathbb
M\cap\mathcal G^g\mathbb M$ ($\alpha$ and $\beta$ depend on the
choice of the compact part).

By comparing the metrics of Subsection~{\ref{applications}}:
$d^{f(g)}_\infty(g,v)=O\bigl( d_c^{f(g)}(g,v) \bigr)$, and by
Local Approximation Theorem  \ref{apth_d_1},
 we obtain the equivalence of~~\eqref{estdiff} to the item~2.
The equivalence of items~2 and~3 is obtained by comparing the
metrics of Subsection~{\ref{applications}}:
$d^g_c(g,v)=O\bigl(d^g_\infty(g,v)\bigr)$ and
$d^g_\infty(g,v)=O\bigl(d_{\infty}(g,v)\bigr)$ as $v\to g$.
\end{proof}

\subsubsection{Chain rule} In this subsubsection, we prove the chain rule.

\begin{thm}[\cite{vod3}]\label{chainrule} Suppose that $\mathbb M, {\mathbb N},\mathbb X$
are Carnot manifolds, $E$ is a set in $\mathbb M$, and
$f:E\to{\mathbb N}$ is a mapping from $E$ into $\mathbb N$
$hc$-differentiable at a point $g\in E$. Suppose also that $F$ is
a set in $\mathbb N$, $f(E)\subset F$ and $\varphi:F\to \mathbb X$
is a mapping from $F$ into $\mathbb X$\,\, $hc$-differentiable at
$p=f(g)\in{\mathbb N}$. Then the composition $\varphi\circ
f:E\to{\mathbb X}$ is $hc$-differentiable at $g$ and
$$
D(\varphi\circ f)(g)=D\varphi(p)\circ Df(g).
$$
\end{thm}

\begin{proof} By hypothesis,
$d^{f(g)}_c(f(v),Df(g)(v))=o\bigl(d^g_c(g,v)\bigr)$ as $v\to g$
and also $d^{\varphi(p)}_c(\varphi(w),D\varphi(p)(w))=o\bigl(
d^p_c(p,w)\bigr)$
 as $w\to p$. We now infer
  \begin{multline*}
d^{\varphi(p)}_c((\varphi\circ f)(v),(D\varphi(p)\circ Df(g))(v))\\
\leq  d^{\varphi(p)}_c(\varphi(f(v)),D\varphi(p)(f(v))) +
d^{\varphi(p)}_c(D\varphi(p)(f(v)),D\varphi(p)(Df(g)(v)))\\
\leq o\bigl( d^p_c(p,f(v))\bigr)+O\bigl(
d^p_c\bigl(f(v),Df(g)(v)\bigr)\bigr)
\\
\leq o\bigl(d^g_c(g,v)\bigr)+O\bigl(o\bigl(d^g_c(g,v)\bigr)\bigr)=
o\bigl(d^g_c(g,v)\bigr)\quad\text{as $v\to g$},
\end{multline*}
since
  \begin{multline*}
d^p_c\bigl(p,f(v)\bigr)\leq d^p_c\bigl(p,Df(g)(v)\bigr)+
d^p_c\bigl(f(v),Df(g)(v)\bigr)
\\ =O\bigl(d^g_c(g,v)\bigr)+o\bigl(d^g_c(g,v)\bigr)=O\bigl(d^g_c(g,v)\bigr)
\quad \text{as } v\to g.
\end{multline*}
(The estimate
$d^p_c\bigl(p,Df(g)(v)\bigr)=O\bigl(d^g_c(g,v)\bigr)$ as $v\to g$
follows from the continuity of the homomorphism $Df(g)$
and~\eqref{commut}.)
\end{proof}

\subsection{$hc$-Differentiability of Curves
 on the Carnot Manifolds}

\subsubsection{Coordinate $hc$-differentiability
criterion}   Recall that a mapping $\gamma:E\to {\mathbb M}$,
where $E\subset \mathbb R$ is an arbitrary set, is called a {\it
Lipschitz mapping} if there exists a constant $L$ such that the
inequality $d_{\infty}(\gamma (y),\gamma (x))\leq L|y-x|$ holds
for all $x,y\in E$.

\begin{defn}\label{def3.1}
A mapping  $\gamma:E\to{\mathbb M}$, where $E\subset\mathbb R$ is
an arbitrary set, is called {\it $hc$-differentiable at a limit
point $s\in E$} to $E$ if there exists a horizontal vector
 $a=\sum\limits_{i=1}^{\dim H_1}\alpha_i\widehat
X^{\gamma(s)}_i(\gamma(s))\in H_{\gamma(s)}\mathbb M$ such that
the local homomorphism
$\tau\mapsto\exp\Bigl(\tau\sum\limits_{i=1}^{\dim
H_1}\alpha_i\widehat X^{\gamma(s)}_i\Bigr)(\gamma(s))\in\mathcal
G^{\gamma(s)}\mathbb M$ as the $hc$-differential of the mapping
$\gamma:E{\to}{\mathbb M}$:
  $d^{\gamma(s)}_c\bigl(\gamma(s+\tau),
\delta^{\gamma(s)}_\tau a\bigr) =o(\tau)$ for $\tau\to0$,
$s+\tau\in E$. The point $\exp\Bigl(\sum\limits_{i=1}^{\dim
H_1}\alpha_i\widehat X^{\gamma(s)}_i\Bigr)(\gamma(s))\in\mathcal
G^{\gamma(s)}\mathbb M$ is called the $hc$-{\it
derivative}\,\footnote{If the $hc$-derivative does not exist in
$\mathcal G^{\gamma(s)}\mathbb M$ then it belongs  in $\mathbb
G_{\gamma(s)}\mathbb M$: we consider the ``preimage'' under
$\theta_{\gamma(s)}$ being equal
$\exp\Bigl(\sum\limits_{i=1}^{\dim H_1}\alpha_i(\widehat
X^{\gamma(s)}_i)^{\prime}\Bigr)(0)$ in all the necessary cases.}.
\end{defn}

Some properties of the introduced notion of $hc$-differentiability
can be obtained from Proposition~\ref{Prop2.1}.
 For instance, the coefficients $\alpha_i$
are defined uniquely: if, in the normal coordinates,
$\gamma(s+\tau)=\exp\Bigl(\sum\limits_{i=1}^{N}\gamma_i(\tau)
\widehat X^{\gamma(s)}_i\Bigr)(\gamma(s))$, $s+\tau\in E$, for
sufficiently small $\tau$ then Proposition~\ref{Prop2.1} implies:

\begin{property}[\cite{vod3}]\label{property3.1}
A mapping $\gamma:[a,b]\to {\mathbb M}$ is $hc$-differentiable at
a point $s\in(a,b)$ if and only if one of the following assertions
holds$:$

$(1)$ $\gamma_i(\tau)=\alpha_i\tau+o(\tau)$, $i=1,\dots,\dim H_1$,
 and $\gamma_i(\tau)=o(\tau^{\deg X_i})$, $i>\dim H_1$, as $\tau\to0$,
$s+\tau\in E$$;$

$(2)$ the vector $\sum\limits_{i=1}^{\dim H_1}\alpha_i\widehat
X^{\gamma(s)}_i(\gamma(s))\in H_{\gamma(s)}\mathbb M$ is the
Riemannian derivative of $\gamma:[a,b]\to {\mathbb M}$ at a point
$s\in(a,b)$, and $\gamma_i(\tau)=o(\tau^{\deg X_i})$, $i>\dim
H_1$, as $\tau\to0$, $s+\tau\in E$.
\end{property}

\subsubsection{$hc$-Differentiability of absolutely
continuous curves} If a curve  $\gamma:[a,b]\to\mathbb M$ is
absolutely continuous in the Riemannian sense then all coordinate
functions $\gamma_i(t)$ are absolutely continuous on the closed
interval $[a,b]$ (it is clear that this property is independent of
the choice of the coordinate system). Therefore the tangent vector
$\dot\gamma(t)$ is defined
 almost everywhere on $[a,b]$. If, moreover,
 $\dot\gamma(t)\in H_{\gamma(t)}\mathbb M$ at the points  $t\in[a,b]$
of Riemannian differentiability then the curve
$\gamma:[a,b]\to\mathbb M$ is called {\it horizontal}.

It is well known that almost all points $t$ of a closed interval
$E=[a,b]$ are Lebesgue points of the derivatives of the horizontal
components, that is, if, in the normal coordinates
$\gamma(t+\tau)=\exp\Bigl(\sum\limits_{j=1}^N\gamma_j(\tau)X_j\Bigr)
(\gamma(t))$ then the {\it horizontal\/} components
$\gamma_j(\sigma)$, $j=1,\dots,\dim H_1$, have the property
\begin{equation}\label{3.1}
 \int\limits_{\{\sigma\in
(\alpha,\beta)\,\mid \,t+\sigma\in
E\}}|\dot\gamma_j(\sigma)-\dot\gamma_j(0)|\,d\sigma=
o(\beta-\alpha)\quad\text{as $\beta-\alpha\to0$}
\end{equation}
on intervals $(\alpha,\beta)\ni 0$. Note that property \ref{3.1}
is independent of the choice of the coordinate system in a
neighborhood of $\gamma(t)$.

\begin{thm}[\cite{vod3}]\label{diffcurve}
Let a curve $\gamma:[a,b]\to\mathbb M$ on a Carnot manifold be
absolutely continuous in the Riemannian sense and horizontal. Then
$\gamma:[a,b]\to\mathbb M$ is $hc$-differentiable almost
everywhere$:$ any point $t\in[a,b]$ which is a Lebesgue point of
the derivatives of its horizontal components is also a point at
which $\gamma$ is $hc$-differentiable. If
$\gamma(t+\tau)=\exp\big(\sum\limits_{j=1}^N\gamma_j(\tau)X_j\big)
(\gamma(t))$ then $hc$-derivative $\dot\gamma(t)$ equals
$$
\exp\biggl(\sum\limits_{j=1}^{\dim H_1}\dot\gamma_j(0)\widehat
X^{\gamma(t)}_j\biggr)(\gamma(t))=\exp\biggl(\sum\limits_{j=1}^{\dim
H_1}\dot\gamma_j(0) X_j\biggr)(\gamma(t)).
$$
\end{thm}

\begin{proof}
Fix a Lebesgue point $t_0\in(a,b)$ of the derivatives of the
horizontal components of the mapping
$\gamma(t_0+\tau)=\exp\Bigl(\sum\limits_{j=1}^N\gamma_j(\tau)X_j\Bigr
) (g)$, $g=\gamma(t_0)$. In this proof, we also fix a normal
coordinate system $\theta_g$ at $g$. To simplify the notation, we
write the vector fields
 $\widetilde{X}^g_i=(\theta^{-1}_g)_*{X}_i$ and
$\widehat{X}'_i{}^g=(\theta^{-1}_g)_*\widehat{X}^g_i$ defined in a
neighborhood of $0\in\mathbb R^N$ without the superscript $g$:
$\widetilde{X}_i=(\theta^{-1}_g)_*{X}_i$ and
$\widehat{X}'_i=(\theta^{-1}_g)_*\widehat{X}^g_i$ respectively.

For proving the $hc$-differentiability of the mapping  $\gamma$ at
$t_0$, we need to establish the estimate
$\gamma_j(\tau)=o(\tau^{\deg X_j})$ as $\tau\to0$ for all $j>\dim
H_1$, $t_0+\tau\in[a,b]$ (see Property~\ref{property3.1}).
Partition the proof of the desired estimate into several steps.

\smallskip

$1^{\text{\sc st}}$ {\sc Step.}  Here we show that the hypothesis
implies the Riemannian differentiability of the mapping $\gamma$
at $t_0$ and $\dot\gamma(t_0)\in H_g\mathbb M$. Put
$\Gamma(t_0+\tau)=\theta^{-1}_g(\gamma(t))=(\gamma_1(\tau),\dots,\gamma_N(\tau))$.
The curve $\Gamma(\tau)$ is absolutely continuous, and its tangent
vector $\dot\Gamma(\tau)$ is horizontal in  a neighborhood of
$0\in T_g\mathbb M$ with respect to the vector fields
$\{\widetilde{X}_i\}$: $\dot\Gamma(\tau)\in
(\theta^{-1}_g)_*H_{\gamma(t_0+\tau)}\mathbb M$ for almost
all~$\tau$. From here, for almost all~$\tau$ sufficiently closed
to~$0$, we infer
\begin{equation}\label{3.2}
\dot\Gamma(\tau)=\sum\limits_{j=1}^N\dot\gamma_j(\tau)\frac{\partial}{\partial
x_j}= \sum\limits_{i=1}^{\dim H_1} a_i(\tau)\widetilde
X_i(\Gamma(\tau)).
\end{equation}

The Riemann tensor pulled back from the manifold $\mathbb M$ onto
a neighborhood of $0\in T_g\mathbb M$ is continuous at the zero.
Therefore, using this continuity, we see that, for any $\tau$,
$t_0+\tau\in[a,b]$, \eqref{3.1} implies
\begin{multline*}
d_c(\gamma (t_0),\gamma (t_0+\tau)) \leq c_1\int\limits_{(0,\tau)}
|\dot\Gamma(\sigma)|_r \, d\sigma\\ \leq
c_2\sum\limits_{j=1}^{\dim H_1}\int\limits_{(0,\tau)}
(|\dot\gamma_j(\sigma)-\dot\gamma_j(0)| + |\dot\gamma_j(0)|)\,
d\sigma=O(\tau)
\end{multline*}
as $\tau\to0$, where $|\dot\Gamma(\sigma)|_r$ stands for the
length of the tangent vector in the pulled-back Riemannian metric.
By Proposition~\ref{equiv} and Remark~\ref{equivlit}, we have
$d_\infty(\gamma (t_0),\gamma (t_0+\tau))=O(d_c(\gamma
(t_0),\gamma (t_0+\tau))$ as $\tau\to0$. Therefore the coordinate
components $\gamma_j(\tau)$ of the mapping $\gamma$ satisfy
\begin{equation}\label{3.3}
\gamma_j(\tau)=O(\tau^{\deg X_j})\quad\text{as $\tau\to0$ for all
$j\geq1$}.
\end{equation}
It follows that the curve $\Gamma(\tau)$ is differentiable at $0$
and
$$\dot\Gamma(0)=(\dot\gamma_1(0),\dots,\dot\gamma_{\dim
H_1}(0),0,\dots,0).$$ Hence, the curve  $\gamma$ is differentiable
in the Riemannian sense at $t_0$ and  $\dot\gamma(t_0)\in
H_{g}\mathbb M$. From \eqref{3.3} we also obtain $\gamma(\tau)\in
B(g,O(\tau))$.

\smallskip
$2^{\text{\sc nd}}$ {\sc Step.}  Corollary {\ref{lem_matrix}} and
the fact that $\gamma(\tau)\in B(g,O(\tau))$ imply that, in a
neighborhood of~$0$,
 the vector
fields $\widetilde X_i$ can be expressed via $\widehat X'_k$ so
that
$$
\widetilde X_i(\Gamma(\tau))=\sum\limits_{k=1}^N
\alpha_{ik}(\tau)\widehat X_k'(\Gamma(\tau)) ,\ \text{where
$\alpha_{ik}(\tau)=\begin{cases} o(\tau^{\deg X_k-\deg X_i})\quad
\text{if}\\
\hskip30pt\deg X_k>\deg X_i,\\
\delta_{ik}+o(1)\quad \text{otherwise} \end{cases} $}
$$
as $\tau\to0$. Now, using expansion~\eqref{mat_tc} of the vector
fields $\widehat X'_i$ in the standard Euclidean basis, for all
points $\tau$ sufficiently close to $0$, from~\eqref{3.2} we now
obtain
\begin{multline}\label{3.4}
\sum\limits_{j=1}^N\dot\gamma_j(\tau)\frac{\partial}{\partial
x_j}= \sum\limits_{i=1}^{\dim H_1} a_i(\tau)\widetilde
X_i(\Gamma(\tau)) =\sum\limits_{k=1}^N \sum\limits_{i=1}^{\dim
H_1} a_i(\tau)\alpha_{ik}(\tau)\widehat X_k'(\Gamma(\tau))
\\
=\sum\limits_{j=1}^N\sum\limits_{k=1}^j \sum\limits_{i=1}^{\dim
H_1} a_i(\tau)\alpha_{ik}(\tau)\hat z^j_k(\Gamma(\tau))
\frac{\partial}{\partial x_j}.
\end{multline}

\smallskip
$3^{\text{\sc rd}}$ {\sc Step.}  For $1\leq j\leq\dim H_1$, we
have $\deg X_j=1$. Then from \eqref{3.3} and~\eqref{mat_tc} we
conclude that $\hat z^j_k(\Gamma(\tau))=\delta_{jk}+O(\tau)$.
Therefore, from \eqref{3.4} we infer
 $\dot\gamma_j(\tau)=\sum\limits_{i=1}^{\dim H_1}
a_i(\tau)\tilde\alpha_{ij}(\tau)$, where, as before,
$\tilde\alpha_{ij}(\tau)=\delta_{ij}+o(1)$. Hence,
$a_i(\tau)=\sum\limits_{n=1}^{\dim H_1}
\dot\gamma_n(\tau)\beta_{ni}(\tau)$, where
$\{\beta_{ni}(\tau)\}$, $n,i=1,\dots,\dim H_1$, is a matrix
inverse to $\{\tilde\alpha_{ij}(\tau)\}$, has the elements
$\beta_{ni}(\tau)=\delta_{ni}+o(1)$. Consequently,
\begin{multline}\label{3.04}
a_i(\tau)=\sum\limits_{i=1}^{\dim H_1}
\dot\gamma_n(\tau)\beta_{ni}(\tau)=\sum\limits_{n=1}^{\dim H_1}
\dot r_n(\tau)\beta_{ni}(\tau)+ \sum\limits_{n=1}^{\dim H_1}
\dot\gamma_n(0)\beta_{ni}(\tau), \\
\text{where
$r_n(\tau)=\int\limits_0^\tau(\dot\gamma_n(\sigma)-\dot\gamma_n(0))\,d\sigma$}.
\end{multline}

\smallskip
$4^{\text{\sc th}}$ {\sc Step.} Fix $\dim H_{l-1}<j\leq\dim
H_{l}$, $1<l\leq M$. For estimating
 $\dot\gamma_j(\tau)$, from \eqref{3.4} we have
 \begin{multline}\label{3.05}
\dot\gamma_j(\tau)=\sum\limits_{k,i,n=1}^{\dim H_1}
\dot\gamma_n(\tau)\beta_{ni}(\tau)\alpha_{ik}(\tau)\hat
z^j_k(\Gamma(\tau))
\\+
\sum\limits_{k=\dim H_1+1}^{j} \sum\limits_{i,n=1}^{\dim H_1}
\dot\gamma_n(\tau)\beta_{ni}(\tau)\alpha_{ik}(\tau)\hat
z^j_k(\Gamma(\tau))= I_j+II_j.
\end{multline}
Since in this case $\alpha_{ik}(\tau)=o(\tau^{\deg X_k-\deg
X_i})$, and $\hat z^j_k(\Gamma(\tau))= O(\tau^{\deg X_j-\deg
X_k})$ by~\eqref{3.3} then all components in the double sum in
\eqref{3.05} have a factor $o(\tau^{l-1})$. Therefore
\begin{equation}\label{3.06}
II_j= \sum\limits_{n=1}^{\dim H_1}
\dot\gamma_n(\tau)o(\tau^{l-1}).
\end{equation}
From another side
\begin{multline}\label{3.6}
I_j= \sum\limits_{n=1}^{\dim H_1} \dot\gamma_n(\tau)\hat
z^j_n(\Gamma(\tau))+ \sum\limits_{k,n=1}^{\dim H_1}
\dot\gamma_n(\tau)o(1)\hat z^j_k(\Gamma(\tau))\\=
\sum\limits_{n=1}^{\dim H_1} \dot \gamma_n(0)\hat
z^j_n(\Gamma(\tau))+ \sum\limits_{n=1}^{\dim H_1} \dot
r_n(\tau)\hat z^j_n(\Gamma(\tau))+ \sum\limits_{k,n=1}^{\dim H_1}
\dot\gamma_n(\tau)o(1)\hat z^j_k(\Gamma(\tau))\\=
\sum\limits_{n=1}^{\dim H_1} \dot \gamma_n(0)\sum\limits_{
|\alpha+e_n|_h=\deg X_j,\,\alpha>0}
F^j_{\alpha,e_n}(g)\Gamma(\tau)^\alpha
\\+
\sum\limits_{n=1}^{\dim H_1} \dot r_n(\tau)O(\tau^{l-1})+
\sum\limits_{n=1}^{\dim H_1} \dot\gamma_n(\tau)o(\tau^{l-1}).
\end{multline}
In the estimation of the increment of  $\gamma_j(\tau)$ on
$[0,\tau]$ by the Newton--Leibnitz formula, the components of
\eqref{3.06} and the last two summands in \eqref{3.6} have order
$o(\tau^{l})$.
 Indeed,
for all $1\leq i\leq\dim H_1$ and $s>0$, from \eqref{3.1} and
\eqref{3.04} we have $|\dot\gamma_n(\tau)|\leq
|\dot\gamma_n(0)|+|\dot r_i(\tau)|$ from \eqref{3.04},
 $
 |r_i(\tau)|\leq\int\limits_0^\tau
|\dot\gamma_i(\sigma)-\dot\gamma_i(0)|\,d\sigma=o(\tau)$ and
$$
 \biggl|\int\limits_0^\tau\dot
 r_i(\sigma)O(\sigma^{s})\,d\sigma\biggr|\leq
|O(\tau^{s})|\int\limits_0^\tau
|\dot\gamma_i(\sigma)-\dot\gamma_i(0)|\,d\sigma=
 o(\tau^{s+1}).
 $$

$5^{\text{\sc th}}$ {\sc Step.} In the remaining double sum in
\eqref{3.6}, the summands with index $\alpha$ for which
$|\alpha+e_n|<\deg X_j$ contain the factor
 $\Gamma(\tau)^\alpha=o(\tau^{l-1})$, since, in this case, the product
$\Gamma(\tau)^\alpha$ necessarily contains the factor
$\gamma_j(\tau)=\dot\gamma_j(0)\tau+o(\tau)=o(\tau)$, $j>\dim
H_1$. Therefore, expression \eqref{3.6} for $\dot\gamma_j(\tau)$
is reduced to the following:
\begin{equation}\label{3.8}
 \dot\gamma_j(\tau)=\sum\limits_{i=1}^{\dim H_1}
\dot\gamma_i(0)\sum\limits_ {\substack{|\alpha+e_n|_h=\deg
X_j,\\|\alpha+e_n|=\deg X_j}}
F^j_{\alpha,e_n}(g)\Gamma(\tau)^\alpha +o(\tau^{l-1}).
\end{equation}
Since also $\Gamma(\tau)=\dot\Gamma(0)\tau+o(\tau)$, we see that
each summand in \eqref{3.8} is equal to $\dot\gamma_i(0)
F^j_{\alpha,e_n}(g)\Gamma(\tau)^\alpha= \tau^{l-1}\dot\gamma_i(0)
F^j_{\alpha,e_n}(g)\dot\Gamma(0)^\alpha+ o(\tau^{l-1})$.
Consequently, leaving only the summands of order $\tau^{l-1}$ in
\eqref{3.8}, we have
\begin{equation}\label{3.9}
\dot\gamma_j(\tau)=\sum\limits_{i=1}^{\dim H_1} \tau^{l-1}
\sum\limits_{|\alpha|=|\alpha|_h=l-1} \dot\gamma_i(0)
F^j_{\alpha,e_n}(g)\dot\Gamma(0)^\alpha+ o(\tau^{l-1}).
 \end{equation}
 Similarly, the second summand in the estimation of the increment
$\gamma_j(\tau)$ is equal to $o(\tau^{l})$. Consequently, for the
validity of the theorem, it is necessary and sufficient that the
double sum in \eqref{3.9} be zero. This was established in
Lemma~\ref{property_sum}.

Thus, we have proved that $\gamma_j(\tau)= o(\tau^{\deg X_j})$ for
all $j>\dim H_1$. Since the horizontal components of~$\gamma$ are
differentiable at $t_0$, by Property~3.1, the estimate
$\gamma_j(\tau)= o(\tau^{\deg X_j})$ for all $j>\dim H_1$ yields
the $hc$-differentiability of $\gamma$ at~$t_0$.
\end{proof}

The method of proving Theorem~\ref{diffcurve} is applicable to a
wider class of mappings and makes it possible to make additional
conclusions about the nature of $hc$-differentiability.

\begin{cor}\label{corollary3.1} Suppose that a curve $\gamma:[a,b]\to\mathbb M$
on a Carnot manifold is Lipschitz with respect to the Riemannian
metric and horizontal, i.e., $\dot\gamma(s)\in
H_{\gamma(s)}\mathbb M$ for almost every $s\in[a,b]$. Then the
curve $\gamma:[a,b]\to\mathbb M$ is $hc$-differentiable almost
everywhere\footnote{\footnotesize In papers {\cite{vod3, V3}}, a
wrong Corollary 3.1 is formulates instead of this.}.
\end{cor}

\begin{proof}
Every Lipschitz curve with respect to the Riemannian metric is
also  absolutely continuous in the Riemannian sense. Thus all
conditions of Theorem~\ref{diffcurve} hold.
\end{proof}

\begin{cor}\label{corollary3.2} Suppose that we have a family of curves
$\gamma:[a,b]\times F\to\mathbb M$ on a Carnot manifold $\mathbb
M$ that is bounded and continuous in the totality of its
variables, where $F$ is a locally compact metric space. Suppose
that, for each fixed $u\in F$, the curve $\gamma(\cdot,u)$ is
differentiable in the Riemannian sense at all points of $[a,b]$
and horizontal, i.e., $\frac{d}{ds}\gamma(s,u)\in
H_{\gamma(s,u)}\mathbb M$ for all $s\in[a,b]$. If the Riemannian
derivative $\frac{d}{ds}\gamma(s,u)$ is bounded and  continuous in
the totality of its variables $s$ and $u$ then its $hc$-derivative
is also bounded and continuous on $[a,b]\times F$. Furthermore,
the convergence $\Delta^{\gamma(s)}_{\tau^{-1}}\gamma(s+\tau,u)$
to $\dot\gamma(s,u)\in \mathcal G^{\gamma(s,u)}\mathbb M$ is
locally uniform in the totality of $s\in[a,b]$ and $u\in F$.
 \end{cor}

\begin{proof} It suffices to prove in all items of the proof of Theorem~\ref{diffcurve} that the smallness
of all quantities converging to zero is locally uniform on
$[a,b]\times F$ (see Proposition~\ref{equiv} for the estimate
$C_0d_\infty(g,v)\leq d_c(g,v)$)
\end{proof}

\begin{cor}\label{corollary3.3} Suppose that a curve $\gamma:[a,b]\to\mathbb M$ on
a Carnot manifold belongs to $C^1$ and its Riemannian tangent
vector $\dot\gamma_i(t)$ is horizontal  for all $t\in[a,b]$. Then
the curve $\gamma:[a,b]\to\mathbb M$ is $hc$-differentiable at all
$t\in[a,b]$. Furthermore, the convergence of
$\Delta^{\gamma(s)}_{\tau^{-1}}\gamma(s+\tau)$ to
$\dot\gamma(s)\in \mathcal G^{\gamma(s)}\mathbb M$ is uniform in
$s\in[a,b]$.
\end{cor}

\begin{proof} For any $x,y\in[a,b]$, the length $L(\gamma\vert_{[x,y]})$
of the curve $\gamma:[x,y]\to\mathbb M$ is defined; moreover,
$d_{\infty}(\gamma (y),\gamma (x))\leq
c_1L(\gamma\vert_{[x,y]})\leq c_1C|y-x|$, where
$C=\max\limits_{t\in[a,b]}|\dot\gamma(t)|$. Thus, the curve
$\gamma:[a,b]\to\mathbb M$ meets the conditions of
Theorem~\ref{diffcurve} at all points of $[a,b]$ and, therefore,
is uniformly $hc$-differentiable by Corollary~\ref{corollary3.2}.
The last assertion of this corollary follows
\end{proof}

\begin{lem}\label{lemma3.1} Every Lipschitz mapping $\gamma:E\to {\mathbb M}$
is differentiable almost everywhere in the Riemannian sense, and
$\dot \gamma(t)\in H_{\gamma(t)}\mathbb M$ at the points of the
Riemannian differentiability of $\gamma$
\end{lem}

\begin{proof}  In the normal coordinates at a point $g=\gamma(t)$, we have
$$
\gamma(t+\tau)=\exp\biggl(\sum\limits_{j=1}^N\gamma_j(\tau)\,X_j\biggr)(g),\quad
t+\tau\in E.
$$
The Lipschitzity of the mapping $\gamma:E\to {\mathbb M}$ and the
properties of $d_\infty$  imply the estimate
$\gamma_j(\tau)=O(\tau^{\deg X_j})$ for all $j\geq1$, $t+\tau\in
E$. Since $\deg X_j\geq2$ for $j>\dim H_1$, the derivative
$\dot\gamma_j(0)$ exists and is zero for all $j>\dim H_1$.
Consequently, the Riemannian differentiability of $\gamma$ at $t$
is equivalent to the differentiability of the horizontal
components $\gamma_j$, $j=1,\dots,n$, of $\gamma$ at~$0$.

Now, the Lipschitz mapping $\gamma:E\to {\mathbb M}$ is also
Lipschitz with respect to the Riemannian metric (see
Proposition~\ref{prop1.2}). Thus, by Rademacher's classical
theorem, the Riemannian derivative $\dot\gamma(t)\in
T_{\gamma(t)}{\mathbb M}$ exists for almost every $t\in [a,b]$.
The above implies that, at every such point, $\dot\gamma(t)\in
H_{\gamma(t)}\mathbb M$
\end{proof}

Since a Lipschitz mapping $\gamma:[a,b]\to {\mathbb M}$ is
absolutely continuous in the Riemannian sense (see the comparison
of the metrics in Proposition~\ref{prop1.2}), from
Lemma~\ref{lemma3.1} and Theorem~\ref{diffcurve} we infer

\begin{cor}\label{corollary3.4} Every Lipschitz mapping
$\gamma:[a,b]\to {\mathbb M}$ is $hc$-differentiable almost
everywhere on $[a,b]$: if $t\in[a,b]$ is a Lebesgue point of the
derivatives of its horizontal components then this point is its
$hc$-differentiability point.
\end{cor}

\subsubsection{$hc$-Differentiability of scalar Lipschitz mappings} In this subsubsection, we establish the
$hc$-differentiability of the Lipschitz mappings $\gamma:E\to
{\mathbb M}$ where $E\subset \mathbb R$ is an arbitrary set.

Recall that $x\in A$, where $A\subset\mathbb R$ is a measurable
set, is the density point of $A$ if
 $$
 |A\cap(\alpha,\beta)|_1=\beta-\alpha+o(\beta-\alpha)\quad\text
 {for $\beta-\alpha\to0$, $x\in(\alpha,\beta)$}
 $$
(here $|\cdot|_1$ stands for the one-dimensional Lebesgue
measure). It is known that almost all points of a measurable set
$A$ are its density points  (for example, see~\cite{F}).

It is explicitly seen from the above proof of Lemma~\ref{lemma3.1}
that the question of $hc$-differentiability for a Lipschitz
mapping depends on the differentiability of its horizontal
components. If a Lipschitz mapping $\gamma:E\to {\mathbb M}$ (we
may assume that $E\subset \mathbb\mathbb R$ is closed) is written
in the normal coordinates:
$\gamma(t+\tau)=\exp\big(\sum\limits_{j=1}^N\gamma_j(\tau)\,X_j\big)(\gamma(t))$,
$t\in E$ is a fixed number, $t+\tau\in E$, then, by
Lemma~\ref{lemma3.1}, its components $\gamma_j(\tau)$,
$j=1,\dots,N$, are differentiable almost everywhere on $E$. It is
known that almost all density points of $E$ are Lebesgue points of
the derivative of the horizontal components, i.e., for intervals
$(\alpha,\beta)\ni \tau$, $t+\tau\in E$, we infer
\begin{equation}\label{3.10}
\int\limits_{\{\sigma\in (\alpha,\beta)\,\mid \,t+\sigma\in
E\}}|\dot\gamma_j(\sigma)-\dot\gamma_j(\tau)|\,d\sigma=
o(\beta-\alpha)\quad\text{for $\beta-\alpha\to0$}
\end{equation}
for all $j=1,\dots,\dim H_1$. Note that property \eqref{3.10} does
not depend on the choice of the coordinate system in a
neighborhood of the point $g=\gamma(t)$.

\begin{thm}[\cite{vod3}]\label{difflipcurve}
Every Lipschitz mapping $\gamma:E\to {\mathbb M}$, $E\subset
\mathbb R$ is closed, is $hc$-differentiable almost everywhere on
$E$$:$ the mapping $\gamma:E\to {\mathbb M}$ is
$hc$-differentiable at every point $t\in E$ such that

\begin{enumerate}
\item $t$ is the density point of $E$$;$ \item there exist
derivatives $\dot\gamma_j(0)$, $j=1,\dots,\dim H_1$, of the
horizontal components of $\gamma$, where
$\gamma(t+\tau)=\exp\Bigl(\sum\limits_{j=1}^N\gamma_j(\tau)\,X_j\Big)(\gamma(t))$,
 $t+\tau\in E$$;$
\item condition \eqref{3.10} is fulfilled at the point $\tau=0$.
\end{enumerate}
\noindent The $hc$-derivative $\dot\gamma(t)$  equals
$$\exp\biggl(\sum\limits_{j=1}^{\dim H_1}\dot\gamma_j(0)\widehat
X^{\gamma(t)}_j\biggr)(\gamma(t))
=\exp\biggl(\sum\limits_{j=1}^{\dim
H_1}\dot\gamma_j(0)X_j\biggr)(\gamma(t)).$$
\end{thm}

\begin{proof}  $1^{\text{\sc st}}$ {\sc Step.} Suppose that $t\in E$ is a point at
which conditions~1--3 of the theorem hold and $g=\gamma(t)$. Since
the result is local, we may also assume that $E$ is included in an
interval $[a, b]\subset\mathbb R$, $t\in [a, b]$, $a,b\in E$,
whose image is included in $\mathcal G^g\mathbb M$ (we may assume
by diminishing the interval $[a, b]$ if necessary that $\gamma([a,
b]\cap E)\subset \mathcal G^{\gamma(\eta)}\mathbb M$ for every
$\eta\in[a, b]\cap E$).

The open bounded set $Z=(a, b)\setminus E$ is representable as the
union of an at most countable collection of disjoint intervals:
$Z=\bigcup_j(\alpha_j,\beta_j)$, where, for convenience of the
subsequent estimates, we put $\alpha_j<\beta_j$ if $t\leq\alpha_j$
and $\beta_j<\alpha_j$ if $\alpha_j< t$. It is known (for example,
see~\cite{fs}), that, in $\mathcal G^{\gamma(\alpha_j)}\mathbb M$,
there exists a horizontal curve $\tilde\sigma_j:[0,b_j]\to
\mathcal G^{\gamma(\alpha_j)}\mathbb M$ joining the points
$\tilde\sigma_j(0)=\gamma(\alpha_j)$ and
$\tilde\sigma_j(b_j)=\gamma(\beta_j)$ and parameterized by the arc
length; moreover,
$b_j=d^{\gamma(\alpha_j)}_c(\gamma(\alpha_j),\gamma(\beta_j)) \leq
cd^{\gamma(\alpha_j)}_\infty(\gamma(\alpha_j),\gamma(\beta_j))=
cd_\infty(\gamma(\alpha_j),\gamma(\beta_j))\leq
cL|\beta_j-\alpha_j|$, where $c$ is independent of $j$ (see the
relation between the metrics in Subsection~{\ref{applications}}).
Consequently, the mapping $\sigma_j:[\alpha_j,\beta_j]\to\mathbb
M$ defined by the rule
$$
[\alpha_j,\beta_j]\ni \eta\mapsto
\sigma_j(\eta)=\tilde\sigma_j\Bigl(\frac{b_j}{|\beta_j-\alpha_j|}|\eta-\alpha_j|\Bigr)
\in \mathcal G^{\gamma(\alpha_j)}\mathbb M
$$
is Lipschitz in the metric $d^{\gamma(\alpha_j)}_c$ with the
Lipschitz constant $cL$ for all $j\in \mathbb N$. Define now the
extension $f:[a, b]\to\mathbb M$ as follows:
$$
f(\eta)=
\begin{cases}  \gamma(\eta),&  \text{if $\eta\in E$,}\\
\sigma_j(\eta),&\text{if $\eta\in(\alpha_j,\beta_j)$}.
\end{cases}
$$

\smallskip
$2^{\text{\sc nd}}$ {\sc Step.} The mapping $f:[a, b]\to\mathbb M$
has the following properties:

(1) $f:[\alpha, \beta]\to\mathbb M$ is a Lipschitz mapping with
respect to the Riemannian metric;

(2) the Riemannian derivative of $f$ exists for almost every
$\eta\in[a, b]$ and is bounded;

(3) the vector $\dot f(\eta)$ belongs to the horizontal space
$H_{\gamma(\eta)}\mathbb M$ for almost every $\eta\in E$;

$(4)$ the mapping $f:[a, b]\to\mathbb M$ has a Riemannian
derivative at $t$ equal to $\dot\gamma(t)$;

\noindent if
$f(t+\tau)=\exp\Big(\sum\limits_{j=1}^Nf_j(\tau)\,X_j\Big)(g)$,
$t+\tau\in [a, b]$, then

$(5)$ $ f_j(\tau)=O(\tau^{\deg X_j})\quad\text{as $\tau\to0$ for
all $j\geq1$}$;

$(6)$  $0$ is a Lebesgue point for the derivatives $\dot
f_j(\tau)$, $j=1,\dots,\dim H_1$.

Indeed, if  $t\leq
\alpha_j<\eta_1<\beta_j<\alpha_k<\eta_2<\beta_k\leq b$ then,
taking the relations between the metrics into account, we obtain
the estimates $\rho (f(\eta_1),f(\eta_2)) \leq \rho
(f(\eta_1),\gamma(\beta_j))+ \rho
(\gamma(\beta_j),\gamma(\alpha_k))+\rho
(\gamma(\alpha_k),f(\eta_2)) \leq
C((\beta_j-\eta_1)+(\alpha_k-\beta_j)+(\eta_2-\alpha_k))=
C|\eta_2-\eta_1|$. The other cases of mutual disposition of
$\eta_1$ and $\eta_2$ with respect to $t$ are considered
similarly. Hence we obtain properties~(1) and~(2).

Next, if $t\leq \alpha_j<t+\tau<\beta_j$ then
$d_\infty(f(t+\tau),f(t)) \leq
C(d_\infty(f(t+\tau),\gamma(\alpha_j))+
d_\infty(\gamma(\alpha_j),\gamma(t))) \leq
C_1\bigl(d^{\gamma(\alpha_j)}_\infty(f(t+\tau),\gamma(\alpha_j))
+(\alpha_j-t)\bigr)= C_2((t+\tau -\alpha_j)
+(\alpha_j-t))=C_2\tau$  by~the triangle inequality, the
construction of~$f$, and the relations between the metrics. From
this we obtain property~(5) and, hence, the differentiability of
all components $f_j$ at~$0$, $j>\dim H_1$:
  $\dot f_j(0)=0$.

Since the derivatives of Lipschitz functions are bounded and $t$
is the density point of $E$, for intervals $(r,s)\ni 0$ we have
\begin{multline}\label{3.11}
 \int\limits_{(r,s)}|\dot f_j(\sigma)-\dot
\gamma_j(0)|\,d\sigma= \int\limits_{\{\sigma\in(r,s)\,\mid
\,t+\sigma\in E\cap[a,b]\}}|\dot \gamma_j(\sigma)-\dot
\gamma_j(0)|\,d\sigma
\\
+\int\limits_{\{\sigma\in(r,s)\,\mid \,t+\sigma\notin
E\cap[a,b]\}}|\dot f_j(\sigma)-\dot \gamma_j(0)|\,d\sigma =
 o(|s-r|)
\end{multline}
as $s-r\to0$ for all $j=1,\dots,\dim H_1$. Hence,
$\int\limits_0^\tau(\dot f_j(\sigma)-\dot\gamma_j(0))\,d\sigma=
f_j(\tau)-\dot\gamma_j(0)\tau=o(\tau)$ and
$\frac{df_j}{d\tau}(0)=\dot\gamma_j(0)$ for all $j=1,\dots,\dim
H_1$. Thus, we have proved properties~(4) and~(6).

Note that the preceding arguments are independent of the
coordinate system. They are based on the following principle: if
$\eta$ is the density point for  $E$, the mapping $f\vert_E$ has a
Riemannian derivative at $\eta\in E$,
 and $\eta\in E$ is a Lebesgue point
for the horizontal coordinate functions of $f\vert_E$ then, with
regard to Lemma~\ref{lemma3.1} and what has been proved above, $f$
has a Riemannian derivative at $\eta$; moreover, the Riemannian
tangent vector belongs to the horizontal space
$H_{\gamma(\eta)}\mathbb M$. This proves property~(3).

$3^{\text{\sc rd}}$ {\sc Step.} Since the Riemannian derivative
$\dot f(\eta)$ of the mapping $f:[a,b]\to\mathbb M$ belongs to the
horizontal space $H_{f(\eta)}\mathbb M$ only at almost every point
$\eta\in E$, a direct application of Theorem~\ref{diffcurve} is
impossible. However, granted the fact that the complement
$[a,b]\setminus E$ has density zero at $t$, the method of its
proof can be adapted also to this case. We now indicate the
changes to the proof of Theorem~\ref{diffcurve} necessary for
obtaining the $hc$-differentiability of $f$ at the point $t$ fixed
above.

Introduce the notation
\begin{equation*}
\Gamma(\tau)=\begin{cases}(\gamma_1(\tau),\dots,\gamma_N(\tau)),\quad
\text{if $t+\tau\in E$},\\
(f_1(\tau),\dots,f_N(\tau)),\quad \text{if $t+\tau\notin E$.}
\end{cases}
\end{equation*}
It has been proved above that
$\dot\Gamma(0)=(\dot\gamma_1(0),\dots,\dot\gamma_{\dim
H_1}(0),0,\dots,0)$.  Deduce~\eqref{3.2} for the points $\tau$
sufficiently close to~$0$ and such that $t+\tau\in E$. At the
points $t+\tau\in (\alpha_j,\beta_j)$, we have
\begin{equation}\label{3.12}
\dot\Gamma(\tau)=\sum\limits_{j=1}^N\dot
f_j(\tau)\frac{\partial}{\partial x_j}= \sum\limits_{i=1}^{\dim
H_1} a_i(\tau)\widehat X_i'{}^{f(\alpha_j)}(\Gamma(\tau)).
 \end{equation}

By Proposition \ref{estimate}, at the points $t+\tau\in
(\alpha_j,\beta_j)$ the relation $f(\tau)\in B(g,O(\tau))$ implies
that, in a neighborhood of~$0$, the vector fields $\widehat
X_i'{}^{f(\alpha_j)}$ are expressed via the vector fields
$\widehat X'_k$ (here we write $\widehat X'_k$ instead of
$\widehat X_k'{}^g$) in the form
$$
\widehat X_i'{}^{f(\alpha_j)}(\Gamma(\tau))=\sum\limits_{k=1}^N
\gamma_{ik}(\tau)\widehat X_k'(\Gamma(\tau)),\ \text{where
$\gamma_{ik}(\tau)= \begin{cases} o(\tau^{\deg X_k-\deg X_i}),\
\text{if}\\
\hskip30pt\deg X_k>\deg X_i,\\
\delta_{ik}+o(1)\ \text{otherwise}
\end{cases} $}
$$
as $\tau\to0$. Really, by \eqref{expan},   we have $\widehat
X_i'{}^{f(\alpha_j)}(\Gamma(\tau))=\sum\limits_{l=1}^N\beta_{il}(\tau)\widetilde
X_l(\Gamma(\tau))$ at points $f(\tau)\in B(g,O(\tau))$, where
\begin{equation}\label{estimate1}
\beta_{il}(\tau)=\begin{cases}o(\tau^{\deg X_l-\deg X_i})\quad&\text{if $\deg X_l>\deg X_i$},\\
\delta_{il}+o(1)&\text{otherwise}
\end{cases}
\end{equation}
as $\tau\to0$, and $\widetilde
X_l(\Gamma(\tau))=\sum\limits_{k=1}^N \alpha_{lk}(\tau)\widehat
X_k'(\Gamma(\tau))$ where
\begin{equation}\label{estimate2}
\alpha_{lk}(\tau)=\begin{cases} o(\tau^{\deg X_k-\deg X_l})\quad
&\text{if $\deg X_k>\deg X_l$},\\
\delta_{ik}+o(1)\quad &\text{otherwise}
\end{cases}
\end{equation}
as $\tau\to0$. It follows $\widehat
X_i'{}^{f(\alpha_j)}(\Gamma(\tau))=
\sum\limits_{k=1}^N\sum\limits_{l=1}^N\beta_{il}(\tau)\alpha_{lk}(\tau)
\widehat X_k'(\Gamma(\tau))$. Now taking into account
\eqref{estimate1} and \eqref{estimate2}, and representing the last
double sum as $\sum\limits_{k\leq
i}\sum\limits_{l=1}^N+\sum\limits_{k>i} \Bigl(\sum\limits_{l\leq
i}+\sum\limits_{i<l\leq k}+\sum\limits_{k<l}\Bigl)$ we obtain the
desired behavior of coefficients  $\gamma_{ik}(\tau)$ as
$\tau\to0$.

Consequently, we have just qualitative situation similar to those
on the
 $3^{\text{\sc rd}}$ {\sc Step} of the proof of Theorem~\ref{diffcurve}.
Thus the theorem follows.
\end{proof}

\subsubsection{$hc$-Differentiability of rectifiable curves}
In this section, we in particular prove that, in a Carnot
manifold, rectifiable curves are $hc$-differentiable almost
everywhere. We obtain this result as a corollary to the more
general assertion about the $hc$-differentiability of a mapping
$f:E\to {\mathbb M}$ from a measurable set $E\subset \mathbb R$
that satisfies the condition
\begin{equation}\label{stepcond}
\varlimsup\limits_{y\to x,\, y\in E}
\frac{d_{\infty}(f(y),f(x))}{|y-x|}<\infty
\end{equation}
for almost all $x\in E$.

\begin{thm}\label{Theorem6.3} Every mapping $f:E\to {\mathbb M}$ satisfying~\eqref{stepcond}
is $hc$-differentiable almost everywhere in~$E$.
\end{thm}

\begin{proof}  Since the result is local, we may assume that $E$ is
bounded. Since, in view of \eqref{stepcond}, the ``upper
derivative'' is finite almost everywhere, it follows that every
point~$x\in E\setminus\Sigma$, where $\Sigma\subset E$ is some set
of measure zero, belongs at least to one of the sets
\begin{equation}\label{5.2}
A_k=\biggl\{x\in E:\frac{{d}_{\infty}(f(x),f(y))}{|x-y|}\leq
k\quad \text{for all} \quad y\in E\cap
(x-k^{-1},x+k^{-1})\biggr\},\, k\in\mathbb N.
\end{equation}
Note that the sequence of sets $A_k$ is monotone: $A_k\subset
A_{k+1}$, $k\in\mathbb N$. Suppose that the measure of $A_{k}$ is
nonzero for some $k\in \mathbb N$. Up to a set of measure zero,
represent $A_{k}$ as the union of a disjoint family of sets
 $A_{k,1}, A_{k,2},\ldots$ of nonzero measure whose
diameters are at most~$1/k$:
$$
A_{k}=Z_{k}\cup A_{k,1}\cup A_{k,2}\cup\ldots,\quad |Z_{k}|=0.
$$
Then the restriction ~ $f_{k,j}=f|_{A_{k,j}}$ meets a Lipschitz
condition for all~$j$; therefore, it is extendable by continuity
to a Lipschitz mapping $\tilde
f_{k,j}:\overline{A}_{k,j}\to{\mathbb M}$.

Verify that if $(E\setminus\Sigma)\cap
(\overline{A}_{k,j}\setminus {A}_{k,j})\ne\emptyset$ then  $\tilde
f_{k,j}:(E\setminus\Sigma)\cap \overline{A}_{k,j}\to{\mathbb M}$
coincides with $f:(E\setminus\Sigma)\cap
\overline{A}_{k,j}\to{\mathbb M}$. In other words, if
$x\in(E\setminus\Sigma)\cap (\overline{A}_{k,j}\setminus
{A}_{k,j})$ then the extension of  $f: {A}_{k,j}\to{\mathbb M}$ by
continuity to the point $x$ equals $f(x)$. Indeed, the chosen
point $x$ belongs $E\setminus\Sigma$ and, therefore, $x\in A_l$
for some $l>k$. Then the inequality of \eqref{5.2}  holds for
$y\in E\cap (x-l^{-1},x+l^{-1})$ with $l$ instead of $k$.  Since
$A_l\cap (x-l^{-1},x+l^{-1})\supset
 A_{k,j}\cap (x-k^{-1},x+k^{-1})$, we have
 $$
 f(x)=\lim\limits_{y\to x,\, y\in A_l}f(y)=
 \lim\limits_{y\to x,\, y\in A_{k,j}}f(y)=\tilde
f_{k,j}(x).
$$

By Theorem~\ref{difflipcurve}, the mapping $\tilde
f_{k,j}:\overline{A}_{k,j}\to{\mathbb M}$ is $hc$-differentiable
almost everywhere in~$\overline{A}_{k,j}$. We are left with
checking the $hc$-differentiability of the mapping~$f:E\to{\mathbb
M}$ at the points of~$hc$-differentiability of the mapping $\tilde
f_{k,j}:\overline A_{k,j}\to{\mathbb M}$ having density one with
respect to~$\overline{A}_{k,j}$.

For brevity, denote the set $A_{k,j}$ by $A$ and denote the
mapping $f_{k,j}$ by $f$. Extend the Lipschitz
mapping~$f:A\to\mathbb M$ by continuity to a Lipschitz
mapping~$\tilde f:\overline A\to\mathbb M$.

Suppose now that a point $a\in A$ is a point of
$hc$-differentiability for~$\tilde f$ and the point density
of~$\overline A$. Recall that, by the definition of~$A$, the
inequality ${d}_{\infty}(f(y),f(z))\le k|y-z|$ holds for all $y\in
A$ and all $z\in (y-k^{-1},y+k^{-1})\cap E$. Note that this
inequality is extendable to $\overline A$ by continuity.
Consequently, the inequality
$$
{d}_{\infty}(\tilde f(y),f(z))\le k|y-z|
$$
holds for all $y\in\overline A$ and all $z\in
(y-k^{-1},y+k^{-1})\cap E$.

If $z\in E$ belongs to the neighborhood $(a-k^{-1},a+k^{-1})$
of~$a$ then, by the well-known property of a density point (see,
for example,~\cite{Stein}), there exists a point $y\in \overline
A$ such that $|y-z|=o(|z-a|$ as $z\to a$. Let $X$ be the
horizontal vector field of the definition of
$hc$-differentiability for the restriction $\tilde f:{\overline
A}\to\mathbb M$ at a point $a$. Then, in a sufficiently small
neighborhood of~$a$, from what was said above we have
\begin{multline*}
\!\!\!\!\! {d}_{\infty}(f(z),\exp((z-a) X)(f(a)) \leq
c^2({d}_{\infty}(f(z),\tilde f(y))+ {d}_{\infty}(\tilde
f(y),\exp(y-a) X)(f(a)))\\+ {d}_{\infty}(\exp(y-a)
X)(f(a)),\exp(z-a) X)(f(a))
\\
\leq c^2(k|y-z|+o(|y-a|)+ \|X\||y-z|)=o(|z-a|)
\end{multline*}
as $z\to a$, $z\in E$. Hence, the mapping $f:E\to{\mathbb M}$ is
$hc$-differentiable at $a$.

Suppose now that $k_1<k_2<k_3\ldots$ is a sequence of naturals
such that the measure of the complement $B_{k_j}=A_{k_j}\setminus
A_{k_{j-1}}$ is nonzero for every $j\geq2$. Obviously, the above
argument applies to each of the sets $B_{k_j}$, $j\geq2$, which
proves the theorem.
\end{proof}

Now we can prove the $hc$-differentiability of rectifiable curves.
Consider a curve (continuous mapping) $\gamma:[a,b]\to {\mathbb
M}$. By a partition $I_n=I_n([a,b])$ of the segment $[a,b]$ we
mean any finite sequence of points $\{s_1,\dots,s_n\}$ with
$a=s_1<\dots <s_n=b$. To every partition $I_n([a,b])$, we assign a
number $M(I_n)$ by setting
$$
M(I_n)=\sum\limits_{i=1}^n
d_{\infty}(\gamma(s_i),\gamma(s_{i+1})).
$$
Put $m_n=\max\{ s_{i+1}-s_i\mid i=1,\dots,n-1\}$.

\begin{defn}[\cite{bbi}]\label{Definition 6.3}
A curve $\gamma:[a,b]\to{\mathbb M}$ is called {\it rectifiable}
if
$$
L([a,b])=\lim\limits_{m_n\to 0}\sup\limits_{I_n} M_n <\infty.
$$
\end{defn}

Making use of a standard argument, we may prove:

\begin{property}\label{Property 6.3}
Suppose that a sequence of curves $\gamma_n:[a,b]\to{\mathbb M}$,
$n\in \mathbb N$, converges pointwise to a curve
$\gamma:[a,b]\to{\mathbb M}$$:$ $\gamma_n(s)\to\gamma(s)$ for
every $s\in[a,b]$. Then the lengths $L_n([a,b])$ of $\gamma_n$
possess the semicontinuity property$:$
$$
L([a,b])\leq\varliminf\limits_{n\to\infty}L_n([a,b]).
$$
\end{property}

\begin{prop}\label{Proposition6.3} Every rectifiable curve
$\gamma:[a,b]\to{\mathbb M}$ meets \eqref{stepcond}.
\end{prop}

\begin{proof} Consider the following set function $\Phi$ defined
on intervals included in~$[a,b]$: the value $\Phi(\alpha,\beta)$
at an interval $(\alpha,\beta)\subset[a,b]$ is equal to
$L([\alpha,\beta])$, the length of the curve
$\gamma:[\alpha,\beta]\to{\mathbb M}$. The set function $\Phi$ is
quasiadditive: the inequality
$$
\sum\limits_i\Phi(\alpha_i,\beta_i)\leq\Phi(\alpha,\beta)
$$
holds for every finite collection of pairwise disjoint intervals
$(\alpha_i,\beta_i)$ with
$(\alpha_i,\beta_i)\subset(\alpha,\beta)$, where
$(\alpha,\beta)\subset[a,b]$ is some interval. It is known (see,
for example, \cite{VU1}), that $\Phi$ has a finite derivative
$$
\Phi'(x)=\lim\limits_{\substack{(\alpha,\beta)\ni x,\\
\beta-\alpha\to0}} \frac{\Phi(\alpha,\beta)}{\beta-\alpha}=
\lim\limits_{\substack{(\alpha,\beta)\ni x,\\
\beta-\alpha\to0}}\frac{L([\alpha,\beta])}{\beta-\alpha}
$$
almost everywhere in~$[a,b]$. Hence,
$$
 \varlimsup\limits_{y\to x}
\frac{d_{\infty}(f(y),f(x))}{|y-x|}\leq
\varlimsup\limits_{\substack{(\alpha,\beta)\ni x,\\
\beta-\alpha\to0}}\frac{d_{\infty}(f(\alpha),f(\beta))}{L([\alpha,\beta])}\cdot
\lim\limits_{\substack{(\alpha,\beta)\ni x,\\
\beta-\alpha\to0}}\frac{L([\alpha,\beta])}{\beta-\alpha}
\leq\Phi'(x)<\infty
$$
for almost all $x\in[a,b]$.
\end{proof}

Theorem~\ref{Theorem6.3} and Proposition~\ref{Proposition6.3}
imply:

\begin{prop}\label{Proposition6.4} Every rectifiable curve
$\gamma:[a,b]\to{\mathbb M}$ is $hc$-differentiable almost
everywhere.
\end{prop}

\begin{rem}\label{Remark 6.2} If  the  Carnot
manifold is a Carnot group our definition of the
$hc$-differentiability of curves coincides with the  $\mathcal
P$-differentiability of curves given by P.~Pansu in~\cite{Pan}.
He proved also~\cite[ Proposition~4.1]{Pan} the  $\mathcal
P$-differentiability almost everywhere of rectifiable curves on
Carnot groups using a different method.
\end{rem}

\subsection{ $hc$-Differentiability of Smooth Mappings
of Carnot Manifolds}\label{diffsmmap}

\subsubsection{Continuity of horizontal derivatives and $hc$-differentiability of mappings} In this subsubsection, we generalize the classical property that the continuity
of the partial derivatives of a function defined on a Euclidean
space guarantees its differentiability.

In what follows, we repeatedly use the following correspondence:
to arbitrary element $a=\exp \Bigl( \sum\limits_{i=1}^Na_i\widehat
X{}^g_i\Bigr)(g)\in\mathcal G^g$ and point $w\in \mathcal G^g$,
assign the element

\begin{equation}\label{not1}
\Delta^w_\varepsilon a=
\exp\biggl(\sum\limits_{j=1}^Na_j\varepsilon^{\deg
X_j}X_j\biggr)(w)
\end{equation}
for those $\varepsilon$ for which the right-hand side
of~\eqref{not1} exists. Note that, by Property~\ref{coindil}, we
have $\Delta^g_\varepsilon a=\delta^g_\varepsilon a$ for all
$a\in\mathcal G^g$.

\begin{thm}\label{theorem3.3}
Suppose that $f:\mathbb M\to{\mathbb N}$ is a Lipschitz mapping of
Carnot manifolds such that, at each point $g\in \mathbb M$, there
exist horizontal derivatives $X_if(g)\in H_{f(g)}\mathbb N$
continuous on $\mathbb M$, $i=1,\dots,\dim H_1$. Then $f$ is
$hc$-differentiable at every point of~$\mathbb M$. The Lie algebra
homomorphism corresponding to the $hc$-differential is uniquely
defined by the mapping
$$
H_g\mathbb M\ni X_{i}(g)\mapsto X_{i}f(g)=\frac{d}{dt}f(\exp
tX_i(g))\vert_{t=0}=\sum\limits_{j=1}^{\dim\widetilde
H_1}b_{ij}Y_j(f(g))\in
 H_{f(g)}{\mathbb N}
$$
of the basis horizontal vectors $X_{i}(g)$, $i=1,\dots,\dim H_1$,
 to horizontal vectors
in $H_{f(g)}\mathbb N$$:$
$$
H\mathcal G^g\mathbb M\ni\widehat X^g_{i}\mapsto
\sum\limits_{j=1}^{\dim\widetilde H_1}b_{ij}\widehat Y^{f(g)}_j\in
H\mathcal G^{f(g)}{\mathbb N}.
$$
\end{thm}

\begin{proof} $1^{\text{\sc st}}$ {\sc Step.} Fix a point $g\in U$ and a compact neighborhood $F\subset \mathcal G^g$ of the local
Carnot group $\mathcal G^g$. For each horizontal vector field
$X_i$, a family of curves $\gamma:[-\varepsilon,\varepsilon]\times
F\to\mathbb N$ is defined: for $u\in F$, put
$\gamma_i(s,u)=f(\exp(s \alpha_i X_i)(u))$, where $\alpha_i\in A$,
$A\subset\mathbb R$ is a bounded neighborhood of $0\in \mathbb R$.
This family of curves meets the conditions of
Corollary~\ref{corollary3.2} and, hence, the convergence
\begin{equation}\label{3.5}
\Delta^{f(u)}_{s^{-1}}\gamma(s,u)\to\delta^{f(u)}_{\alpha_i}\exp (
[X_if](u))(f(u))\in\mathcal G^{f(u)}
\end{equation}
is uniform on $F\times A$ and the $hc$-derivative
$\delta^{f(u)}_{\alpha_i}\exp ( X_if(u))(u)$ is continuous with
respect to
 $(u,\alpha_i)\in F\times A$. Denote by $x_i$ the~``horizontal
basis element'' $\exp(X_i)(g)=\exp(\widehat X^g_i)(g)\in\mathcal
G^{g}$ and, for all $1\leq i\leq \dim H_1$, denote by $a_i$ the
horizontal derivative $\exp (X_if(g))(f(g))$.

It is known \cite{fs} that any element $v\in F$ can be represented
(nonuniquely) in the form
\begin{equation}\label{3.17}
\delta^g_{\alpha_1}x_{j_1}\cdot\dots\cdot
\delta^g_{\alpha_S}x_{j_S},\quad 1\leq j_i\leq \dim H_1,
\end{equation}
where $S$ is independent of the choice of the point and the
numbers $\alpha_i$ are bounded by a common constant. Together with
the mapping
$$
[0,\varepsilon)\ni t\mapsto \hat
v_i(t)=\delta^g_{t\alpha_1}x_{j_1}\cdot\dots\cdot
\delta^g_{t\alpha_i}x_{j_i},\quad 1\leq j_k\leq \dim H_1,\quad
1\leq k\leq i\leq S,
$$
consider the mapping (see \eqref{not1})
\begin{align*}
 [0,\varepsilon)\ni t\mapsto
v_i(t)&=\Delta^{v_{i-1}(t)}_{_{t\alpha_i}}x_{j_i}=\exp(t\alpha_iX_{j_i})(v_{i-1}(t)),\
2\leq
i\leq S,\ \text{where}\\
 v_1(t)&=\Delta^{g}_{_{t\alpha_1}}x_{j_1}=\exp(t\alpha_1X_{j_1})(g).
\end{align*}
By Theorem~\ref{mainresult}, $d_\infty(v_i(t),\hat v_i(t))=o(t)$
as $t\to0$ uniformly in $g\in F$ and $\alpha_i\in A$, $i\leq S$.
Since the mapping $f$ is Lipschitz on $F$, the limits
$\lim\limits_{t\to0}\Delta^{f(g)}_{t^{-1}}f(\hat v_{S}(t))$ and
$\lim\limits_{t\to0}\Delta^{f(g)}_{t^{-1}}f(v_{S}(t))$ exist
simultaneously. Consequently, it suffices to prove the existence
of the second limit.

$2^{\text{\sc nd}}$ {\sc Step.} For proving this, by \eqref{3.5},
we infer that
\begin{multline*}
 w_1(t)=f(v_{1}(t))=\exp\biggl(\sum\limits_{k=1}^{\widetilde N}z^1_k(t)Y_k\biggr)(f(g))\\
 \text {has $hc$-derivative}\quad
\delta^{f(g)}_{\alpha_1}a_{j_1}\in\mathcal G^{f(g)}\quad\text {at
$t=0$}.
\end{multline*}
Here $Y_k$, $k=1,\ldots,\widetilde N$, is a local basis on
$\mathbb N$ around the point $f(g)$. Assume that the mapping
\begin{multline*}
t\mapsto
w_i(t)=f(v_{i}(t))=\exp\biggl(\sum\limits_{k=1}^{\widetilde
N}z^i_k(t)Y_k\biggr)(f(v_{i-1}(t)))
 \\
 \text {has $hc$-derivative}\quad
\delta^{f(g)}_{\alpha_1}a_{j_1}\cdot\ldots\cdot\delta^{f(g)}_{\alpha_i}a_{j_i}\in\mathcal
G^{f(g)},\quad\text {at $t=0$}, \quad 2\leq i<S.
\end{multline*}
Our next goal is to show that  $hc$-derivative of the mapping
$t\mapsto w_{i+1}(t)=f(v_{i+1}(t))=
\exp\biggl(\sum\limits_{k=1}^{\widetilde
N}z^{i+1}_k(t)Y_k\biggr)(f(v_{i}(t)))$ equals
$\delta^{f(g)}_{\alpha_1}a_{j_1}\cdot\ldots\cdot\delta^{f(g)}_{\alpha_i}a_{j_i}\cdot\delta^{f(g)}_{\alpha_{i+1}}a_{j_{i+1}}$.
Together with the mapping  $w_{i+1}(t)$, consider the mapping
$$
t\mapsto\widehat w_{i+1}(t)=
\exp\biggl(\sum\limits_{k=1}^{\widetilde
N}z^{i+1}_k(t)\widehat{Y}^g_k\biggr)(f(v_{i}(t))).
$$
By Theorem~~\ref{mainresult} we have $d^{f(g)}_c(w_{i+1}(t),
\widehat w_{i+1}(t))=o(t)$ as  $t\to0$. Therefore, the relation
$d^{f(g)}_c\bigl(w_{i+1}(t),\delta^{f(g)}_t\bigl(\delta^{f(g)}_{\alpha_1}a_{j_1}\cdot\ldots\cdot\delta^{f(g)}_{\alpha_{i+1}}a_{j_{i+1}}\bigr)\bigr)=o(t)$
as $t\to0$ holds if and only if $d^{f(g)}_c\bigl(\widehat
w_{i+1}(t)),\delta^{f(g)}_t\bigl(\delta^{f(g)}_{\alpha_1}a_{j_1}\cdot\ldots\cdot\delta^{f(g)}_{\alpha_{i+1}}a_{j_{i+1}}\bigr)\bigr)=o(t)$
as $t\to0$. By Property~\ref{Prop2.1}, this is equivalent to the
relation
$$
d^{f(g)}_c\bigl(\delta^g_{t^{-1}}\widehat
w_{i+1}(t)),\delta^{f(g)}_{\alpha_1}a_{j_1}\cdot\ldots\cdot\delta^{f(g)}_{\alpha_{i+1}}a_{j_{i+1}}\bigr)=o(1)\quad\text{as}\quad
i\to\infty.
$$
Note that, by the continuity of the group operation in $\mathcal
G^g$, we always have the convergence
$$
\delta^g_{t^{-1}}\widehat
w_{i+1}(t))\to\delta^{f(g)}_{\alpha_1}a_{j_1}\cdot\dots\cdot
\delta^{f(g)}_{\alpha_{i+1}}a_{j_{i+1}} \quad\text{as $t\to0$.}
$$
Thus, by induction, the $hc$-derivative of the mapping
$[0,\varepsilon)\ni t\mapsto f(v_{S}(t))$ at $0$ is equal to
$\delta^{f(g)}_{\alpha_1}a_{j_1}\cdot\dots\cdot
\delta^{f(g)}_{\alpha_S}a_{j_S}$; moreover, the convergence is
uniform in $v\in F$ and $\alpha_{i}$, $1\leq i\leq S$.
Consequently, granted the equality $v_{S}(t)=\delta^g_t v$, we
infer
\begin{equation}\label{3.18}
 d_c^{f(g)}\bigl(f\bigl(\delta^g_{t}v\bigr),
 L\bigl(\delta^g_{t}v\bigr)\bigr)=
 o\bigl(d^g_c\bigl(g,\delta^g_{t}v\bigr)\bigr)
= o(t)
\end{equation}

uniformly in $v\in F$, where $L$ stands for the correspondence
$$
\mathcal G^g\ni v=\delta^g_{\alpha_1}x_{j_1}\cdot\dots\cdot
\delta^g_{\alpha_S}x_{j_S}\mapsto
\delta^{f(g)}_{\alpha_1}a_{j_1}\cdot\dots\cdot
\delta^{f(g)}_{\alpha_S}a_{j_S}\in\mathcal G^{f(g)}.
$$
For finishing the proof, it remains to check that the
correspondence $L:\mathcal G^g\to \mathcal G^{f(g)}$ is a
homomorphism of the local Carnot groups.

$3^{\text{\sc rd}}$ {\sc Step.}  Note that $L(v)$ is the
$hc$-derivative at $0$ of the mapping $t\mapsto
f\bigl(\delta^g_{t}v\bigr)$ for a fixed $v\in \mathcal G^g$ (see
\eqref{3.18}), which is obviously independent of representation
\eqref{3.17}. Consequently, $L:\mathcal G^g\to \mathcal G^{f(g)}$
is a mapping of the local groups. Clearly, this mapping is
continuous.

Demonstrate that it is a group homomorphism. Consider  a second
element $\overline v= \delta^g_{\beta_1}x_{j_1}\cdot\dots\cdot
\delta^g_{\beta_S}x_{j_S}$, $1\leq j_i\leq \dim H_1$, such that
\begin{equation}\label{3.19}
v\overline v=\delta^g_{\alpha_1}x_{j_1}\cdot\dots\cdot
\delta^g_{\alpha_S}x_{j_S}\cdot
\delta^g_{\beta_1}x_{j_1}\cdot\dots\cdot
\delta^g_{\beta_S}x_{j_S}\in\mathcal G^g\quad\text{and}\quad
L(v)\cdot L(\overline v)\in\mathcal G^{f(g)}.
\end{equation}

By \eqref{3.18}, the value $L(v\overline v)$ is independent of the
representation of an element $v\overline v$ as the product
\eqref{3.19}. Hence, applying the conclusions of the previous step
to $v\overline v$ and its representation \eqref{3.19}, we see that
$$
L(v\overline v)= \delta^{f(g)}_{\alpha_1}a_{j_1}\cdot\dots\cdot
\delta^{f(g)}_{\alpha_S}a_{j_S}\cdot
\delta^{f(g)}_{\beta_1}a_{j_1}\cdot\dots\cdot
\delta^{f(g)}_{\beta_S}a_{j_S}=L(v)\cdot L(\overline v).
$$
Thus, the mapping $L:\mathcal G^g\to{\mathcal G}^{f(g)}$ is a
continuous group homomorphism. By the well-known properties of the
Lie group theory~\cite{W}, the mapping $L$ is a homomorphism of
the local Lie groups.

Now, from \eqref{3.18} it can be deduced that $L$  commutes with a
dilation, $L\circ\delta^g_t=\delta^{f(g)}_t\circ L$, $t>0$.
Furthermore, since $X_if(g)\in H_{f(g)}{\mathbb M}$, the
homomorphism $L$ is the {\it $hc$-differential of the mapping
$f:{\mathbb M}\to{\mathbb N}$} at~$g$. The Lie algebra
homomorphism corresponding to $L$ is a mapping of horizontal
subspaces.
\end{proof}

\begin{cor}[\cite{vod3}]\label{cortheorem3.3}
Assume that we have a basis $\{X_i\}$, $i=1,\dots, N$, on a~Carnot
manifold $\mathbb M$ for which  Assumption {\ref{vectorfields}} or
conditions of Remark~{\ref{mainrem}} hold with some
$\alpha\in(0,1]$. Suppose that $f:\mathbb M\to{\mathbb N}$ is a
mapping of Carnot manifolds such that, at each point $g\in \mathbb
M$, there exist horizontal derivatives $X_if(g)\in H_{f(g)}\mathbb
N$ continuous on $\mathbb M$, $i=1,\dots,\dim H_1$. Then $f$ is
$hc$-differentiable at every point of~$\mathbb M$. The Lie algebra
homomorphism corresponding to the $hc$-differential is uniquely
defined by the mapping
$$
H_g\mathbb M\ni X_{i}(g)\mapsto X_{i}f(g)=\frac{d}{dt}f(\exp
tX_i(g))\vert_{t=0}=\sum\limits_{j=1}^{\dim\widetilde
H_1}b_{ij}Y_j(f(g))\in
 H_{f(g)}{\mathbb N}
$$
of the basis horizontal vectors $X_{i}(g)$, $i=1,\dots,\dim H_1$,
 to horizontal vectors
in $H_{f(g)}\mathbb N$$:$
$$
H\mathcal G^g\mathbb M\ni\widehat X^g_{i}\mapsto
\sum\limits_{j=1}^{\dim\widetilde H_1}b_{ij}\widehat Y^{f(g)}_j\in
H\mathcal G^{f(g)}\mathbb N.
$$
\end{cor}

\begin{proof}  The hypothesis implies that $f$
is a locally Lipschitz mapping: \linebreak$\tilde
d_\infty(f(x),f(y))\leq Cd_\infty(x,y)$, $x$, $y$ belong to some
compact neighborhood of $U$. To verify this, it suffices to join
points $x, y\in U$ by the horizontal curve $\gamma$ of
Subsection~{\ref{applications}} whose length is controlled by the
$hc$-distance $d_\infty(x,y)$ and observe that $f\circ\gamma$ is a
horizontal curve whose length is controlled by the length of the
initial curve. From this, Corollary~ \ref{equiv} and
Remark~\ref{equivlit}
 we infer
$\tilde d_\infty(f(x),f(y))\leq C_1 L(f\circ\gamma)\leq C_2
L(\gamma)\leq C_3d_\infty(x,y)$.
\end{proof}

\subsubsection{Functorial property of tangent cones} The definition of the tangent cone depends on the local
basis. The question arises on the connection between two tangent
cones found from two different bases. The last Theorem~3.3
implies:

\begin{cor}[\cite{vod2,vod3}]\label{cor3.5} Suppose that we have two local bases
$\{X_i\}$ and $\{Y_i\}$, $i=1,\dots, N$, on a~Carnot manifold  for
both of which  Assumption {\ref{vectorfields}} or conditions of
Remark~{\ref{mainrem}} hold with some $\alpha\in(0,1]$, and that
two collections
 $X_1,\dots$, $X_{\dim H_1}$ and
$Y_1,\dots,$ $Y_{\dim H_1}$ generate the same horizontal
subbundle~$H_1$. Then the tangent cone $\mathcal G^g$ defined by
the $\{X_i\}$'s is isomorphic to the local Carnot group
$\widetilde{\mathcal G}^g$, determined by the $\{Y_i\}$'s:
$(\widetilde\delta^g_{t^{-1}}\circ\delta^g_{t})(v)$ converges to
an isomorphism of local Carnot groups  $\mathcal G^g$ and
$\widetilde{\mathcal G}^g$ as $t\to0$ uniformly in
 $v\in \mathcal G^g$. $($Here
$\widetilde\delta^g_{t}$ is the one-parameter dilation group
defined by the vector fields $\{Y_i\}$.$)$

The isomorphism of the Lie algebras corresponding to the
$hc$-differential is defined uniquely by giving the mapping
$$
H_g{\mathbb M}\ni X_{i}(g)\mapsto X_{i}(g)=\sum\limits_{j=1}^{\dim
H_1}b_{ij}Y_j(g)\in
 H_g{\mathbb M}
$$
of the basis vectors $X_{i}(g)$, $i=1,\dots,\dim H_1$, of the
horizontal space~$H_g{\mathbb M}$ to horizontal vectors of the
space $H_g{\mathbb M}$$:$
$$
H\mathcal G^g{\mathbb M}\ni\widehat X_{i}^g\mapsto
\sum\limits_{j=1}^{\dim H_1}b_{ij}\widehat Y_j^g\in
H\widetilde{\mathcal G}^g\mathbb M.
$$
\end{cor}

\begin{proof} Denote by $\mathbb M^X$ the Carnot manifold
 $\mathbb M$ with local basis  $\{X_i\}$ and denote by $\mathbb M^Y$
 the Carnot manifold $\mathbb M$ with local basis
  $\{Y_i\}$, $i=1,\dots, N$.
 Let also the symbol $i:\mathbb M^X\to{\mathbb M}^Y$  stand for the identity mapping
from $\mathbb M$ into ${\mathbb M}$.  Clearly,  $i$ meets the
conditions of Corollary~\ref{cortheorem3.3}. Then $i$ is
$hc$-differentiable at $g$ and, by Corollary~\ref{cortheorem3.3},
the ``difference ratios''
$\widetilde\delta^g_{t^{-1}}(\delta^g_{t}(w))$ converge uniformly
to a~homomorphism  $Di(g):\mathcal G^g\to\widetilde{\mathcal G}^g$
as $t\to0$. Applying the same argument to the inverse mapping
$i^{-1}$ and Theorem~\ref{chainrule}, we infer that $Di(g)$ is an
isomorphism of the local Carnot groups (of the local tangent cones
at $g$ with respect to different local bases).
\end{proof}

\begin{rem} In \cite{AM,B,greshn,MM1} above statement
is proved by other methods under additional assumptions on the
smoothness of the basis vector fields.
\end{rem}

\subsubsection{Rademacher  Theorem}
The aim of this part is to formulate Rademacher type theorems on
the differentiability of Lipschitz mappings of Carnot manifolds.
This theorem was proved in \cite{vod3} by means of the theory
expounded above. The way of proving this result is based on the
methods of~\cite{V4}, where the $\mathcal P$-differentiability of
Lipschitz mappings of Carnot groups defined on measurable sets was
proved in details.

Let $\mathbb M$, $\mathbb N$ be two Carnot manifolds and let
$E\subset\mathbb M$ be an arbitrary set. A mapping $f:E\to{\mathbb
N}$ is called a Lipschitz mapping if
$$
\tilde d_{\infty}(f(x),f(y))\leq Cd_{\infty}(x,y),\quad x,y\in E,
$$
for some constant $C$ independent of $x$ and~$y$. The least
constant in this relation is denoted by $\operatorname{Lip} f$.

The following result extends the theorems on the $\mathcal
P$-differentiability on Carnot groups~\cite{Pan,V4,VU4} (see also
\cite{M1}) to Carnot manifolds.

\begin{thm}[\cite{vod3}]\label{Radth}
 Let $E$ be a set in $\mathbb M$ and let
$f:E\to\mathbb N$ be a Lipschitz mapping from $E$ into $\mathbb
N$. Then $f$ is $hc$-differentiable on $E$.

The homomorphism of the Lie algebras corresponding to the
$hc$-differential is defined uniquely by the mapping
$$
H_g\mathbb M\ni X_{i}(g)\mapsto X_{i}f(g)=\frac{d}{dt}f(\exp
tX_i(g))\vert_{t=0}=\sum\limits_{j=1}^{\dim\widetilde
H_1}a_{ij}Y_j(f(g))\in
 H_{f(g)}\mathbb N
$$
of the horizontal basis vectors $X_{i}(g)$, $i=1,\dots,\dim H_1$,
to horizontal vectors of the space $H_{f(g)}\mathbb N$$:$
$$
H\mathcal G^g\mathbb M\ni\widehat X_{i}^g\mapsto
\sum\limits_{j=1}^{\dim\widetilde H_1}a_{ij}\widehat Y_j^{f(g)}\in
H\mathcal G^{f(g)}\mathbb N.
$$
\end{thm}

\subsubsection{Stepanov Theorem}

As a corollary to Theorem~\ref{Radth}, we obtain a generalization
of Stepanov's theorem:

\begin{thm}[\cite{vod3}]\label{Stth}
 Let $E\subset\mathbb M$ be a set in $\mathbb M$
and let $f:E\to{\mathbb N}$ be a mapping such that
$$
\varlimsup_{x\to a,x\in E}
\frac{\tilde{d}_\infty(f(a),f(x))}{d_\infty(a,x)}<\infty
$$
for almost all~$a\in E$. Then $f$ is $hc$-differentiable almost
everywhere on $E$ and the $hc$-differential is unique.

The homomorphism of the Lie algebras corresponding to the
$hc$-differential is defined uniquely by the mapping
$$
H_g\mathbb M\ni X_{i}(g)\mapsto X_{i}f(g)=\frac{d}{dt}f(\exp
tX_i(g))\vert_{t=0}=\sum\limits_{j=1}^{\dim\widetilde
H_1}a_{ij}Y_j(f(g))\in
 H_{f(g)}\mathbb N
$$
of the basis horizontal vectors $X_{i}(g)$, $i=1,\dots,\dim H_1$,
 to horizontal vectors of the space
$H_{f(g)}\mathbb N$$:$
$$
H\mathcal G^g\mathbb M\ni\widehat X_{i}^g\mapsto
\sum\limits_{j=1}^{\dim\widetilde H_1}a_{ij}\widehat Y_j^{f(g)}\in
H\mathcal G^{f(g)}\mathbb N.
$$
\end{thm}

\section{Application: The Coarea Formula}\label{coareacarnot}

\subsection{Notations}
All the above results on geometry and differentiability are
applied in proving the sub-Riemannian analog of the well-known
coarea formula for some classes of contact mappings of Carnot~---
Carath\'{e}odory spaces.

\begin{notat}
Denote by $N_i$ the topological dimensions of $\mathbb M_i$ and
denote by $\nu_i$ the Hausdorff dimensions of $\mathbb M_i$,
$i=1,2$. Assume that
$$
T\mathbb M_1=\bigoplus\limits_{j=1}^{M_1}(H_j/H_{j-1}),\
H_0=\{0\},\text{ and }T\mathbb
M_2=\bigoplus\limits_{j=1}^{M_2}(\widetilde{H}_j/\widetilde{H}_{j-1}),\
\widetilde{H}_0=\{0\},
$$
where $H_1\subset T\mathbb M_1$ and $\widetilde{H}_1\subset
T\mathbb M_2$ are {\it horizontal} subbundles. The subspace
$H_j\subset T\mathbb M_1$ $(\widetilde{H}_j\subset T\mathbb M_2)$
is spanned by $H_1$ $(\widetilde{H}_1)$ and all commutators of
order not exceeding $j-1$, $j=2,\ldots, M_1$ $(M_2)$.

Denote the dimension of
${H}_j/{H}_{j-1}~(\widetilde{H}_j/\widetilde{H}_{j-1})$ by
$n_{j}~(\tilde{n}_j)$, $j=1,\ldots,M_1~(M_2)$.

Here the number $M_1~(M_2)$ are such that
${H}_{M_1}/{H}_{M_1-1}\neq
0~(\widetilde{H}_{M_2}/\widetilde{H}_{M_2-1}\neq 0)$, and
${H}_{M_1+1}/{H}_{M_1}=0~(\widetilde{H}_{M_2+1}/\widetilde{H}_{M_2}=0)$.
The number $M_1~(M_2)$ is called the {\it depth} of $\mathbb
M_1~(\mathbb M_2)$.
\end{notat}

\begin{assump}\label{assump_dim}
Suppose that
\begin{enumerate}
\item $N_1\geq N_2$;

\item $\dim H_i\geq\dim \widetilde{H}_i$, $i=1,\ldots, M_1$;

\item the basis vector fields $X_1,\ldots, X_{N_1}$ (in the
preimage) are $C^{1,\alpha}$-smooth, $\alpha>0$, and
$\widetilde{X}_1,\ldots,\widetilde{X}_{N_2}$ (in the image) are
$C^{1,\varsigma}$-smooth, $\varsigma>0$, or conditions of
Remark~{\ref{mainrem}} hold for $\alpha>0$ in the preimage and
$\varsigma>0$ in the image.
\end{enumerate}
\end{assump}

\begin{rem}
Note that, if there exists at least one point where the
$hc$-differential $\widehat{D}\varphi$ is non-degenerate, then the
condition $\dim H_1\geq\dim \widetilde{H}_1$ implies
$n_i\geq\tilde{n}_i$, $i=2,\ldots,M_1$ (compare with the above
assumption).
\end{rem}

\begin{notat}
Denote by $Z$ the set of points $x\in\mathbb M_1$ such that
\linebreak$\operatorname{rank}(D\varphi(x))<N_2$.
\end{notat}

\subsection{Lay-out of the Proof}
The key point in proving the non-holonomic coarea formula is to
investigate the interrelation of ''Riemannian`` and Hausdorff
measures on level sets (see below). The research on the comparison
of ''Riemannian`` and Hausdorff dimensions of submanifolds of
Carnot groups can be found in paper by Z.~M.~Balogh, J.~T.~Tyson
and B.~Warhurst~{\cite{btw}}. See other results on sub-Riemannian
geometric measure theory in works by L.~Ambrosio, F.~Serra Cassano
and D.~Vittone~{\cite{ascv}}, L.~Capogna, D.~Danielli, S.~D.~Pauls
and J.~T.~Tyson~{\cite{cdpt}}, D.~Danielli, N.~Garofalo and
D.-M.~Nhieu~{\cite{dgn}}, B.~Franchi, R.~Serapioni and F.~Serra
Cassano~{\cite{fssc1, fssc2}}, B.~Kirchheim and F.~Serra
Cassano~{\cite{ksc}}, V.~Magnani~{\cite{mblow}},
S.~D.~Pauls~{\cite{pauls}} and many other.

The purpose of Section~{\ref{coareacarnot}} is to explain the
ideas of proof of the coarea formula for sufficiently smooth
contact mappings $\varphi:\mathbb M_1\to\mathbb M_2$ of Carnot
manifolds. Note that, all the obtained results are new even for
the particular case of a mappings of Carnot groups.

\begin{rem}
For proving Theorems {\ref{tang_meas}}, {\ref{Leb_meas}},
{\ref{derivativemeas}}, and {\ref{degen}}, the smoothness $C^1$
(in Riemannian sense) for mappings $\varphi:\mathbb M_1\to\mathbb
M_2$ is sufficient. For proving Theorem {\ref{charset}}, the
(Riemannian) smoothness $C^{2,\varpi}$, $\varpi>0$, of $\varphi$
is sufficient.
\end{rem}

As it is mentioned above, for the first time, a  non-holonomic
analogue of the coarea formula is proved in paper of P.~Pansu
{\cite{P}}. The main idea of this work (which is used in many
other ones) is to prove the coarea formula via the Riemannian one:
\begin{multline}\label{riemsub}
{\eqref{coarea_e}}\Longrightarrow\int\limits_{U}{\mathcal
J}^{Sb}_{N_2}(\varphi,x)\,d\mathcal
H^{\nu_1}(x)\\=\int\limits_{\mathbb M_2}\,d\mathcal
H^{\nu_2}(z)\int\limits_{\varphi^{-1}(z)}\frac{\mathcal
J^{Sb}_{N_2}(\varphi,u)}{\mathcal J_{N_2}(\varphi,x)}\,d{\mathcal
H}^{N_1-N_2}(u)\overset{?}=\int\limits_{\mathbb M_2}\,d\mathcal
H^{\nu_2}(z)\int\limits_{\varphi^{-1}(z)}\,d{\mathcal
H}^{\nu_1-\nu_2}(u)
\end{multline}
Here $N_1, N_2$ are topological dimensions, and $\nu_1,\nu_2$ are
Hausdorff dimensions of preimage and image, respectively; it is
well-known that, in sub-Riemannian case, topological and Hausdorff
dimensions differ. It easily follows from {\eqref{riemsub}}, that
the key point in this problem is to   investigate the
interrelation of ''Riemannian`` and Hausdorff measures on Carnot
manifolds theirselves and on level sets of $\varphi$, and of
Riemannian and sub-Riemannian coarea factors. It is well known
that the question on interrelation of measures on Carnot manifolds
is quite easy, while both the investigation of geometry of level
sets and the calculation of sub-Riemannian coarea factor are
non-trivial. The main problems are connected with peculiarities of
a sub-Riemannian metric. In particular, the non-equivalence of
Riemannian and sub-Riemannian metrics can be seen in the fact that
''Riemannian`` radius of a sub-Riemannian ball of a radius $r$
varies from $r$ to $r^M$, $M>1$, where the constant $M$ depends on
the Carnot manifold structure. Thus, a question arises immediately
on how ''sharp`` the approximation of a level  by its tangent
plain is (since the ''usual`` order of tangency $o(r)$ is
obviously insufficient here: a level may ''jump`` from a ball
earlier then it is  expected). Also a question arises on existence
of a such sub-Riemannian metric suitable for the description of
the geometry of an intersection of a ball and a level set. But
even if we answer these questions, one more question appears: what
is the relation of the Hausdorff dimension of the image and
measure of the intersection of a ball and a level set.

We have solved all the above problems. First of all, the points in
which the differential is non-degenerate, are divided into two
sets: regular and characteristic.
\begin{defn}
The set
$$
\chi=\{x\in\mathbb M_1\setminus
Z:\operatorname{rank}\widehat{D}\varphi(x)<N_2\}
$$
is called the {\it characteristic} set. The points of $\chi$ are
called {\it characteristic}.
\end{defn}

\begin{defn}
The set
$$
\mathbb D=\{x\in\mathbb
M_1:\operatorname{rank}\widehat{D}\varphi(x)=N_2\}
$$
is called the {\it regular} set. If $x\in\mathbb D$, then we say
that, $x$ is a {\it regular} point.
\end{defn}

We define a number $\nu_0(x)$ depending on $x\in\mathbb M$ that
shows whether a point is regular or characteristic.

\begin{defn}
Consider the number $\nu_0$ such that
\begin{multline*}
\nu_0(x)=\min\Bigl\{\nu: \exists\{X_{i_1},\ldots, X_{i_{N_2}}\}\\
\bigl(\operatorname{rank}([X_{i_j}\varphi](x))_{j=1}^{N_2}=N_2\bigl)\Rightarrow\Bigl(\sum\limits_{j=1}^{N_2}\operatorname{deg}
X_{i_j}=\nu\Bigr)\Bigr\}.
\end{multline*}
\end{defn}

It is clear that ${\nu_0}|_{\chi}>\nu_2$ and ${\nu_0}|_{\mathbb
D}=\nu_2$.

We also define such sub-Riemannian quasimetric~$d_2$, that makes
the calculation of measure of the intersection of a sub-Riemannian
and a tangent plain to a level set possible:

\begin{defn}
Let $\mathbb M$ be a Carnot manifold of topological dimension $N$
and of depth $M$, and let $x=\exp\Bigl(\sum\limits_{i=1}^{N_1}x_i
X_i\Bigr)(g)$. Define the distance $d_2(x,g)$ as follows:
\begin{multline*}
d_2(x,g)=\max\Bigl\{\Bigl(\sum\limits_{j=1}^{n_1}|x_j|^2\Bigr)^{\frac{1}{2}},\\
\Bigl(\sum\limits_{j=n_1+1}^{n_1+n_2}|x_j|^2\Bigr)^{\frac{1}{2\cdot\operatorname{deg}
X_{n_1^1+1}}},\ldots,\Bigl(\sum\limits_{j=N-n_{M}+1}^{N}|x_j|^2\Bigr)^{\frac{1}{2\cdot\operatorname{deg}
X_{N}}}\Bigr\}.
\end{multline*}
The similar metric $d_2^u$ is introduced in the local Carnot group
$\mathbb G_u\mathbb M$.
\end{defn}

The construction of $d_2$ is based on the fact that a ball in this
quasimetric $\operatorname{Box}_{2}$ asymptotically equals a
Cartesian product of Euclidean balls:
$$
\operatorname{Box}_{2}(x,r)\approx B^{n_1}(x,r)\times
B^{n_2}(x,r^2)\times\ldots\times B^{n_M}(x,r^M), \ M>1,
$$
where $N,n_i$, $i=1,\ldots, M$, are (topological) dimensions of
balls. The latter fact makes the calculation of above-mentioned
measure possible (while in the case when we replace balls by
cubes, it is quite complicated since cubes have different shapes
of sections).

Using properties of this quasimetric, we calculate the $\mathcal
H^{N_1-N_2}$-measure of the intersection of a tangent plain to a
level set and a sub-Riemannian ball in the quasimetric $d_2$.

\begin{thm}\label{tang_meas}
 Fix $x\in\varphi^{-1}(t)$. Then, the
${\mathcal H}^{N_1-N_2}$-measure of the intersection
$T_0[(\varphi\circ\theta_x)^{-1}(t)]\cap\operatorname{Box}_2(0,r)$
is equivalent to
$$
C(1+o(1))r^{\nu_1-\nu_0(x)}
$$
where $C$ does not depend on $r$, and $o(1)\to0$ as $r\to0$.
\end{thm}

While investigating the approximation of a surface by its tangent
plain, we introduce a ''mixed`` metric possessing some Riemannian
and sub-Riemannian properties.

\begin{defn}
For $v,w\in\operatorname{Box}_2(0,r)$ put
${d_2^0}_E(v,w)=d_2^0(0,w-v)$, where $w-v$ denotes the Euclidean
difference.
\end{defn}
 This definition implies that
$\operatorname{Box}_{2}(0,r)$ coincides with a ball
$\operatorname{Box}_{2E}(0,r)$ centered at $0$ of radius $r$ in
the metric ${d_2^0}_E$.

We prove that in regular points the tangent plain approximates the
level set quite sharp with respect to this metric, and from here
we deduce the possibility of calculation of the Riemannian measure
of a level set and a sub-Riemannian ball intersection. Notable is
the fact that this measure can be expressed via Hausdorff
dimensions of the preimage and the image: it is equivalent to
$r^{\nu_1-\nu_2}$ (see below):

\begin{thm}\label{Leb_meas}
Suppose that $x\in\varphi^{-1}(t)$ is a regular point. Then$:$

{\bf(I)} In the neighborhood of $0=\theta_x^{-1}(x)$, there exists
a mapping from
$T_0[(\varphi\circ\theta_x)^{-1}(t)]\cap\operatorname{Box}_2(0,r(1+o(1)))$
to $\psi^{-1}(t)\cap\operatorname{Box}_2(0,r)$, such that both
$d_2$- and $\rho$-distortions with respect to $0$ equal $1+o(1)$,
where $o(1)$ is uniform on $\operatorname{Box}_2(0,r)$$;$

{\bf (II)} The ${\mathcal H}^{N_1-N_2}$-measure of the
intersection $\varphi^{-1}(t)\cap\operatorname{Box}_2(x,r)$ equals
$$
\bigl|g|_{\ker
D\varphi(x)}\bigr|\cdot\prod\limits_{k=1}^{M_1}\omega_{n^1_k-n^2_k}\cdot
\frac{|{D}\varphi(x)|}{|\widehat{D}\varphi(x)|}r^{\nu_1-\nu_2}(1+o(1)),
$$
where $g$ is a Riemann tensor, $\widehat{D}\varphi$ is the
$hc$-differential of $\varphi$, and $o(1)\to0$ as $r\to0$.
\end{thm}

From these results and obtained properties, using a result of
{\cite{VU1}}, we deduce the interrelation of two measures in
regular points of a level sets.

\begin{thm}[Measure Derivative on Level Sets]
\label{derivativemeas} Hausdorff measure $\mathcal
H^{\nu_1-\nu_2}$ of the intersection
$\operatorname{Box}_2(x,r)\cap\varphi^{-1}(\varphi(x))$, where $x$
is a regular point, and
\linebreak$\operatorname{dist}(\operatorname{Box}_2(x,r)\cap\varphi^{-1}(\varphi(x)),\chi)>0$,
asymptotically equals $\omega_{\nu_1-\nu_2}r^{\nu_1-\nu_2}$. The
derivative $D_{{\mathcal H}^{N_1-N_2}}{\mathcal
H}^{\nu_1-\nu_2}(x)$ equals
$$
\frac{1}{\bigl|g|_{\ker
D\varphi(x)}\bigr|}\cdot\frac{\omega_{\nu_1-\nu_2}}{\prod\limits_{k=1}^{M_1}{\omega_{n_k-\tilde{n}_k}}}\cdot
\frac{|\widehat{D}\varphi(x)|}{|{D}\varphi(x)|}.
$$
\end{thm}

Finally, we introduce the notion of the sub-Riemannian coarea
factor via the values of the $hc$-differential of $\varphi$.

\begin{defn}
The {\it sub-Riemannian coarea factor} equals
$$
{\mathcal J}^{SR}_{N_2}(\varphi,x)=|\widehat{D}\varphi(x)|\cdot
\frac{\omega_{N_1}}{\omega_{\nu_1}}\frac{\omega_{\nu_2}}{\omega_{N_2}}
\frac{\omega_{\nu_1-\nu_2}}{\prod\limits_{k=1}^{M_1}\omega_{n_k-\tilde{n}_k}}.
$$
\end{defn}

We consider and solve problems connected with the characteristic
set. The case of characteristic points is a little more
complicated since in characteristic points a surface may jump from
a sub-Riemannian ball, consequently, we cannot estimate the
measure of the intersection of the surface and the ball via the
one of the tangent plain and the ball. Note also that in all the
other works on sub-Riemannian coarea formula, the preimage has a
group structure, which is essentially used in proving the fact
that the Hausdorff measure of characteristic points on each level
set equals zero (see also the paper {\cite{Bal}} by Z.~M.~Balogh,
dedicated to properties of the characteristic set). In the case of
a mapping of two Carnot manifolds, there is no group structure
neither in image, nor in preimage. Moreover, the approximation of
Carnot manifold by its local Carnot group is insufficient for
generalization of methods developed before. That is why we
construct new ''intrinsic`` method of investigation of properties
of the characteristic set. First of all, in all the characteristic
points the $hc$-differential is degenerate. We solve this problem
with the following assumption.

\begin{property}\label{assumpchi}
Suppose that $x\in\chi$, and
$\operatorname{rank}\widehat{D}\varphi(x)=N_2-m$. Let also
$\widehat{D}\varphi(x)$ equals zero on $n_1-\tilde{n}_1+m_1$
horizontal (linearly independent) vectors, $n_2-\tilde{n}_2+m_2$
(linearly independent) vectors from $H_2/H_1$,
$n_k-\tilde{n}_k+m_k$ (linearly independent) vectors from
$H_k/H_{k-1}$, $k=3,\ldots, M_2$. Then, on the one hand, since
$\operatorname{rank}\widehat{D}\varphi(x)=N_2-m$, we have
$\sum\limits_{i=1}^{M_1}m_k=m$. On the other hand,
$\operatorname{rank}D\varphi(x)=N_2$. Consequently, there exist
$m$ (linearly independent) vectors $Y_1,\ldots, Y_m$ of degrees
$l_1,\ldots, l_{M_2}$ (which are minimal) from the kernel of the
$hc$-differential $\widehat{D}\varphi$, such that
$D\varphi(x)(\operatorname{span}\{H_{M_2},Y_1,\ldots,
Y_m\})=T_{\varphi(x)}\mathbb M_2$.

In this subsection, we will assume that, among the vectors
$Y_1,\ldots, Y_m$, $m_1$ of them of the degree $l_1$ have the
horizontal image, $m_2$ of them of the degree $l_2\geq l_1$ have
image belonging to $\widetilde{H}_2$, and $m_k$ of them of the
degree $l_k$, $l_k\geq l_{k-1}$, have image belonging to
$\widetilde{H}_k$, $k=3,\ldots, M_2$.
\end{property}

By another words, the ''extra`` vector fields on which the
$hc$-differential of $\varphi$ is degenerate in characteristic
points, possess the following property: if in $H_k/H_{k-1}(x)$ the
quantity of such ''extra`` vectors equals $m_k>0$, then there
exist $m_k$ vectors from $H_{l_k}/H_{l_k-1}(x)$ such that their
images have the degree $k$, they are linearly independent with
each other and with the images of $H_{l_k-1}(x)$, $l_k\geq
l_{k-1}$. We develop new ''intrinsic`` method of investigation of
the properties of the characteristic set.

\begin{ex}
\label{exchi} The condition described in Assumption
{\ref{assumpchi}}, is always valid for the following $\mathbb M_1$
and $\mathbb M_2$:

\begin{enumerate}
\item $\mathbb M_1$ is an arbitrary Carnot--Carath\'{e}odory
space, and $\mathbb M_2=\mathbb R$;

\item $M_1$ is an arbitrary Carnot--Carath\'{e}odory space of the
topological dimension $2m+1$, $\mathcal G^u\mathbb M_1=\mathbb
H^{m}$ for all $u\in\mathbb M_1$, $\mathbb M_2=\mathbb R^k$,
$k\leq 2m$;

\item $M_1=M_2$, $\dim H_1\geq\dim\widetilde{H}_1$,
$\dim(H_i/H_{i-1})=\dim(\widetilde{H}_i/\widetilde{H}_{i-1})$,
$i=2,\ldots, M_1$;

\item $M_1=M_2+1$, $\dim H_i=\dim\widetilde{H}_i$, $i=1,\ldots,
M_2$.
\end{enumerate}
\end{ex}

In particular, in Theorem {\ref{tang_meas}} it is shown, that in
the characteristic points $\mathcal H^{N_1-N_2}$-measure of the
intersection of a sub-Riemannian ball and the tangent plain to the
level set is equivalent to $r$ to the power
$\nu_1-\nu_0(x)<\nu_1-\nu_2$. Next, we show, that $\mathcal
H^{N_1-N_2}$-measure of the intersection of the level set and the
sub-Riemannian ball centered at a characteristic point is
infinitesimally big in comparison with $r^{\nu_1-\nu_2}$, i.~e.,
is equivalent to $\frac{r^{\nu_1-\nu_2}}{o(1)}$ (but it is not
necessarily equivalent to $r^{\nu_1-\nu_0(x)}$). From here we
deduce that, the intersection of the characteristic set with each
level set has zero $\mathcal H^{\nu_1-\nu_2}$-measure.

\begin{thm}
[Size of the Characteristic Set]\label{charset} The Hausdorff
measure \linebreak${\mathcal
H}^{\nu_1-\nu_2}(\chi\cap\varphi^{-1}(t))=0$ for all $z\in\mathbb
M_2$.
\end{thm}

We also show that the degenerate set of the differential does not
influence both parts of the coarea formula.

\begin{thm}\label{degen}
For ${\mathcal H}^{\nu_2}$-almost all $t\in\mathbb M_2$, we have
$$
{\mathcal H}^{\nu_1-\nu_2}(\varphi^{-1}(t)\cap Z)=0.
$$
\end{thm}

Finally, we  deduce the sub-Riemannian coarea formula.

\begin{thm}
For any smooth contact mapping $\varphi:\mathbb M_1\to\mathbb M_2$
possessing Property~{\ref{assumpchi}}, the coarea formula holds:
\begin{equation*}
\int\limits_{\mathbb M_1}{\mathcal
J}_{N_2}^{Sb}(\varphi,x)\,d{\mathcal H}^{\nu_1}(x)
=\int\limits_{\mathbb M_2}\,d{\mathcal
H}^{\nu_2}(t)\int\limits_{\varphi^{-1}(t)}\,d{\mathcal
H}^{\nu_1-\nu_2}(u).
\end{equation*}
\end{thm}

As an application, using the result of the paper by R.~Monti and
F.~Serra Cassano~{\cite[Theorem~4.2]{MS}} for
$\operatorname{Lip}$-functions defined on a
Carnot--Carath\'{e}odory space $\mathbb M$ of the Hausdorff
dimension $\nu$, we deduce that the De Giorgi perimeter coincides
with $\mathcal H^{\nu-1}$-measure on almost every level of a
smooth function $\varphi:\mathbb M$:

\begin{thm}
For $C^{2,\alpha}$-functions $\varphi:\mathbb M\to\mathbb R$,
$\alpha>0$, where $\dim_{\mathcal H}\mathbb M=\nu$, the De Giorgi
perimeter coincides with $\mathcal H^{\nu-1}$-measure on almost
every level.
\end{thm}

\section{Appendix}\label{proofdiff}

\subsection{Proof of Lemma {\ref{ODE}}}

\begin{proof}
It is well known, that the solution $y(t,u)$ of the ODE
{\eqref{diffequation}} equals
$y(t,u)=\lim\limits_{n\to\infty}y_n(t,u)$, where
$$
y_0(t,u)=\int\limits_0^tf(y(0),u)\,d\tau,\ \text{  and  }\
y_n(t,u)=\int\limits_0^tf(y_{n-1}(\tau,u),u)\,d\tau.
$$
This convergence is uniform in $u$, if $u$ belongs to some compact
set.

From the definition of this sequence it follows, that $y_n(t)\to
y(t)$ as $n\to\infty$ in $C^1$-norm.

{\bf 1.} We show, that every $y_n(t,u)\in C^{\alpha}(u)$ for each
$t\in[0,1]$. We have
\begin{multline*}
\max\limits_{t}|y_n(t,u_1)-y_n(t,u_2)|\\
\leq\int\limits_0^1
|f(y_{n-1}(\tau,u_1),u_1)-f(y_{n-1}(\tau,u_2),u_2)|\,d\tau\\
\leq \int\limits_0^1
|f(y_{n-1}(\tau,u_1),u_1)-f(y_{n-1}(\tau,u_1),u_2)|\,d\tau\\
+\int\limits_0^1
|f(y_{n-1}(\tau,u_1),u_2)-f(y_{n-1}(\tau,u_2),u_2)|\,d\tau\\
\leq
H(f)|u_1-u_2|^{\alpha}+L\max\limits_t|y_{n-1}(t,u_1)-y_{n-1}(t,u_2)|\\
\leq H(f)\sum\limits_{m=0}^{n-1}L^{m}|u_1-u_2|^{\alpha}
+L^{n}\max\limits_t|y_{0}(t,u_1)-y_{0}(t,u_2)|\\
\leq H(f)\sum\limits_{m=0}^{\infty}L^{m}|u_1-u_2|^{\alpha},
\end{multline*}
where $H(f)$ is a constant, such that $|f(u_1)-f(u_2)|\leq
H(f)|u_1-u_2|^{\alpha}$. Note that the constant
$H=H(f)\sum\limits_{m=0}^{\infty}L^{m}<\infty$ since $L<1$, and it
does not depend on $n\in\mathbb N$.

Suppose that $u$ belongs to some compact set $U$. Then
\begin{multline*}
|y(t,u_1)-y(t,u_2)|\\
\leq
|y(t,u_1)-y_n(t,u_1)|+|y_n(t,u_1)-y_n(t,u_2)|+|y(t,u_2)-y_n(t,u_2)|\\
\leq H|u_1-u_2|^{\alpha}+2\varepsilon
\end{multline*}
for every $\varepsilon=\varepsilon(n)>0$. Since the convergence is
uniform in $u\in U$, and $\varepsilon(n)\to0$ as $n\to\infty$,
then $|y(t,u_1)-y(t,u_2)|\leq H|u_1-u_2|^{\alpha}$, and $y\in
C^{\alpha}(u)$ locally.

To show, that $\frac{\partial y}{\partial v_i}(t,v,u)\in
C^{\alpha}(u)$ locally, $i=1,\ldots, N$, we obtain our estimates
in the simplest case of $N=1$.

{\bf 2.} Note that the mappings $\{y_n\}_{n\in\mathbb N}$ converge
to $y$ in $C^{1}$-norm, and this convergence is uniform, if $u$
belongs to some compact set $U$.

Let $u\in U$, $v\in W(0)\subset\mathbb R^N$. Then similarly to the
case {\bf 1}, we see, that if the H\"{o}lder constant of
$y_{n}^{\prime}$ does not depend on $n\in\mathbb N$, then
$y^{\prime}\in C^{\alpha}(u)$.
\begin{multline}\label{1st}
\max\limits_{t,v}\Bigl|\frac{dy_n}{dv}(t,v,u_1)-\frac{dy_n}{dv}(t,v,u_2)\Bigr|\\
\leq\max\limits_{t,v}\Bigl|\frac{d}{dv}\int\limits_0^tf(y_{n-1}(\tau,v,u_1),v,u_1)-f(y_{n-1}(\tau,v,u_2),v,u_2)\,d\tau\Bigr|\\
\leq\max\limits_{t,v}\Bigl|\frac{d}{dv}\int\limits_0^tf(y_{n-1}(\tau,v,u_1),v,u_1)-f(y_{n-1}(\tau,v,u_2),v,u_1)\,d\tau\Bigr|\\
+\max\limits_{t,v}\Bigl|\frac{d}{dv}\int\limits_0^tf(y_{n-1}(\tau,v,u_2),v,u_1)-f(y_{n-1}(\tau,v,u_2),v,u_2)\,d\tau\Bigr|.
\end{multline}
For the first summand we have
\begin{multline}\label{1.5}
\max\limits_{t,v}\Bigl|\frac{d}{dv}\int\limits_0^1f(y_{n-1}(\tau,v,u_1),v,u_1)-f(y_{n-1}(\tau,v,u_2),v,u_1)\,d\tau\Bigr|\\
\leq\max\limits_v\int\limits_{0}^1\Bigl|\frac{d}{dv}(f(y_{n-1}(\tau,v,u_1),v,u_1)-f(y_{n-1}(\tau,v,u_2),v,u_1))\Bigr|\,d\tau\\
\leq\max\limits_v\int\limits_{0}^1\Bigl|\frac{df}{dy}\frac{dy_{n-1}}{dv}(\tau,v,u_1)
-\frac{df}{dy}\frac{dy_{n-1}}{dv}(\tau,v,u_2)\Bigr|\,d\tau\\
+\max\limits_v\int\limits_{0}^1\Bigl|\frac{\partial f}{\partial
v}(y_{n-1}(\tau,v,u_1))-\frac{\partial f}{\partial
v}(y_{n-1}(\tau,v,u_2))\Bigr|\,d\tau.
\end{multline}
Then, we get
$$
\max\limits_v\int\limits_{0}^1\Bigl|\frac{\partial f}{\partial
v}(y_{n-1}(\tau,v,u_1))-\frac{\partial f}{\partial
v}(y_{n-1}(\tau,v,u_2))\Bigr|\,d\tau \leq
C(f)H(y)|u_1-u_2|^{\alpha},
$$
since each $y_{m}$ is H\"{o}lder. The first summand in
{\eqref{1.5}} is evaluated in the following way:
\begin{multline}\label{2nd}
\max\limits_v\int\limits_{0}^1\Bigl|\frac{df}{dy}\frac{dy_{n-1}}{dv}(\tau,v,u_1)
-\frac{df}{dy}\frac{dy_{n-1}}{dv}(\tau,v,u_2)\Bigr|\,d\tau\\
\leq\max\limits_v\int\limits_{0}^1\Bigl|\frac{df}{dy}(u_1)\frac{dy_{n-1}}{dv}(\tau,v,u_1)
-\frac{df}{dy}(u_1)\frac{dy_{n-1}}{dv}(\tau,v,u_2)\Bigr|\,d\tau\\
+\max\limits_v\int\limits_{0}^1\Bigl|\frac{df}{dy}(u_1)\frac{dy_{n-1}}{dv}(\tau,v,u_2)
-\frac{df}{dy}(u_2)\frac{dy_{n-1}}{dv}(\tau,v,u_2)\Bigr|\,d\tau\\
\leq
L\max\limits_{t,v}\Bigl|\frac{dy_{n-1}}{dv}(t,v,u_1)-\frac{dy_{n-1}}{dv}(t,v,u_2)\Bigr|\\
+\max\limits_{u,v}\int\limits_{0}^1\Bigl|\frac{dy_{n-1}}{dv}(\tau,v,u)\Bigr|\,d\tau\cdot
H(Df)|u_1-u_2|^{\alpha}.
\end{multline}
Next, we estimate
\begin{multline*}
\max\limits_{u,v}\int\limits_{0}^1\Bigl|\frac{dy_{m}}{dv}(\tau,v,u)\Bigl|\,d\tau
\leq\max\limits_{t,u,v}\Bigl|\frac{dy_{m}}{dv}(t,v,u)\Bigr|\\
=\max\limits_{t,u,v}
\Bigl[L\Bigl|\frac{dy_{m-1}}{dv}\Bigr|+\Bigl|\frac{\partial
f}{\partial
v}\Bigr|\Bigr]\leq\max\limits_{t,u,v}\Bigl|\frac{\partial
f}{\partial
v}\Bigr|\Bigl[\sum\limits_{k=0}^{\infty}L^k\Bigr]<\infty.
\end{multline*}
Thus, in the first summand of {\eqref{1st}} we have
$$
L\max\limits_{t,v}\Bigl|\frac{dy_{n-1}}{dv}(t,v,u_1)-\frac{dy_{n-1}}{dv}(t,v,u_2)\Bigr|+C|u_1-u_2|^{\alpha},
$$
where $0<C<\infty$ does not depend on $n\in\mathbb N$. The second
summand in {\eqref{1st}} is
\begin{multline*}
\max\limits_{t,v}\Bigl|\frac{d}{dv}
\int\limits_0^tf(y_{n-1}(\tau,v,u_2),v,u_1)-f(y_{n-1}(\tau,v,u_2),v,u_2)\,d\tau\Bigr|\\
\max\limits_{v}\int\limits_0^1\Bigl|\frac{\partial f}{\partial
v}(y_{n-1},v,u_1)-\frac{\partial f}{\partial
v}(y_{n-1},v,u_2)\Bigr|\,d\tau\leq C(f)|u_1-u_2|.
\end{multline*}

Thus,
\begin{multline*}
\max\limits_{t,v}\Bigl|\frac{dy_n}{dv}(t,v,u_1)-\frac{dy_n}{dv}(t,v,u_2)\Bigr|\\
\leq
L\max\limits_{t,v}\Bigl|\frac{dy_{n-1}}{dv}(t,v,u_1)-\frac{dy_{n-1}}{dv}(t,v,u_2)\Bigr|+K|u_1-u_2|^{\alpha}\\
\leq k\sum\limits_{k=0}^{\infty}L^k|u_1-u_2|^{\alpha},
\end{multline*}
and $\frac{dy_n}{dv}\in C^{\alpha}(u)$ locally. Hence,
$\frac{dy}{dv}\in C^{\alpha}(u)$ locally.
\end{proof}

\subsection*{Acknowledgment}

We are grateful to Valeri\v{\i} Berestovski\v{\i}, for his collaboration
 in proving of Gromov's theorem on nilpotentization of vector fields (see Remark~{\ref{Ber}}).
We also want to thank  Pierre Pansu for fruitful discussion of our
results concerning the geometry of Carnot--Carath\'{e}odory spaces
and the sub-Rieman\-nian coarea formula, and H.~Martin Reimann for
interesting discussions of differentiability theorems.

\medskip

{The research was partially supported by the Commission of the
European Communities (Specific Targeted Project ``Geometrical
Analysis in Lie groups and Applications'', Contract
number~028766), the Russian Foundation for Basic Research
(Grant~06-01-00735), the State Maintenance Program for Young
Russian Scientists and the Leading Scientific Schools of Russian
Federation (Grant~NSh-5682.2008.1). }


\begin{thebibliography}{555}
\bibitem{agr} A.~A.~Agrachev, \textit{Compactness for sub-Riemannian length minimizers and subanalyticity.} Rend. Semin. Mat. Torino, \textbf{56} (1998).
\bibitem{agrgam} A. A. Agrachev and R. Gamkrelidze. \textit{Exponential representation of flows
and chronological calculus.} Math. USSR-Sb. \textbf{107} (4) (1978), 487--532 (in Russian).
\bibitem{agrg} A. A. Agrachev and J.-P. Gauthier. \textit{On subanalyticity of Carnot-Carath\'eodory
distances.}  Ann. Inst. H. Poincar\'{e} Anal. Non Lin\'{e}aire, \textbf{18} (3), 2001.
\bibitem{AM} A.~A. Agrachev, A. Marigo, \textit{Nonholonomic tangent
spaces{\rm:} intrinsic construction and rigid dimensions.}
Electron. Res. Announc. Amer. \textbf{9} (2003), 111--120.
\bibitem{agrs1}
A. A. Agrachev and A. V. Sarychev, \textit{Filtrations of a Lie algebra of vector fields and
nilpotent approximations of control systems.} Dokl. Akad. Nauk SSSR \textbf{285} (1987), 777--781.
\bibitem{agrs2}
A. A. Agrachev and A. V. Sarychev, \textit{Strong minimality of abnormal geodesics for
2-distributions.} J. Dyn. Control Syst. \textbf{1} (2) (1995).
\bibitem{agrs3}
A. A. Agrachev and A. V. Sarychev, \textit{Abnormal sub-Riemannian
geodesics: Morse index and rigidity.} Ann. Inst. Henri
Poincar\'{e}, Analyse Non Lin\'{e}aire \textbf{13} (6) (1996),
635--690.
\bibitem{agrs4}
A. A. Agrachev and A. V. Sarychev, \textit{On abnormal extremals
for Lagrange variational problems.} J. Math. Syst. Estim. Cont.
\textbf{8} (1) (1998), 87--118.
\bibitem{agrs5}
A. A. Agrachev and A. V. Sarychev, \textit{Sub-Riemannian metrics: Minimality of abnormal geodesics versus subanalyticity.} ESAIM Control Optim. Calc. Var. \textbf{4} (1999).
\bibitem{AMal}
H.~Airault, P.~Malliavin, \textit{Int\'{e}gration g\'{e}ometric
sur l'espace de Wiener}, Bull. Sci. Math. \textbf{112} (1988), 3--52.
\bibitem{AK} L. Ambrosio, B. Kirchheim, \textit{Rectifiable sets in metric and
Banach spaces.} Math. Ann. \textbf{318} (2000), 527--555.
\bibitem{ascv} L.~Ambrosio, F.~Serra Cassano and D.~Vittone, \textit{Intrinsic regular hypersurfaces in Heisenberg groups.} J.
Geom. Anal. \textbf{16} (2) (2006), 187--232.
\bibitem{Bal}{Z. M. Balogh},
\textit{Size of characteristic sets and functions with prescribed
gradients}, Crelle's Journal \textbf{564} (2003), 63--83.
\bibitem{btw} Z.~M.~Balogh, J.~T.~Tyson, B.~Warhurst, \textit{Sub-Riemannian Vs. Euclidean Dimension Comparison and Fractal Geometry on Carnot Groups}, preprint
\bibitem{B} A. Bella\"{\i}che, \textit{Tangent Space in Sub-Riemannian
Geometry.}   Sub-Riemannian geometry,  Birkh\"auser,
Basel, 1996, 1--78.
\bibitem{birolim1} {M. Biroli, U. Mosco,} \textit{Formed de
Dirichlet et estimationes structurelles dans les mileux discontinues.}
{C. R. Acad.\ Sci. Paris}, \textbf{313} (1991), 593--598.
\bibitem{birolim4} {M. Biroli, U. Mosco,} \textit{Sobolev inequalities
on homogeneous spaces.} {Pot. Anal.} \textbf{4} (1995), 311--324.
\bibitem{blu} A. Bonfiglioli, E. Lanconelli, F. Uguzzoni,
\textit{Stratified Lie Groups and Potential Theory for Their
Sub-Laplacians.} Springer, 2007.
\bibitem{buckleykl1} {S.~M.~Buckley, P.~Koskela, G.~Lu,}
\textit{Subelliptic Poincar\'e inequalities: the case $p<1$.}
{Publ.\ Mat.} \textbf{39} (1995), 313--334.
\bibitem{bbi}
{D.~Yu.~Burago, Yu.~D.~Burago, S.~V.~Ivanov}, \textit{A Course in Metric Geometry.} Graduate Studies in Mathematics, {\bf 33},
American Mathematical Society, Providence, RI, 2001.
\bibitem{cap1} L. Capogna,
\textit{Regularity of quasi-linear equations in the Heisenberg group},
Comm. Pure Appl. Math. \textbf{50} (9) (1997) 867-889.
\bibitem{cap2}
L. Capogna,
\textit{Regularity for quasilinear equations and $1$-quasiconformal
maps in Carnot groups}, Math. Ann. \textbf{313} (2)  (1999), 263-295.
\bibitem{capognadg2} {L.~Capogna, D.~Danielli, N.~Garofalo,}
\textit{An imbedding theorem and the Harnack inequality for nonlinear
subelliptic equations.} {Comm.\ Partial Diff.\ Equations} \textbf{18}
(1993), 1765--1794.
\bibitem{capognadg2.5} {L.~Capogna, D.~Danielli, N.~Garofalo,}
\textit{Subelliptic mollifiers and characterization of Rellich and Poincar\'e
domains.} {Rend.\ Sem.\ Mat.\ Univ.\ Polit.\ Torino} \textbf{54} (1993),
361--386.
\bibitem{capognadg4} {L.~Capogna, D.~Danielli, N.~Garofalo,}
\textit{The geometric Sobolev embedding for vector fields and the
isoperimetric inequality.} {Comm. Anal.\ Geom.\ } \textbf{2} (1994),
203--215.
\bibitem{capognadg5} {L.~Capogna, D.~Danielli, N.~Garofalo,}
\textit{Subelliptic mollifiers and a basic pointwise estimate of
Poincar\'e type.} {Math. Zeit.} \textbf{226} (1997), 147--154.
\bibitem{capognadg6} {L.~Capogna, D.~Danielli,  N.~Garofalo,}
\textit{Capacitary estimates and the local behavior of solutions to nonlinear
subelliptic equations.} Amer. J. Math. \textbf{118} (1996), 1153--1196.
\bibitem{cdpt}
L.~Capogna, D.~Danielli, S.~D.~Pauls and J.~T.~Tyson, \textit{An introduction to the Heisenberg group and the
sub-Riemannian isoperimetric problem.} Progress in Mathematics \textbf{259}. Birkh\"{a}user, 2007.
\bibitem{chernikovv}
{V.~M.~Chernikov, S.~K.~Vodop'yanov,}
\textit{Sobolev Spaces and hypoelliptic equations I,II.}
{Siberian Advances in Mathematics.} \textbf{6} (3) (1996) 27--67;
\textbf{6} (4), 64--96. Translation from:
Trudy In-ta matematiki RAN. Sib. otd-nie. 29 (1995), 7--62.
\bibitem{chow} {W. L. Chow},
\textit{\"Uber Systeme von linearen partiellen Differentialgleichungen
erster Ordung}, Math. Ann. \textbf{117} (1939), 98-105.
\bibitem{cgl}{G.Citti, N.Garofalo, E.Lanconelli},
\textit{Harnack's inequality for sum of squares of
vector fields plus a potential}, Amer. J. Math. {\bf 115} (3) (1993) 699-734.
\bibitem{CS} G. Citti, A. Sarti, \textit{A cortical based model of perceptual
completion in the roto-translation space.} Lecture Notes of
Seminario Interdisciplinare di Matematica \textbf{3} (2004),
145--161.
\bibitem{daniellign}
{D.~Danielli, N.~Garofalo, D.-M.~Nhieu,}
\textit{Trace inequalities for Carnot--Carath\'eodory spaces and
applications
to quasilinear subelliptic equations}, preprint.
\bibitem{dgn}D.~Danielli, N.~Garofalo and D.-M.~Nhieu, \textit{Non-doubling Ahlfors measures, perimeter measures, and
the characterization of the trace spaces of Sobolev functions in Carnot--Carath\'{e}odory spaces.} Mem. Amer. Math. Soc. {182} \textbf{857} (2006).
\bibitem{el1} Ya. Eliashberg, \textit{Contact 3-Manifolds Twenty Years Since J. Martinet's Work.} Ann. Inst Fourier (Grenoble) \textbf{42} (1992), 1-12.
\bibitem{el2} Ya. Eliashberg, \textit{New Invariants of Open Symplectic and Contact Manifolds.} J. Amer. Math. Soc. \textbf{4} (1991), 513-520.
\bibitem{el3} Ya. Eliashberg, \textit{Classification of overtwisted contact structures on 3-manifolds.} Invent. Math. \textbf{98} (1989), 623-637.
\bibitem{el4} Ya. Eliashberg, \textit{Legendrian and transversal knots in tight contact -manifolds.} Topological Methods in Modern Mathematics (1993), 171-193.
\bibitem{eg} L. C. Evans, R. F. Gariepy, \textit{Measure theory and fine
properties of functions.} CRC Press, Boca Raton, 1992.
\bibitem{fd1} H. Federer, \textit{Curvature measures.} Trans. Amer. Math. Soc. \textbf{93} (1959),
418--491.
\bibitem{F} H. Federer, \textit{Geometric Measure Theory.} NY: Springer, 1969.
\bibitem{FF} H. Federer, W. H. Fleming \textit{Normal and Integral Currents.} Ann. Math. \textbf{72} (2) (1960),
458--520.
\bibitem{feffermanp} C.~Fefferman, D.~H.~Phong, \textit{Subelliptic
eigenvalue problems.} {Proceedings of the conference in
harmonic analysis in honor of Antoni Zygmund}, Wadsworth Math.\ Ser.,
Wadsworth, Belmont, California, 1981, 590--606.
\bibitem{Foll}{G. B. Folland},
\textit{A fundamental solution for a subelliptic operator},
Bull. Amer. Math. Soc. {\bf 79} (1973), 373-376.
\bibitem{Foll1}{G. B. Folland},
\textit{Subelliptic estimates and function spaces on nilpotent Lie groups},
Ark. Mat. {\bf 13} (2) (1975) 161--207.
\bibitem{fs} G.~B. Folland, E.~M. Stein, \textit{Hardy spaces on
homogeneous groups.}  Princeton Univ. Press, 1982.
\bibitem{franchi1} {B.~Franchi,} \textit{Weighted Sobolev--Poincar\'e
inequalities and pointwise inequalities for a class of degenerate
elliptic equations.} {Trans.\ Amer.\ Math.\ Soc.\ } \textbf{327}
(1991), 125--158.
\bibitem{franchigaw} {B.~Franchi, S.~Gallot,  R.~Wheeden,} \textit{Sobolev and isoperimetric inequalities for degenerate
metrics.} {Math.\ Ann.\ } \textbf{300} (1994), 557--571.
\bibitem{franchigw} {B.~Franchi, C.~E.~Guti\'{e}rrez,
R.~L.~Wheeden,} \textit{Weighted Sobolev--Poincar\'{e} inequalities for Grushin
type operators.} {Comm.\ Partial Differential Equations} \textbf{19} (1994),
523--604.
\bibitem{franchil1} {B.~Franchi, E.~Lanconelli,} \textit{H\"older
regularity theorem for a class of non uniformly elliptic operators
with measurable coefficients.} {Ann.\ Scuola Norm.\ Sup.\ Pisa}
\textbf{10} (1983), 523--541.
\bibitem{franchil2} {B.~Franchi, E.~Lanconelli,} \textit{An
imbedding theorem for Sobolev spaces related to non smooth vector
fields and Harnack inequality.} {Comm.\ Partial Differental
Equations} \textbf{9} (1984), 1237--1264.
\bibitem{franchilw1} {B.~Franchi, G.~Lu, R.~Wheeden,}
\textit{Representation formulas and weighted Poincar\'e inequalities for
H\"ormander vector fields.}
{Ann.\ Inst.\ Fourier (Grenoble)} \textbf{45} (1995), 577-604.
\bibitem{franchilw3}
{B.~Franchi, G.~Lu, R.~Wheeden,}
\textit{A relationship between Poincar\'e type inequalities and
representation formulas in spaces of homogeneous type.}
{Int. Mat. Res. Notices} (1) (1996), 1--14.
\bibitem{franchis1} {B.~Franchi, R.~Serapioni,} \textit{Pointwise
estimates for a class of strongly degenerate elliptic operators: a
geometric approach.} {Ann.\ Scuola Norm.\ Sup.\ Pisa} \textbf{14}
(1987), 527--568.
\bibitem{fssc1}
B.~Franchi, R.~Serapioni and F.~Serra Cassano, \textit{Regular hypersurfaces, intrinsic perimeter and implicit
function theorem in Carnot groups.} Comm. Anal. Geom. \textbf{11} (5) (2003), 909--944.
\bibitem{fssc2} B.~Franchi, R.~Serapioni and F.~Serra Cassano, \textit{Regular submanifolds, graphs and area formula in
Heisenberg groups}, Adv. Math. \textbf{211} (1) (2007), 152--203.
\bibitem{fssc3}
B. Franchi, R. Serapioni, F. Serra Cassano, \textit{Rectifiability
and Perimeter in the Heisenberg group}, Math. Ann. \textbf{321}
(3) (2001), 479--531.
\bibitem{fssc4}
B. Franchi, R. Serapioni, F. Serra Cassano, \textit{On the
structure of finite perimeter sets in step 2 Carnot groups},
J.~Geom. Anal.~\textbf{13} (3) (2003), 421--466.
\bibitem{garattini} {R.~Garattini},
\textit{Harnack's inequality on homogeneous spaces}, Annali di
Matematica Pura ed Applicata \textbf{179} (1) (2001), 1--16.
\bibitem{garofalobook} N. Garofalo, \textit{Analysis and Geometry of Carnot--Carath\'{e}odory Spaces, With Applications to PDE's}, Birkh\"{a}user, in preparation.
\bibitem{garofalol} {N.~Garofalo, E.~Lanconelli,}
\textit{Existence and nonexistence results for
semilinear equations on the Heisenberg group.} {Indiana Univ.\
Math.\ J.} \textbf{41} (1992), 71--98.
\bibitem{GN1}{N. Garofalo, D.-M. Nhieu},
\textit{Isoperimetric and Sobolev Inequalities for Carnot-Carath\'eodory
Spaces and the Existence of Minimal Surfaces}, Comm. Pure Appl. Math.
\textbf{49} (1996), 1081-1144.
\bibitem{GN2}{N. Garofalo, D.-M. Nhieu},
\textit{Lipschitz continuity, global smooth approximation
and extension theorems for Sobolev functions in
Carnot-Carath\'eodory spaces}, Jour. Anal. Math., \textbf{74} (1998), 67-97.
\bibitem{gq} M. Giaquinta, G. Modica, J. Sou\v cek, \textit{Cartesian currents in
the calculus of variations. V.~I,~II.}  Springer-Verlag, Berlin,
1998.
\bibitem{good} R.~W.~Goodman, \textit{Nilpotent Lie groups{\rm:} structure and applications to
 analysis.}  Springer-Verlag, Berlin-Heidelberg-New York,
1976. Lecture Notes in Mathematics, vol. 562.
 \bibitem{greshn} A. V. Greshnov, \textit{Metrics of Uniformly Regular
Carnot--Carath\'{e}odory Spaces and Their Tangent Cones.} Sib. Math. Zh.
\textbf{47} (2) (2006), 259--292.
\bibitem{greshn1} A. V. Greshnov, \textit{Local Approximation
of Equiregular Carnot--Carath'eodory Spaces by its Tangent Cones.} Sib. Math. Zh.
\textbf{48} (2) (2007), 290--312.
\bibitem{G} M. Gromov, \textit{Carnot--Carath\'{e}odory Spaces Seen From
Within.}   Sub-Riemannian geometry, Birkh\"auser, Basel, 1996, 79--318.
\bibitem{G2001}
M. Gromov, \textit{Metric Structures for Riemannian and Non-Riemannian Spaces.} Birkh\"{a}user, 2001.
\bibitem{hk} P. Haj{\l}asz, P.~Koskela, \textit{Sobolev Met Poincar\'{e}.} Memoirs of the American Mathematical Society
\textbf{145} (2000), no.~688.
\bibitem{hajlaszs} {P.~Haj\l{}asz, P.~Strzelecki,}
\textit{Subelliptic $p$-harmonic maps into spheres and the ghost
of Hardy spaces}, {Math.\ Ann.} \textbf{312} (1998), 341--362.
\bibitem{H} J. Heinonen, \textit{Calculus on Carnot groups.}  {Fall school in
analysis,} (Jyv\"askyl\"a, 1994). Jyv\"askyl\"a, University
of Jyv\"askyl\"a, 1994, pp. 1--32.
\bibitem{HP} R. K. Hladky, S. D. Pauls, \textit{Minimal surfaces in the
roto-translation group with applications to a neuro-biological
image completion model.} arXiv:math.DG/0509636, 27 Sep. 2005.
\bibitem{hor} L.~H\"ormander, \textit{Hypoelliptic second order differential equations.}
 Acta Math.
 \textbf{119}  (1967),  147--171.
 \bibitem{J} F. Jean, \textit{Uniform estimation of sub-riemannian balls.} Journal on Dynamical and Controle Systems
\textbf{7} (4) (2001), 473--500.
\bibitem{jerison} {D.~Jerison,}
\textit{The Poincar\'{e} inequality for
vector fields satisfying H\"{o}rmander's condition.}
{Duke Math. J.}   \textbf{53} (1986), 503--523.
\bibitem{jost1} {J. Jost,} \textit{Equilibrium maps between metric spaces.}
{Calc.\ Var.\ } \textbf{2} (1994), 173--205.
\bibitem{jost2} {J. Jost,} \textit{Generalized harmonic maps between
metric spaces}, in: {Geometric analysis and the calculus of
variations} (J. Jost, ed.), International Press, 1966, 143--174.
\bibitem{jost3} {J. Jost,} \textit{Generalized Dirichlet forms and
harmonic maps.} {Calc.\ Var.} \textbf{5} (1997), 1--19.
\bibitem{jost4} {J. Jost,} \textit{Nonlinear Dirichlet forms},
preprint.
\bibitem{jostx} {J. Jost, C.~J.~Xu,} \textit{Subelliptic harmonic maps.}
{Trans.\ Amer.\ Math.\ Soc.} \textbf{350} (1998), 4633-4649.
\bibitem{jur} V.~Jurdjevic, \textit{Geometric Control Theory}, Cambridge Studies in Mathematics \textbf{52}. Cambridge University Press, 1997.
\bibitem{Km1} M. B. Karmanova, \textit{Metric Differentiability of Mappings and
Geometric Measure Theory.} Doklady Mathematics \textbf{71} (2)
(2005), 224--227.
\bibitem{Km2} M. B. Karmanova, \textit{Rectifiable Sets and the Coarea Formula for Metric-Valued Mappings.} Doklady Mathematics \textbf{73} (3)
(2005), 323--327.
\bibitem{Km3} M. Karmanova, \textit{Geometric Measure Theory Formulas on
Rectifiable Metric Spaces.} Contemporary Mathematics \textbf{424}
(2007), 103--136.
\bibitem{Km4} M. Karmanova, \textit{Rectifiable Sets and Coarea Formula for Metric-Valued Mappings},
Journal of Functional Analysis \textbf{254} (5) (2008),
1410--1447.
\bibitem{ksc} B.~Kirchheim and F.~Serra Cassano, \textit{Rectifiability and parameterization of intrinsic regular surfaces in the
Heisenberg group.} Ann. Sc. Norm. Super. Pisa Cl. Sci. (5) \textbf{3} (4) (2004), 871--896.
\bibitem{kr} A. S. Kronrod, \textit{On functions of two variables.} Uspekhi
Matematicheskikh Nauk (N.~S.) \textbf{5} (1950), 24--134.
\bibitem{LeoRig}{G. P. Leonardi, S. Rigot},
\textit{Isoperimetric sets on Carnot groups}, Houston Jour. Math.
\textbf{29} (3) (2003), 609-637.
\bibitem{fx} F.~Lin, X.~Yang, \textit{Geometric measure theory
--- an introduction.} Science Press, Beijing~a.~o., 2002.
\bibitem{LS}{W. Liu, H. J. Sussman},
\textit{Shortest paths for sub-Riemannian metrics on rank-two  distributions},  Mem. Amer. Math. Soc. \textbf{118} (564) (1995).
\bibitem{Lu}{G. Lu},
\textit{Weighted Poincar\'e and Sobolev inequalities for vector fields
satisfying H\"{o}rmander's condition and applications},
Rev. Mat. Iberoamericana {\bf 8} (3) (1992), 367-439.
\bibitem{M1} V. Magnani, \textit{The coarea formula for real-valued Lipschitz maps on stratified groups.} Math. Nachr. \textbf{27} (2)
(2001), 297--323.
\bibitem{mblow} V.~Magnani, \textit{A blow-up theorem for regular hypersurfaces on nilpotent groups.} Man. math. \textbf{110} (2003), 55--76.
\bibitem{M00} V. Magnani, \textit{Elements of Geometric Measure Theory on
sub-Riemannian groups.} Tesi di Perfezionamento. Pisa: Scuola
Normale Superiore (Thesis), 2002.
\bibitem{M05} V. Magnani, \textit{Blow-up of regular submanifolds in Heisenberg
groups and applications.} Cent. Eur. J. Math. \textbf{4} (1)
(2006), 82--109.
\bibitem{Mal}
P.~Malliavin, \textit{Stochastic Analysis}, Springer, NY, 1997.
\bibitem{marchi}
{S. Marchi,}
\textit{H\"older continuity and Harnack inequality for De Giorgi classes
related to H\"ormander vector fields.}
{Ann.\ Mat.\ Pura Appl.} \textbf{168} (1995), 171--188.
\bibitem{MM} G. A. Margulis, G. D. Mostow, \textit{The differential of
quasi-conformal mapping of a Carnot--Carath\'e\-odory spaces.}
Geometric and Functional Analysis \textbf{5} (2) (1995), 402--433.
\bibitem{MM1} G. A. Margulis, G. D. Mostow, \textit{Some remarks on the
definition of tangent cones in a Carnot--Carath\'eodory space.}
Journal D\'Analyse Math. \textbf{80} (2000), 299--317.
\bibitem{Me} G.~Metivier, \textit{Fonction spectrale et valeurs propres
      d'une classe d'op\'{e}rateurs  non elliptiques.}
Commun. Partial Differential Equations \textbf{1} (1976), 467--519.
\bibitem{Mit} J. Mitchell, \textit{On Carnot--Carath\'eodory metrics.}
J. Differential Geometry \textbf{21} (1985), 35--45.
\bibitem{Montgom1}{R. Montgomery},
\textit{Abnormal minimizers}, SIAM J. Control Optim. \textbf{32} (6) (1994), 1605-1620.
\bibitem{M} R. Montgomery, \textit{A Tour of Subriemannian Geometries, Their
Geodesics and Applications.} Providence, AMS, 2002.
\bibitem{Mon}{R. Monti},
\textit{Some properties of Carnot-Carath\'eodory balls in the Heisenberg group},
Atti Accad. Naz. Lincei Cl. Sci. Fis. Mat. Natur. Mem. s.9 \textbf{11} (2000), 155-167.
\bibitem{Mon1}{R. Monti},
\textit{Distances, boundaries and surface measures in Carnot-Carath\'eodory
spaces}, PhD thesis
\bibitem{MS}
{R. Monti and F. Serra Cassano}, \textit{Surface measures in
Carnot-Carath\'edory spaces}, Calc. Var. Partial Differential
Equations \textbf{13} (3) (2001), 339--376.
\bibitem{NRS}{A. Nagel, F. Ricci, E. M. Stein},
\textit{Fundamental solutions and harmonic analysis on nilpotent groups},
Bull. Amer. Math. Soc. (N.S.) \textbf{23} (1) (1990), 139-144.
\bibitem{NRS1}{A. Nagel, F. Ricci, E. M. Stein},
\textit{Harmonic analysis and fundamental solutions on nilpotent Lie groups},
Analysis and partial differential equations, 249-275,
Lecture Notes in Pure and Appl. Math. {\bf 122},
Dekker, New York, 1990.
\bibitem{nsw} A. Nagel, E. M. Stein, S. Wainger, \textit{Balls and
metrics defined by vector fields I: Basic properties.} Acta Math.
\textbf{155} (1985), 103--147.
\bibitem{ot} M. Ohtsuka, \textit{Area Formula.} Bull. Inst. Math. Acad. Sinica
\textbf{6} (2) (2) (1978), 599--636.
\bibitem{P} P. Pansu, \textit{Geometrie du group d'Heisenberg.} Univ. Paris VII,
1982.
\bibitem{Pan01}{P. Pansu},
\textit{Une in\'egalit\'e isoperimetrique sur le groupe de Heisenberg},
C.R. Acad. Sc. Paris, \textbf{295} S\'erie I (1982), 127-130.
\bibitem{Pan1}{P. Pansu},
\textit{Croissance des boules et des g\'eod\'esiques ferme\'es dans les
nilvari\'et\'e}, Ergod. Dinam. Syst. {\bf 3} (1983) 415-445.
\bibitem{Pan} P. Pansu, \textit{M\'{e}triques de Carnot--Carath\'{e}odory et quasiisom\'{e}tries
des espaces sym\'{e}triques de rang un.} Ann. Math. (2) \textbf{129} (1) (1989), 1--60.
\bibitem{pauls} S.~D.~Pauls, \textit{A Notion of Rectifiability Modeled on Carnot Groups},
Indiana University Mathematics Journal \textbf{53} (2004) 49-82.
\bibitem{post} M. M. Postnikov, \textit{Lectures in Geometry. Semester V:
Lie Groups and Lie Algebras.} Moscow, ''Nauka``, 1982.
\bibitem{rash} {P. K. Rashevsky},
\textit{Any two point of a totally nonholonomic space may be connected
by an admissible line}, Uch. Zap. Ped. Inst. im. Liebknechta. Ser.
Phys. Math. \textbf{2} (1938), 83-94.
\bibitem{rs} L. P. Rothschild, E. M. Stein, \textit{Hypoelliptic differential
operators and nilpotent groups.} Acta Math. \textbf{137} (1976),
247--320.
\bibitem{sanchez}
{A. S\'anchez-Calle,}
\textit{Fundamental solutions and geometry of sums of squares of vector
fields.}
{Invent. Math.} \textbf{78} (1984), 143--160.
\bibitem{Stein} E.~M.~Stein, \textit{Harmonic analysis{\rm:} real-variables methods, orthogonality, and oscillatory integrals.} Princeton, NJ, Princeton University Press, 1993.
\bibitem{Stri}{R. S. Strichartz},
\textit{Sub-Riemannian geometry}, J. Diff. Geom. \textbf{24} (1986) 221-263.
Corrections: J. Diff. Geom. {\bf 30} (1989), 595-596.
\bibitem{sturm4} {K. T. Sturm,} \textit{Analysis on local Dirichlet
spaces III. The parabolic Harnack inequality.} {J. Math.\ Pures
Appl.\ } \textbf{75} (1996), 273--297.
\bibitem{VerG} A.~M.~Vershik, V.~Ya.~Gershkovich, \textit{Nonholonomic dynamical systems, geometry of distributions and variational problems.} Dynamical systems. VII. Encycl. Math. Sci.,
vol.~16.  1994, pp.~1--81 (English Translation from: Itogi Nauki Tekh., Ser. Sovrem. Probl. Mat., Fundam. Napravleniya vol.~16, Moscow, VINITI 1987, pp. 7--85).
\bibitem{V4} S.~K.~Vodop$'$yanov, \textit{$\mathcal P$-differentiability on Carnot groups in different topologies and related topics.}
Proceedings on Analysis and Geometry (S. K. Vodopyanov, ed.)
 Novosibirsk,  Sobolev Institute Press, 2000, pp. 603--670.
 \bibitem{V1} S. K. Vodopyanov, \textit{Theory of Lebesgue Integral: Lecture Notes
on Analysis.} Novosibirsk: NSU Publishing, 2003.
\bibitem{vod1} S. K. Vodopyanov, \textit{Differentiability of Curves in
Carnot Manifold Category.} Dokl.~AN. \textbf{410} (4) (2006),
1--6.
\bibitem{vod2} S. K. Vodopyanov, \textit{Differentiability of mappings of Carnot
Manifolds and Isomorphism of Tangent Cones.} Doklady Mathematics
\textbf{74} (3) (2006), 844--848.
\bibitem{vod3} S. K. Vodopyanov, \textit{Geometry of Carnot--Carath\'{e}odory
Spaces and Differentiability of Mappings.} Contemporary
Mathematics \textbf{424} (2007), 247--302.
\bibitem{V3} S. K. Vodopyanov, \textit{Differentiability of Mappings in Carnot
Manifold Geometry.} Sib. Mat. J. \textbf{48} (2) (2007), 251--271.
\bibitem{vg} S. K. Vodopyanov, A. V. Greshnov, \textit{On the
Differentiability of Mappings of Carnot--Carath\'{e}odory Spaces.}
Doklady Mathematics \textbf{67} (2) (2003), 246--250.
\bibitem{vk1} S. K. Vodopyanov, M. B. Karmanova, \textit{Local Geometry
of Carnot Manifolds Under Minimal Smoothness.} Doklady Mathematics
\textbf{75} (2) (2007), 240--246.
\bibitem{vk2} S. K. Vodopyanov, M. B. Karmanova, \textit{Sub-Riemannian geometry
under minimal smoothness of vector fields.} Doklady Mathematics
\textbf{}  (2008), to appear.
\bibitem{vk3} S. K. Vodopyanov, M. B. Karmanova, \textit{Coarea Formula
for Smooth Contact Mappings of Carnot Manifolds.} Doklady Mathematics
\textbf{76} (3) (2007), 908--912.
\bibitem{VU4} S. K. Vodopyanov, A. D. Ukhlov, \textit{Approximately differentiable transformations and change of variables on nilpotent groups.}
Siberian Math. J. \textbf{37} (1) (1996), 62--78.
\bibitem{VU1} S. K. Vodopyanov, A. D. Ukhlov, \textit{Set functions and their
applications in the theory of Lebesgue and Sobolev spaces. I.}
Siberian Adv. Math. \textbf{14} (4) (2004), 78--125. \textit{II.}
Siberian Adv. Math. \textbf{15} (1) (2005), 91--125.
\bibitem{W}F.~W.~Warner, \textit{Foundations of differentiable manifolds and Lie groups.} New York a.~o., Springer-Verlag, 1983. Graduate Texts in Mathematics, vol. 94.
\bibitem{XZ}{C. J. Xu, C. Zuily},
\textit{Higher interior regularity for quasilinear subelliptic sy\-stems},
Calc. Var. Partial Differential Equations \textbf{5} (4) (1997) 323-343.
\end{thebibliography}
\end{document}